\documentclass[11pt,twoside]{article}
\usepackage{amsfonts,amssymb,euscript,a4,epsfig}
\input{macro}
\newcommand{\bfigure}{\begin{figure}}
\newcommand{\efigure}{\end{figure}}

\newcommand{\ff}{\emph{focus-focus}}
\newcommand{\D}{\mathcal{D}}
\newcommand{\C}{\mathcal{C}}
\newcommand{\macaption}[1]{\caption{\textsf{#1}}}
\newcommand{\brho}{\mbold{\rho}}
\newcommand{\out}{\textrm{\tiny out}}
\newcommand{\SY}{\mathcal{S}}

\title{\textsf{\BS\  conditions for Integrable Systems\\ with
    critical manifolds of \ff\ type}}

\author{V\~u Ng\d oc San}
\date{September 17, 1998}

\begin{document}

\maketitle

{\small
  \begin{center}
 
Mathematics Institute,\\
Budapestlaan 6, University of Utrecht, 3508 TA Utrecht, The
Netherlands.\\ \&\\ 
   Institut Fourier\\
Unit\'e mixte de recherche CNRS-UJF 5582\\
B.P. 74, 38402 Saint-Martin d'H\`eres, France\\
e-mail : \texttt{vu-ngoc@math.uu.nl}

  \end{center}
}

\vfill
\begin{abstract}
We present a detailed study, in the semi-classical regime $h \fleche
0$, of microlocal properties of systems of two commuting $h$-\pdo s
$P_1(h)$, $P_2(h)$ such that the joint principal symbol $p=(p_1,p_2)$
has a special kind of singularity called a \ff\ singularity. Typical
examples include the quantum spherical pendulum or the quantum
Champagne bottle.

In the spirit of Colin de Verdi\`ere and Parisse
\cite{colin-p,colin-p2,colin-p3}, we show that such systems have a
universal behavior described by \emph{singular quantization
conditions} of \BS\ type.

These conditions are used to give a precise description of the
\emph{joint spectrum} of such systems, including the phenomenon of
quantum monodromy and different formulations of the counting function
for the joint eigenvalues close to the singularity, in which a
logarithm of the semi-classical constant $h$ appears. Thanks to
numerical computations done by M.S. Child for the case of the
Champagne bottle, we are able to accurately illustrate our statements.
\end{abstract}

\vfill

\newpage
\setlength{\parskip}{1ex plus 0.5ex minus 0.2ex}
\section{Introduction}
In the long history of \cis s, an important object was discovered
quite recently (Duistermaat \cite{duistermaat})~: the monodromy of the
system, whose non-triviality prevents the construction of global
action variables.  The question about the impact of this invariant on
the spectrum of quantum integrable system was raised by Cushman and
Duistermaat in \cite{duist-cushman}; an answer is proposed in
\cite{san-mono}.  The issue is to describe the joint spectrum of two
commuting $h$-\pdo s $P_1$, $P_2$ in a region close to a critical
point of the $\Cinf$ joint principal symbol
\[ p\egdef (p_1,p_2), \] 
when the underlying Liouville integrable system 
${p_1,p_2}$ has non trivial monodromy.

However, because of a well-known drawback of both the usual $WKB$ 
construction and the standard \BS\ quantization conditions, the 
descriptions had to keep a reasonable distance away from the singular 
value of $p$.  Here reasonable may be small but means fixed, as 
$h$ tends to zero.  In this way, an increasing number of eigenvalues 
(as $h\fleche 0$) remained out of control.

On the other hand, recent achievements in semi-classical analysis of
\schr\ operators near a critical point of the potential, often via the
use of microlocal normal forms, suggested that this problem should be
solvable. I am referring here for instance to the work of
Helffer-Robert \cite{helffer-robert}, Helffer-Sj\"ostrand
\cite{helffer-sjostrand}, Sj\"ostrand \cite{sjostrand2}, M\"arz
\cite{marz}, Brummelhuis-Paul-Uribe \cite{brummelhuis}, \dots and in
particular to the articles by Colin de Verdi\`ere and Parisse
\cite{colin-p,colin-p2,colin-p3}, in which the case of a local maximum
of a smooth ($\Cinf$) potential for a 1-dimensional \schr\ operator is
treated. Their method rested upon a smooth normal form theorem and on
the study of the hyperbolicity of the classical Hamiltonian flow.

It turns out that any two degree of freedom quantum integrable system
with a \ff\ singularity exhibits at the same time a non-trivial
monodromy and a hyperbolic behavior of the Hamiltonian flow. Moreover,
such a singularity admits a smooth normal form, due to Eliasson
\cite{eliasson-these}, that has a semi-classical analogue
\cite{san}. This allows us to settle an analysis combining a
geometrical description of the underlying classical \cis\  with a
microlocal analysis near the singularity.

The main result of this paper is the statement of the \emph{singular}
\BS\  quantization conditions (theorem \ref{theo:global}),
which are uniform in a neighborhood of the critical value. In other
words, these conditions are able to uniformly describe an increasing
number of joint eigenvalues, and as such they contain the description
of the quantum monodromy far from the singularity as well as the
asymptotic distribution of the joint eigenvalues near the singularity.

The organization of this paper is as follows. In the first parts, the
aim is to settle the ingredients needed to have a global picture of the
problem. The necessary background concerning Liouville integrable
systems is recalled (section \ref{sec:classique}), and the general
notion of \emph{semi-classical} integrable systems is explained in
section \ref{sec:semi-classique}, where we point out the role of both
the principal and sub-principal symbols. The version of the microlocal
tools that we will be using is given in section \ref{sec:analyse}; its
application to \cis s is developed in the next
section (\ref{sec:CBSreg}). Besides known results concerning the WKB
construction, the most fundamental results explaining our methods are
probably propositions \ref{prop:dim1reg} and \ref{prop:holonomy}.

The rest of the paper is devoted to the case of a \ff\ singularity in
a 4-dimensional cotangent bundle. Section \ref{sec:focus-focus}
contains the microlocal analysis and the geometrical description of
the monodromy that lead to our main result (theorem
\ref{theo:global}).  Finally, we derive in section \ref{sec:spectrum}
the structure of the joint spectrum. Theorem \ref{theo:compter} shows
that one can count the number of joint eigenvalues without knowing
their precise distribution near the critical value, and theorem
\ref{theo:gaps}, on the contrary, studies this distribution, by
estimating the gaps between the joint eigenvalues.  These statements
are illustrated by the example of the Champagne bottle.

\noindent\textbf{Acknowledgments.} 
First of all, I would like to thank Y. Colin de Verdi\`ere whose
numerous ideas and suggestions made this work possible. I also had
stimulating discussions with J.J.Duistermaat, particularly concerning
the monodromy. I would like to thank him for this and for his
hospitality at Utrecht. Finally, it was a pleasure to talk with
M.S.Child, whose viewpoint as a chemist was quite interesting for me
(see \cite{child}), and who accepted to provide me with a load of
numerics concerning the quantum ``Champagne bottle''. Let him be
thanked for this.

\noindent My research is supported by a Marie Curie Fellowship
Nr. ERBFMBICT961572.

\newpage
\setlength{\parskip}{0ex}
{\small \tableofcontents}
\setlength{\parskip}{1ex plus 0.5ex minus 0.2ex}

\section{Classical \cis s}
\label{sec:classique}
\subsection{Definition}
Let $(M,\omega)$ be a symplectic manifold of dimension $2n$. A
completely (or ``Liouville'') integrable system on $M$ is the data of
$n$ functions $f_1,\dots,f_n$ in involution with respect to the
symplectic Poisson bracket, with the requirement that their
differentials $df_i$ are almost everywhere independent. The function
\[F~: M\ni m \mapsto (f_1(m),\dots,f_n(m))\in \RM^n\]
is called the momentum map. It is indeed a momentum map for the Abelian
infinitesimal action of $\RM^n$ into the Lie algebra $\ham{}(M)$ of
Hamiltonian vector fields on $M$, given by the generators
$(\ham{1},\dots,\ham{n})$. Here we have denoted by $\ham{i}$ the
Hamiltonian vector field $\ham{f_i}$ associated to $f_i$. The flows of
these vector fields yield a local Abelian action of $\RM^n$ on $M$,
simply referred to in the sequel by ``the action'' or ``the flow'' of
the system.

Note that each $f_i$, and hence every function of the form
$g(f_1,\dots,f_n)$ is locally constant under this action. The momentum
map $F$ thus defines a singular fibration by invariant leaves. These
leaves $\Lambda_c\egdef F^{-1}(c)$, $c\in\RM^n$, are generically
Lagrangean submanifolds of $M$. We shall always assume that $F$ is
\emph{proper}, so that each $\Lambda_c$ is compact.

\subsection{Known results} 
Local and semi-global descriptions of \cis s are provided by the
following facts~:
\begin{itemize}
\item when the vector fields $(\ham{1},\dots,\ham{n})$ are independent
  at a point $m$, the \emph{Darboux-Ca\-ra\-th{\'e}o\-do\-ry theorem}
  states that the functions $f_i$ can be taken as ``momentum''
  coordinates of a symplectic chart $\{(x,\xi)\}$ on a neighborhood
  $U$ of $m$. In these coordinates, the foliation $\Lambda_c\cap U$
  for $c$ close to $F(m)$ is given by $\xi=c$. The $x_i$
  variables are therefore local coordinates for $\Lambda_c$.
\item when $c$ is a regular value of $F$ -- this includes the previous
  case, of course -- the \emph{Arnold-Liouville} theorem endows an
  invariant neighborhood $\Omega$ of $\Lambda_c$ with action-angle
  coordinates~: $\Omega$ is symplectomorphic to a neighborhood of the
  zero section of the cotangent bundle $T^*(\T^n)$ of the $n$-torus
  $\T^n=\RM^n/\ZM^n$, and the fibration $F$ gets transformed
  into a smooth function of the momentum variables $\xi_i$ only. As a
  consequence, $\Lambda_c$, as well as the neighboring fibers, are
  Lagrangean tori equipped with an affine structure~: the one given by
  the $\RM^n$ action which, in these coordinates, is linear.
\end{itemize}

In the first case, the $\RM^n$ action is locally free at $m$, and
gives rise to a local diffeomorphism between $(\RM^n,0)$ and
$(\Lambda_c,m)$. By transporting the standard Lebesgue measure via
this diffeomorphism we obtain an invariant measure on a neighborhood
of $m$ in $\Lambda_c$ which in Darboux-Ca\-ra\-th{\'e}o\-do\-ry
coordinates $(x,\xi)$ is nothing else than $|dx|$.

In the second case, this can be done globally to give a natural
invariant measure $\rho_c$ on the whole torus $\Lambda_c$. Letting
$m_c$ be the total mass of $\Lambda_c$, we can
identify $m_c^{-1}\rho_c$ with the Lebesgue (or
Haar) measure on the torus $\T^n$ via action-angle
coordinates. $\rho_c$ smoothly depends on $c$ and is called the
\emph{Liouville} measure.

\begin{rema}
 This terminology sometimes also applies to the symplectic measure
$\omega^{\wedge n}/n!$. Note that the push-forward of this measure by
$F$ gives a measure $\mu_c$ on $\RM^n$ that satisfies
\[ \frac{\omega^{\wedge n}}{n!} = \rho_c\otimes\mu_c, \]
and that is long known to be of particular importance when it comes to
counting eigenvalues (see section \ref{sec:spectrum}). Because of the
article \cite{duist-heckman} $\mu_c$ is sometimes called
Duistermaat-Heckman's measure.
\end{rema}

Since we will be dealing only with pseudo-differential quantization,
we will always assume that $M$ is an open subset of a cotangent bundle
$T^*X$, which implies that the symplectic form is exact~:
$\omega=d\alpha_0$, where $\alpha_0$ is the standard Liouville
1-form. If we let $i_c : \Lambda_c\hookrightarrow M$ be the inclusion,
the fact that $\Lambda_c$ is Lagrangean implies that the 1-form
$i_c^*\alpha_0$ is closed and thus gives rise to an element of the
cohomology $H^1(\Lambda_c,\RM)$. With slight abuse of notation, we
will still call it $[\alpha_0]$.

\section{Semi-classical integrable systems}
\label{sec:semi-classique}
The aim of this section is to define what we consider to be the proper
semi-classical quantization of \cis s (in accordance with
\cite{colinII,charbonnel}), and to present a geometric interpretation
of their sub-principal terms. 

\subsection{Definition}

Let $X$ be a differentiable manifold of dimension $n$, equipped with a
half-density $|dx|^{1/2}$, and let $\Omega$ be an open subset of
$T^*X$. The spaces of \pdo\  that are used here are defined in section
\ref{sec:analyse}; however, we shall here only use the fact that they
are defined up to the order $O(h^2)$ by two functions of $\Omega$~:
their principal and sub-principal symbols.

Throughout this work, \pdo s are always \emph{classical} and formally 
self-adjoint, so that their principal and sub-principal symbols are 
\emph{real}-valued.

A set $\{P_1(h),\dots,P_n(h)\}$ of $n$ \pdo s in $\Psi^0(\Omega)$ of
order zero is called a \emph{semi-classical integrable system} if~:
\begin{itemize}
\item the principal symbols $p_1,\dots,p_n$ form a \cis.
\item $\forall i,j \quad [P_i(h),P_j(h)]=0.$
\end{itemize}
Note that the second condition already implies that the principal
symbols are in involutions.

Because of this definition any semi-classical integrable system has an
underlying classical \cis, which means that the main geometric
ingredient of such a system is the Lagrangean fibration $\Lambda_c$
given by the principal symbols. Of particular importance for us will
be the \emph{principal Lagrangean} $\Lambda_0$.

\begin{rema}
  \label{rem:commute}
  All the results of sections 2-6 would still hold if the commutation 
  property is weakened to $[P_i(h),P_j(h)]=\oh$. Section 7 however 
  requires the exact commutation of the operators, in order to define 
  their joint spectrum. Similarly, a number of results before that 
  section would still be valid if the self-adjointness requirement is 
  dropped and replaced by the assumption that both the principal and 
  sub-principal symbols are real-valued. 
\end{rema}
%In section \ref{sec:analyse}, we develop the necessary notions of
%microlocal analysis we use here. Our viewpoint will always be that
%such a system defines a space of \emph{\mi\ solutions}
%(\ref{sec:solutions}) $u_h$ of the
%equations
%\[ P_ju_h \sim 0, \quad i=1,\dots,n. \] 
%Now, since the \emph{wave front set} of such a solution is included in
%$\Lambda_0$, which is compact, we can always assume that the
%$P_{j,h}$'s are compactly supported. This essentially means that we
%won't need any assumption about their behavior at infinity.

%The aim of this work is to present some quantization conditions as
%necessary and sufficient conditions to the existence of \mi\ solutions
%of the system, when $\Lambda_0$ carries a singular point of a certain
%type (see section \ref{sec:focus-focus}). Theses conditions, as well as
%the celebrated \BS\  conditions for smooth tori (section
%\ref{sec:CBSreg}), are expressed at the principal level through three
%geometric ingredients~: the Liouville 1-form $\alpha_0$, the Maslov
%cocycle $\mu$, and the sub-principal form $\kappa$, to be defined in
%the following paragraphs.

\subsection{The sub-principal form}
\subsubsection{Definition}
\label{sec:sous-principale}
Apart from the principal symbols, the data of a semi-classical
integrable system defines another set of functions, namely the
sub-principal symbols $r_1,\dots,r_n$. We assume here that these
functions are real-valued. While the principal symbols gave rise to
the momentum map and its Lagrangean fibration, the sub-principal
symbols will be viewed as characteristic of an infinitesimal
deformation of the Lagrangean leaves.
\begin{defi}
  The \emph{sub-principal form} $\kappa_c$ of a semi-classical
  integrable system is the differentiable 1-form on $\Lambda_c$
  defined at non-singular points of $F$ by~:
  \[ \kappa_c(\ham{i})=-r_i. \]
\end{defi}
Recall that $\ham{i}$ denote the Hamiltonian vector field of $p_i$. At
such a non-singular point, $(\ham{1},\dots,\ham{n})$ is a basis of the
tangent space of $\Lambda_c$. 

This form also appears in \cite{charb-popov}.

The first property of $\kappa_c$ is the fact that it is
\emph{closed}. Indeed, the $r_i$'s are not just any functions. The
Weyl rule (section \ref{sec:analyse}, formula (\ref{equ:weylcommut}))
applied to the commutation property of the \pdo s $P_i$ is, at the
order $h^2$, equivalent to~:
 \[ \{r_i,p_j\} = \{r_j,p_i\}. \]
In a symplectic chart given by the  Darboux-Carath\'eodory theorem,
this reads~:
\[ \deriv{r_i}{x_j} = \deriv{r_j}{x_i},\] and $\kappa_c$ is the closed
1-form $\kappa_c=-\sum r_idx_i$. \cqfd

\subsubsection{Deformation of Lagrangean submanifolds}
\label{sec:deformation}
Let $\mathcal{L}(M)$ be the set of all Lagrangean submanifolds of a
$2n$ dimensional symplectic manifold $M$. By the Darboux-Weinstein
theorem, we can identify a tubular neighborhood of a Lagrangean
submanifold $\Lambda_0$ with a neighborhood of the zero section of
$T^*\Lambda_0$. Any other Lagrangean submanifold within that
neighborhood can then be identified with a closed 1-form on
$\Lambda_0$ (see \cite{weinstein-symplectic}). This gives a ``chart''
for $\mathcal{L}(M)$, making it formally a differentiable infinite
dimensional manifold, whose tangent space at $\Lambda_0$ is naturally
identified with the space of closed 1-forms on $\Lambda_0$.

One can be more specific. The exponential map for vector fields gives a
diffeomorphism between a tubular neighborhood of $\Lambda_0$ and the
normal bundle $\frac{T_{\Lambda_0}M}{T\Lambda_0}$. The latter is
identified with $T^*\Lambda_0$ by means of the symplectic form~:
let $x\in\Lambda_0$; to any $X\in T_xM$ is associated the cotangent
vector 
\[ \tilde{\omega}(X)=i_X\omega_{\restr T_x\Lambda_0}. \]
Since $\Lambda_0$ is Lagrangean, the kernel of $\tilde{\omega}$ is
exactly $T_x\Lambda_0$.

In this way, we can easily describe the tangent space of
$\mathcal{L}(M)$~: an infinitesimal variation of $\Lambda_0$ is by
definition a vector field transversal to $\Lambda_0$, that is, a
section of $\frac{T_{\Lambda_0}M}{T\Lambda_0}$. By the above
isomorphism, it can be identified with a 1-form on $\Lambda_0$. That
infinitesimal variation is performed within the space of Lagrangean
submanifolds if and only if this 1-form is closed, that is, the
deformation vector field is locally Hamiltonian (see also
\cite{weinstein-phases}).

Let us apply this description to the case where a path in
$\mathcal{L}(M)$ is given as a 1-parameter family of level sets of $n$
independent functions $p_i^t$ in involution~:
\[ \Lambda^t = \{ p_i^t=0. \} \]
\begin{lemm}
  \label{lem:kappa}
  The infinitesimal variation 1-form $\kappa=\ddt{\Lambda^t}_{|t=0}$
  is given by 
  \[ \kappa(\ham{p_i^0}) = - \deriv{p_i^t}{t}_{|t=0}. \]
\end{lemm}
\demo Let $x\in\Lambda^0$. Complete $(p_1^0,\dots,p_n^0)$ into a
Darboux-Carath\'eodory chart, so that the tangent vectors
$\deriv{}{p_j^0}$ form a basis of the normal bundle. Taylor's formula
$p_i^t \sim p_i^0 + t\deriv{p_i^t}{t}_{|t=0}$ shows that at $t=0$, the
deformation vector field for $\Lambda^t$ is given by
\[ X(x) = \sum_{j} -  \deriv{p_i^t(x)}{t}_{|t=0} \deriv{}{p_j^0}. \]
To this transversal vector field corresponds the 1-form $\kappa =
i_X\omega$. On the basis $(\ham{p_1^0},\dots,\ham{p_n^0})$ of
$T\Lambda^0$, it is given by $\kappa(\ham{p_i^0}) =
\omega(X,\ham{p_i^0}) = dp_i^0.X = - \deriv{p_i^t}{t}_{|t=0}$. \cqfd

\begin{rema} For each $t$, the functions $p_1^t,\dots,p_n^t$ define
  a \cis\  in a neighborhood of $\Lambda^t$. The variation 1-form can
  be computed not only for $\Lambda^t$, but also for the neighboring
  leaves $\Lambda^t_c\egdef \cap_i (p_i^t)^{-1}(c)$. Then
  $\kappa_c=\ddt{\Lambda_c^t}_{\restr t=0}$ is given on $\Lambda_c$ as
  in lemma \ref{lem:kappa}~:
 \[ \kappa_c(\ham{p_i^0}) = - \deriv{p_i^t}{t}_{|t=0}. \]
The fact that for all $c$, $\kappa_c$ is closed can be checked by
  differentiating the equality $\{p_i^t,p_j^t\}=0$ at $t=0$.
\end{rema}
\begin{rema} It is natural now to be willing to consider the level sets
\[ p_i^0 + t\deriv{p_i^t}{t}_{\restr t=0}  = 0\]
as a ``linearization'' of the family $\Lambda_t$. Unfortunately, these
level sets are in general not Lagrangean. More precisely, the
symplectic form restricted to them is of order $t^2$. However, one can
prove that there exists a set of functions
$\tilde{r}_1,\dots,\tilde{r}_n$ such that the functions
\[ \tilde{p}^t_i\egdef p^0_i+t\tilde{r}_i \]
are indeed in involution, and such that the corresponding deformation
1-forms $\tilde{\kappa}_c$ satisfies
\[ [\tilde{\kappa}_c] = [\kappa_c] \in H^1(\Lambda_c^0).\]
The value of $\tilde{r}_1$ on $\Lambda_c$ is constant and obtained by
\emph{averaging} the sub-principal terms $\deriv{p_i^t}{t}_{\restr
t=0}$ over $\Lambda_c$.
%we have the following result~:
%\begin{lemm}
%  \label{lem:poincareglobal}
%  There exists a \emph{linear} 1-parameter family of $n$ independent
%  functions $\tilde{p}_i^t$ in involution, with variation 1-form
%  $\tilde{\kappa}$, such that
%  \[ \tilde{p}_i^0=p_i^0, \quad \textrm{ and } \quad
%  [\tilde{\kappa}] = [\kappa]. \]
%%  There exists a 1-form $\kappa'$ cohomologous to $\kappa$ in
%%  $H^1(\Lambda_0)$ and such that the level sets
%%\[ \Lambda_t' = \{ p_i^0 - t\kappa'(\ham{p_i^0})  = 0 \} \] are Lagrangean.
%\end{lemm}
%\demo Let $(x,\xi)$ be action-angle coordinates in a neighborhood of
%$\Lambda_0$ and denote by $-r_i(')$ the coordinates of $\kappa(')$
%with respect to the basis $(\deriv{}{x_1},\dots,\deriv{}{x_n})$ of
%$T\Lambda_0$. 

%The lemma will be proved as soon as we find a function
%$f$ such that the 1-form $\kappa'=\kappa+df$ on $\Lambda_0$ satisfies
%$\forall i,j,\quad \{r_i',p_j^0\}=0$. On the basis
%$(\ham{p_i^1},\dots,\ham{p_i^n})$, this reads~:
%\[ \deriv{f}{x_i} = r_i - r_i',  \textrm{ and } \deriv{r_i'}{x_j}=0.\]
%The first equation has a solution if and only if
%$\int_0^1(r_i-r_i')dx_i=0$. We can now set $r_i' = \int_0^1
%r_idx_i$. One easily verifies that such a $r_i'$ does not depend on
%$x$. This in turn implies that the $r_i'$'s commute with each other,
%and thus answer the question. \cqfd
\end{rema}

\subsubsection{Deformation of the Action integral}
Let $\Lambda_t$ be a smooth family of Lagrangean submanifolds of $M$,
and let $\gamma_t(\theta), \theta\in S^1$ be a smooth family of loops
such that each $\gamma_t$ is drawn on $\Lambda_t$. Suppose that on a
neighborhood of the image of $\gamma_0$, the symplectic form $\omega$
is exact~: $\omega=d\alpha$. Then this holds for $t$ small enough, and
we can define the action integral~:
\[ A(\gamma_t) = \int_{\gamma_t} \alpha. \]
\begin{lemm}
  \label{lem:action}
  The variation 1-form $\kappa=\ddt{\Lambda_t}_{\restr
    t=0}$ on $\Lambda_0$ is characteristic of the infinitesimal
  variation of the action, in the following sense~:
  \[ \ddt A(\gamma_t)_{\restr t=0} = \int_{\gamma_0}\kappa. \]
\end{lemm}
\demo We want to prove that the 1-forms
$\ddt\left(\gamma_t^*(\alpha)\right)_{\restr t=0}$ and
$\gamma_0^*(\kappa)$ are cohomologous on $S^1$. We have
\[ \ddt\left(\gamma_t^*(\alpha)\right) =
\gamma_t^*\left(\lie_{\deriv{\gamma_t}{t}} \alpha\right) =\]
\[ = \gamma_t^*\left(i_{\deriv{\gamma_t}{t}} d\alpha +
  d(i_{\deriv{\gamma_t}{t}} \alpha)\right). \] At $t=0$, the vector
field $\deriv{\gamma_t}{t}$ splits into two components $X_\tau$ and
$X_\nu$, the first one  being tangent to $\Lambda_0$. The other
one, in the normal bundle $\frac{T_{\Lambda_0}M}{T\Lambda_0}$, is by
definition the deformation vector field of the family
$\Lambda_t$. Therefore,
\[ \gamma_0^*\left(i_{\deriv{\gamma_t}{t}} d\alpha\right) =
\gamma_0^*(i_{X_\tau}\omega) +
\gamma_0^*(i_{X_\nu}\omega) = 0 + \gamma_0^*(\kappa), \]
which gives the result. \cqfd

\begin{rema}
  Of course, the value $t=0$ plays a arbitrary role~: if we define
  $\kappa_{t} = \ddt\Lambda_t$, then $\ddt A(\gamma_t) =
  \int_{\gamma_t}\kappa_t$.
\end{rema}

\subsubsection{$\kappa$ as a semi-classical deformation}
\label{sec:kappa}
Returning to the hypothesis of paragraph (\ref{sec:sous-principale}),
we are now able to give a geometrical interpretation of the
sub-principal form.

The total symbol $p_i(h)$ of $P_i(h)$ is considered as a
semi-classical deformation of the principal symbol $p_i^0$. Then
around any regular point of the principal symbol, the sub-principal
form is equal to the deformation 1-form -- as $h$ tends to zero -- of
any family of Lagrangean submanifolds defined as level sets of
functions of the form $p_i^0+hr_i + O(h^2)$.

Of course, the level sets of the total symbols themselves need no be
Lagrangean. Recall however that in the case where $0$ is a regular value
of the principal symbol we can average the sub-principal symbols
without changing the cohomology class of $\kappa_c$ in such a way that
the level sets of $p_i^0+hr_i$ are indeed Lagrangean.

As we shall see in section \ref{sec:CBSreg}, the \BS\ quantization 
conditions that apply to such a situation depend on the sub-principal 
symbols only through the cohomology class of $\kappa_c$.  This shows 
that articles like \cite{colinII,charbonnel} that would rather assume 
that the sub-principal symbols are either equal to zero or at least 
constant on each $\Lambda_c$ are not far from the general case.

%Let us get back to the notations of (\ref{sec:sous-principale}), so that
%we see the total symbol $p_i(h)$ of $P_i(h)$ as a semi-classical
%deformation of the principal symbol $p_{i,0}$. We suppose here that
%the principal Lagrangean is non-singular. We also assume that the
%sub-principal symbols are real-valued. The level sets $\cap
%(p_i(h))^{-1}(0)$ need not be Lagrangean, but the previous lemmas give
%the following interpretation of $\kappa$~:
%\begin{prop}
%  \label{prop:deformation}
%  The sub-principal form $\kappa$ is the infinitesimal variation
%  1-form, as $h\fleche 0$, of any family of Lagrangean submanifolds
%  defined as level sets of functions of the form $p_{i,0}+hr_i +
%  O(h^2)$.

%  Moreover, there exists, in the same cohomology class as $\kappa$, a
%  1-form $\kappa'$ that gives a Lagrangean linearization of the total
%  symbol $(p_1(h),\dots,p_n(h))$.
%\end{prop}
%Let us give a last result stressing the relevance of the cohomology
%class $[\kappa]$~:
Another nice property of the cohomology class of $\kappa_c$ is its
invariance under conjugation by unitary \pdo s.
\begin{prop}
  Let $P_i(h)$ be a semi-classical integrable system with a
  non-singular principal Lagrangean $\Lambda_0$. Let $V(h)$ be a
  classical unitary \pdo\ on a neighborhood of $\Lambda_0$. Let
  $P_i'=V^{-1}P_iV$ be the transformed integrable system
  and $\kappa'$ the new sub-principal form.

  Then $[\kappa']=[\kappa]$.
\end{prop}
\demo Since $P$ and $P'$ have same principal symbols, we can write
$P'=P+hQ$, with $Q=(Q_1,\dots,Q_n)$. The intertwining property
$VP'=PV$ reads~:
\[ [P_j,V] = hVQ_j .\]
At the principal level, this gives $\frac{1}{i}\{p_j,v\} = vq_j$,
which, writing $v=e^{ic}$, yields~:
\[ \{p_j,c\} = q_j. \] 
The sub-principal symbols $r_j^{(')}$ of $P^{(')}$ satisfy
$r_j'=r_j+q_j$, that is,
\[ \kappa' = \kappa + dc, \] where $dc$ denotes the differential of
$c$ as a function on $\Lambda_0$. \cqfd 

Note that the unitarity of $V$ is only used up to $O(h^2)$ in order to
ensure that the sub-principal symbols of the transformed system are
real-valued.

\section{Microlocal analysis}
\label{sec:analyse}

The aim of this section is to present a minimal version of the \mi\
analysis needed in order to give a precise and usable definition of
what we call ``\mi\ solutions'' of pseudo-differential systems
(section \ref{sec:solutions}).

The results can be retrieved from
\cite{colin-p,courscolin,robert,FIO1}; however, note that we focus
here a little bit more on the \emph{microlocal} point of view, and
assumptions about ``behavior at infinity'' are usually irrelevant in
this theory.

\subsection{Symbols}
\label{sec:symbols} We recall here the basic definitions of symbols
and classical symbols. One reason for including here this standard
material, besides the sake of completeness, is that our work mainly
deals with local properties, in a neighborhood of a point or at least
in a neighborhood of a compact region of the phase space; and for this
purpose, a very simple notion of 'symbols' can be introduced, which is
the following~:
\begin{defi} A family of complex-valued functions $(p(h))_{h\in H}$ on
  a manifold $Z$, where $H$ is a subset of $\RM$ having $0$ as an
  accumulation point, is called a \emph{compactly supported symbol} if
  there exists a compact set $K\subset Z$ and an integer $m\in\ZM$
  such that every $p(h)$ is a $\Cinf$ function with support in $K$ such
  that~: 
  \begin{equation}
    \label{equ:estimation}
    \forall \alpha, \exists C_\alpha, \textrm{ s.t. } \forall h\in H,
    \sup_{z\in Z}|\partial_z^\alpha p(z;h)| \leq C_\alpha h^m.
  \end{equation}
  We denote by $\SY^m_0(Z)$ the vector space of all such symbols.
  If we drop the condition on the support we obtain the space of all
  \emph{symbols} of order $m$, denoted by $\SY^m(Z)$.
\end{defi}
The notation $\SY^m_{(0)}(Z)$ shall be used in any assertion that holds 
for $\SY^m_{0}(Z)$ as well as for $\SY^m(Z)$.
\begin{rema}
  In the notation $\SY^m_{(0)}(Z)$, we forgot the dependence in $H$. Somewhat
  later we will have to be more careful, but this omission is harmless
  at this point.
\end{rema}

Note that $\SY^m_{(0)}(Z)\subset \SY^{m'}_{(0)}(Z)$ for $m\geq m'$ and,
due to the Leibniz formula, the space $\SY^*_{(0)}(Z)=\bigcup_{m\in\ZM}
\SY^m_{(0)}(Z)$ of symbols of any order is a graded algebra for the
usual multiplication. It is also clear that $\SY^m_{(0)}(Z)=h^m
\SY^0_{(0)}(Z)$. If $Z$ sits in a symplectic manifold, we can also endow
$\SY^*_{(0)}(Z)$ with the Poisson bracket $\{p_1(h),p_2(h)\}$ turning it
into a graded Poisson algebra.
\begin{defi}
  A symbol $(p(h))_{h\in H}\in \SY^m_{(0)}(Z)$ is said to
  be \emph{classical} if it admits an asymptotic expansion of the
  form~:
  \[ p(z;h) \sim \sum_{k\geq m} p_k(z)h^k, \] in the sense that
  $\forall N\geq m$, $p(h)-\sum_{k=m}^N p_kh^k \in \SY^{N+1}_{(0)}(Z)$.
\end{defi}
\begin{rema}
  If we apply the above definitions to symbols independent of
  $z$, we get the notion of a ``constant symbol''. We will let
  $\CM_h\inject \SY^*(Z)$ denote the algebra of such constants. With
  that respect, the spaces $\SY^*(Z)$ and $\SY^*_0(Z)$ can be
  considered as modules over $\CM_h$.
\end{rema}

% for any domain $D$ containing
%$\overline{\Omega}$ in its interior, and for any
%$p(h)\in\SY^0(\Omega)$, there is a $\tilde{p}(h)\in
%\SY^0_0(\Omega)$ such that $p(h)=\tilde{p}(h)$ on $D$.

A functional calculus can be performed on symbols of non-negative 
order; we will only need the following~: if $p(h)\in \SY^m_{(0)}(Z)$ 
with $m\geq 0$, then $\exp p(h)\in \SY^0(Z)$. 
% It is of course not 
%compactly supported.  One can easily check that for any cut-off 
%function $\chi$ supported in $K$ we have $\chi\exp p(h)\in \SY^0_0(Z)$, 
%and 
Conversely, let $q(h)$ be an element of $\SY^0_{(0)}(Z)$; if $\Omega$ is
any proper simply-connected open subset of $Z$ on which $q(h)$ is
\emph{elliptic} (\emph{i.e.} $|q(h)|>c>0$ on $\Omega$, for all $h$),
then one can define a symbol $p(h)\in \SY^0_0(Z)$ such that $q(h)=\exp
p(h)$ on $\Omega$.

Everything can be restated in the ``classical'' category.

\begin{defi}
  We say that a symbol is $O(h^\infty)$ if it is in $\SY^m(Z)$ for all
  $m>0$. We denote by $\SY^\infty(Z)$ the subspace of all such
  symbols. Similarly, $\SY_0^\infty(Z)$ is the intersection of all
  $\SY^m_0(Z)$.
\end{defi}
Now, the spaces that will be really of interest for us are the
quotient spaces $\SY^*_{(0)}(Z)/\SY^\infty_{(0)}(Z)$. Elements in the same
class will be called \emph{microlocally equal} on $Z$. This word is
easily justified when $Z$ is a subset of a $T^*X$; its extension to
any $Z$ should raise no problem. The space of classical symbols modulo
microlocal equality can then be isomorphically identified with the
space of formal expansions of the form $\sum_{k\geq m}
p_k(z)h^k$. Note also that the classical constant symbols modulo
microlocal equality form the field of numerical formal series in
$h$. In general, as it is easily verified, the quotient of $\CM_h$ by
microlocal equality is a field that we will denote by $\bar{\CM_h}$.

In order to get meaningful information from the ``quantization
conditions'' that we are going to derive later, we need to consider
symbols depending uniformly on a parameter $E\in\RM^n$. By uniformly,
we mean that for every compact subset $B$ of $\RM^n$, the estimate
(\ref{equ:estimation}) is valid uniformly for $E\in B$. In the classical
case, this is true for instance if each $p_k$ depends continuously on $E$.

\subsection{Pseudo-differential operators}
For the general theory, we refer to \cite{robert,courscolin}. Proofs
can usually be derived from the homogeneous theory of H\"ormander
\cite{hormanderIV}. We wish here to present a microlocal version of
the standard classes of pseudo-differential operators. Let $X$ be a
differentiable manifold of dimension $n$, and $M=T^*X$. In every local
coordinates, a compactly supported pseudo-differential operator
$P(h)$ on an open subset $U$ of $X$ is an operator with smooth kernel
of the form~:
\begin{equation}
  \label{equ:quantif}
   K_h(x,y) = \frac{1}{(2\pi h)^n}\int_{\RM^n}
   e^{\frac{i}{h}(x-y).\xi}a(x,y,\xi;h)d\xi,  
\end{equation}
where $a(h)$ is a symbol in $\SY^*_0(U\times U\times\RM^n)$.  Note that 
$K_h$ has compact support in $U\times U$, hence these operators are 
continuous linear operators (for fixed $h$) from $\Cinf(U)$ to 
$\Cinf_0(U)$, or even from distributions $\mathcal{D}'(U)$ to 
$\Cinf_0(U)$.  We shall assume that a $\Cinf$ half density 
$|dx|^{\frac{1}{2}}$ is given on $X$ so that each function $u$ on $X$ 
is associated with the half density $u|dx|^{\frac{1}{2}}$.

Pseudo-differential operators have a Weyl symbol, which is a
$h$-dependent function $p(h)=\sigma_W(P(h))$ on $M$, such that $P(h)$ can
be retrieved from $p(h)$ by the so-called ``Weyl quantization''
scheme~:
\begin{equation}
  \label{equ:Weyl}
  (Op^W(p(h))u)(x) = \frac{1}{(2\pi h)^n}\int e^{\frac{i}{h}(x-y).\xi}
  p(\frac{x+y}{2},\xi;h)u(y)dyd\xi.
\end{equation}
%On can show (\cite[courscolin,???]) that $p(h)$ is given by $K(h)$
%through the formula~:
%\[ p_h(x,\xi) = \frac{1}{(2\pi h)^n}\int
%K_h(x+y,x-y)e^{\frac{i}{h}2y.\xi}dy. \]
 
Now let $p(h)\in \SY^m_0(\Omega)$, and $\chi\in\Cinf_0$ equal to $1$ in
some neighborhood of $0\in\RM^n$. We have
$a(x,y,\xi;h)=p(\frac{x+y}{2},\xi;h)\chi(x-y)\in \SY^m_0(U\times
U\times\RM^n)$, which allows us to form a \pdo\  $Op^W_\chi(p(h))$ via
(\ref{equ:quantif}). Then one can show that the Weyl symbol of the latter
\pdo, although perhaps not of compact support, is nevertheless
microlocally equal to $p(h)$ on $\Omega$. That motivates the following
definitions~:
\begin{defi}
  Let $\Omega$ be an open subset of $M$, in which we have symplectic
  coordinates $(x,\xi)$. We suppose here that $x$ varies in $U$, an
  open subset of $X$ containing $\pi(\Omega)$ where $\pi$ is the
  natural projection $T^*X\fleche X$.  A compactly supported \pdo\  $P(h)$
  on $U$ will be called $h$-smoothing or ``$O(h^\infty)$'' on $\Omega$
  or element of $\Psi^\infty_0(\Omega)$ if its Weyl symbol belongs to
  the space $\SY^\infty(\Omega)$.
\end{defi}
The quantization procedure $Op^W_\chi(p(h))$ is therefore independent
of $\chi$ modulo $\Psi^\infty_0(\Omega)$.
\begin{defi}
  We denote by $\Psi^m_0(\Omega)$ the space of compactly supported
  \pdo s of the form $Op^W_\chi(p(h)) + \Psi^\infty_0(\Omega)$, for any
  $p(h)$ and $\chi$ defined as above.
\end{defi}
\begin{rema} 
Dealing with \pdo s, we shall always assume that the
symbols involved are classical.
\end{rema}

This definition leads to an isomorphism between the spaces
$\Psi^m_0(\Omega)/\Psi^\infty_0(\Omega)$ and
$\SY^m_0(\Omega)/\SY^\infty_0(\Omega)$ (for the latter is isomorphic to
$\frac{\SY^m_0 + \SY^\infty}{\SY^\infty}$).  This isomorphism, still denoted
by $\sigma_W$, depends on the local coordinates, but we recall that
the first two terms in the asymptotic expansions of the Weyl symbols
are intrinsically defined on $T^*X$, provided we let \pdo s act on
half-densities (the proof of this was given in \cite[prop 5.2.1]{FIO2}
for homogeneous \pdo s, \emph{i.e.} without a small parameter $h$, and
applies to our situation with no essential change); these two terms
are respectively called the \emph{principal} and \emph{sub-principal}
symbols, and are compactly supported.

It also ensures that theses classes of \pdo s are stable with respect to
the operations of composition and taking
adjoints.

We are now able to define the notion of the restriction of a \pdo;
namely, if $N$ is any subset of $\Omega$, two elements $P(h)$ and $Q(h)$
of $\Psi^m_0(\Omega)$ will be called microlocally equal on $N$, and
written ``$P\sim Q$ on $N$'' if $\sigma_W(P)=\sigma_W(Q)$ on
$N$. Then we can perform inversion of elliptic operators in the
following sense: $P(h)\in\Psi^m_0(\Omega)$ is elliptic at a point $m$
if its principal symbol does not vanish at $m$. Then there is a
neighborhood $N$ of $m$ and a $Q(h)\in\Psi^{-m}_0(\Omega)$ such that
$PQ\sim QP\sim I$ on $N$, where $I$ is a \pdo\  satisfying
$I\sim Op^W(1)$ on $N$.

\begin{rema} 
The space of Weyl symbols is naturally equipped with the Lie
algebra structure given by the symplectic Poisson bracket. On the
other hand, the space of \pdo s has a natural Lie operator algebra
structure. Though it is known that there is no hope for finding a
``quantization'' isomorphism that fully respects this Lie algebra
structure, the Weyl quantization has the not-so-bad following
behavior (see \cite{flato})~: Let $p_1(h)$ and $p_2(h)$ be
symbols of order 0, and $P_1(h)$ and $P_2(h)$ their
Weyl-quantization. Then $[P_1,P_2]$ is a \pdo\ of order 1 and of
Weyl symbol 
\begin{equation}
  \label{equ:weylcommut}
  \frac{h}{i}\left(\{p_1,p_2\} + O(h^2)\right).
\end{equation}
In particular the sub-principal symbol of $[P_1,P_2]$ is easily
computed in terms of those of $P_1$ and $P_2$~:
\[ \ssub(\frac{i}{h}[P_1,P_2]) = \{p_{1,0},\ssub(P_2)\} +
\{\ssub(P_1),p_{2,0}\}, \]
where $p_{j,0}$ is the principal symbol of $P_j$.
\end{rema}
\begin{rema}
The above classes of operators would be purely abstract if we could
not relate them to usual $h$-admissible \pdo s in the sense of
\cite{robert} and others, for which additional assumptions concerning
the behavior of the symbols at infinity are needed. Such a link is
here easy: if $P(h)\in\Psi^m_0(\Omega)$, where $\Omega$ is some open
subset of $\RM^n$, we can form a strongly admissible \pdo\  on $\RM^n$ by
extending the Weyl symbol to zero outside its support, and using a
$Op^W_\chi$ quantization. Conversely, if $Q(h)$ is an admissible \pdo\  in
the sense of \cite{robert}, it has a Weyl symbol, that we can cut to
make it compactly supported in $\Omega$, and get a $P(h)$ again. If
$Q(h)$ was already obtained by a $\tilde{P(h)}$ in $\Psi^m_0(\Omega)$,
then by definition the new $P(h)$ will microlocally coincide with
$\tilde{P(h)}$ on the set where the cut-offs are both equal to 1.  The
interest of this microlocalization procedure will be made clear in the
next section (prop \ref{prop:microlocalisation}).
\end{rema} 

Finally, compactly supported \pdo s on manifolds can then be defined to
be locally finite sums of operators admitting the above description in
a local coordinate chart. The principal and sub-principal symbols
remain well-defined functions on $T^*X$.

\subsection{Fourier integral operators}
\label{sec:FIO}
Similarly to the space $\Psi^*_0(\Omega)$, microlocal classes of
Fourier integral operators can be constructed. We will here skip most
of the details. Let $X$, $Y$ be differentiable manifolds, and
$\Omega_X$ and $\Omega_Y$ open subsets of $T^*X$ and $T^*Y$
respectively. Suppose we have a symplectic diffeomorphism $\chi$ from
$\Omega_Y$ to $\Omega_X$. Let $\Lambda$ be the graph of $\chi$ in
$\Omega_X\times \Omega_Y$. It is an immersed Lagrangean submanifold of
$T^*X\times\overline{T^*Y}$, where $\overline{T^*Y}$ denotes the
cotangent $T^*Y$ equipped with $-\omega$, the opposite of the standard
symplectic form. Thus, $\Lambda$  admit parameterization by
non-degenerate phase functions $\phy$~: $\Lambda = \cup
\Lambda_\phy$. For details about this construction, refer to
\cite{duistermaat-oscillatory} or \cite{wein}. Following
\cite{courscolin}, we define $\Psi^m_0(\Omega_X,\Omega_Y,\chi)$ to be
the space of operators whose kernel is microlocally equal to a
compactly supported classical $h$-oscillatory integral on $\Lambda$
(see paragraph \ref{sec:oscillatory}). On each $\Lambda_\phy$, it has a
principal symbol defined as a function on $\Lambda_\phy$. A \fio\  $U(h)$
is said to be elliptic at a point $(x,\xi)$ or, equivalently, at the
point $\chi^{-1}(x,\xi)$, if its principal symbol does not vanish at
$(\chi^{-1}(x,\xi),(x,\xi))\in\Lambda_\phy$.

These operators behave as expected with respect to composition and
microlocal inversion of elliptic operators. Most important for us,
they allow to find normal forms for pseudo-differential operators
thanks to the celebrated theorem of Egorov~: if $P(h)$ is a \pdo\  on
$\Omega_X$ of principal symbol $p$ and $U(h)$ an \fio\  elliptic at a
point $(x,\xi)$, then $U^{-1}PU$ is, near $\chi^{-1}(x,\xi)$, a
\pdo\  of principal symbol $p\circ\chi$.

\subsection{$h$-admissible functionals}
In the above, we have let \pdo s act on $h$-independent distribution
half-densities, which is of course not enough for our purposes. Let
$U$ be a regular domain in $\RM^n$. If
$u\in\mathcal{D}'(U,\Omega_{\frac{1}{2}})$, then for any \pdo\  $P(h)$,
$Pu$ is a compactly supported distribution, hence in the Sobolev
space $H^{-s}$ for some integer $s$; this holds for any fixed $h$. In
order to get asymptotic information for all $h\in H$, we are lead to
the following definition (see \cite[prop IV-8]{robert})~:
\begin{defi}
  \label{defi:admissible}
  A family $(u_h)_{h\in H}$ of distribution half densities is called
  \emph{admissible} on $\Omega$ if for any \pdo\ 
  $P(h)\in\Psi^0_0(\Omega)$, there is an $N\in\ZM$ and some $s\in\ZM$
  such that $h^N (Pu_h)$ is uniformly bounded in $H^{-s}(U)$ for all $h\in
  H$. We denote by $\mathcal{D}_h(\Omega)$ this space of admissible
  functionals. 
\end{defi}

$\mathcal{D}_h(\Omega)$ is then a $\CM_h$-module that is by
definition stable under the action of $\Psi^m_0(\Omega)$. 
\begin{prop}
  If $u_h$ is admissible on $\Omega$ then  for any
  $P(h)\in\Psi^0_0(\Omega)$ there is an $N'\in\ZM$ such
  that $h^{N'}Pu_h$ is bounded in $L^2(U)$.
\end{prop}
\demo[(see \mbox{\cite[p.195-196]{robert}})] One can find an elliptic
  operator $D(h)$, that is a zero-order usual \pdo\ on $U$, that
  uniformly maps $L^2(U)$ into some power of $h$ times $H^{-s}(U)$
  (extend to $\RM^n$, and take $D=h^{s}(-\Delta+|x|^2+1)^{s/2}$). Let
  $D^{-1}$ be a right inverse of $D$, and $P(h)\in\Psi^0_0(\Omega)$ be
  as in definition \ref{defi:admissible}. So there is a $N'$ and a
  $C<\infty$ such that $\|D^{-1}(h)P(h)u_h\|_{L^2(U)} < Ch^{-N'}$. Now
  let $\Omega'\supset\overline{\Omega}$ and let
  $I(h)\in\Psi^0_0(\Omega')$ have symbol equal to $1$ on a
  neighborhood of $\Omega$. Then
  $(I(h)-1)P(h)\in\Psi^\infty_0(\Omega)$, and by the symbolic
  calculus, $I(h)D(h)\in\Psi^0_0(\Omega')$. Therefore, $I(h)D(h)$ is
  uniformly $L^2$-continuous (\cite[th.II-36]{robert}), and we can
  write $\|P(h)u_h\|_{L^2}=\|I(h)D(h)D^{-1}(h)P(h)u_h\|_{L^2} + \oh
  \leq C'h^{-N'}$. \cqfd

A natural notion of ``microlocal equality'' for admissible functionals
can now be defined~:
\begin{defi}
  \label{defi:equality}
  $\mathcal{D}^\infty_h(\Omega)$ denotes the space of
  $u_h\in\mathcal{D}_h(\Omega)$ such that for any \pdo\ 
  $P(h)\in\Psi^*_0(\Omega)$ we have~:
  \[ \|Pu_h\|_{L^2(U)} = O(h^\infty). \]
  Two admissible functionals $u_h$ and $v_h$ are called microlocally
  equal on $\Omega$ if they belong to the same class modulo
  $\mathcal{D}^\infty_h(\Omega)$. Following \cite{colin-p}, we will
  write ``$u_h\sim v_h$ on $\Omega$'' in that case.  If $m\in T^*X$,
  we say that $u_h\sim v_h$ at the point $m$ if there exists an open
  neighborhood $\Omega$ of $m$ such that $u_h\sim v_h$ on $\Omega$.
\end{defi}
Note that, because the above definition only involves estimates on
compact subsets of $\Omega$, $u_h\sim v_h$ on $\Omega$ if and only if
$u_h\sim v_h$ at each point of $\Omega$. Moreover, to test microlocal
equality at a point, it is sufficient to pick up an elliptic operator
that satisfies the required estimate. There again we have obvious
although useful properties~:
\begin{prop}
    \begin{itemize}
  \item if $P(h)\in\Psi^\infty_0(\Omega)$, then for all
    $u_h\in\mathcal{D}_h(\Omega)$, $Pu_h\sim 0$ on $\Omega$;
  \item if $u_h\sim 0$ on $\Omega$, then for all
    $P(h)\in\Psi^*_0(\Omega)$, $Pu_h\sim 0$ on $\Omega$.
  \end{itemize}
\end{prop}
This is essentially due to the $L^2$-continuity of \pdo\ of order $0$.
In other words, the quotient space
$\Psi^*_0(\Omega)/\Psi^\infty_0(\Omega)$ acts naturally on
$\mathcal{D}_h(\Omega)/\mathcal{D}^\infty_h(\Omega)$.
This says that the action of a \pdo\  on an admissible functional is
microlocally given by the (formal) Weyl symbol $\sigma_W$.
Note also that the space
$\mathcal{D}_h(\Omega)/\mathcal{D}^\infty_h(\Omega)$ has the structure
of a $\bar{\CM_h}$-vector space.

In particular, if $N$ is any subset of $\Omega$, the definition
\ref{defi:equality} gives a natural microlocal equality on $N$, which
is compatible with the notion of restriction of a \pdo\ previously
introduced:~ two \pdo s microlocally equal on $N$ have the same action
on two admissible functionals microlocally equal on $N$.

The same remark applies for the following interesting result~:
\begin{prop}
  \label{prop:microlocalisation}
  If $Q(h)$ is a global admissible \pdo\ in $\RM^n$ in the sense of
  \cite{robert}, $\Omega$ an open subset of $\RM^{2n}$, $u_h$ an
  admissible functional in $\Omega$, and $P(h)\in\Psi^*_0(\Omega)$ a
  microlocalization of $Q$ in a compact $K$ of $\Omega$, then
  $Qu_h$ is admissible in $\Omega$ and~:
 \[ Qu_h \sim Pu_h \textrm{ on } K. \]
\end{prop}

Alternatively, one can relate this notion of microlocal equality for
admissible functionals to the so-called \emph{semi-classical wave
front set}, as introduced in \cite{robert,courscolin}, and denoted by
$WF_h(u_h)$. This notion is not stabilized yet, in the sense that,
depending on authors, it includes or not uniform estimates at infinity
(see remark in \cite[p.1541]{colin-p}). Anyhow, the simpler part of
$WF_h$, that is to say its intersection with $T^*X$, is defined as
follows~:

suppose $u_h$ is admissible on $T^*X$. Then $WF_h(u_h)\cap T^*X$ is
the complement subset in $T^*X$ of the biggest open subset $\Omega$
such that $u_h\sim 0$ on $\Omega$.

The definition is of course based upon the homogeneous wave front set
introduced by H\"ormander (see \cite[sec.2.5]{FIO1}). It has a useful
local characterization in terms of the semi-classical Fourier
transform~:
\begin{lemm}[\cite{robert,courscolin}]
  \label{lem:local}
  An admissible functional $u_h$ microlocally vanishes at the point
  $(x_o,\xi_o)\in T^*X$ (equivalently, $(x_o,\xi_0)\not\in WF_h(u_h)$)
  if and only if there exists a function $\phy\in\Cinf_0(X)$ with
  $\phy(x_o)\neq 0$ such that, in some local coordinates~:
  \[ \fourier(\phy u_h)(\xi) = O(h^\infty) \]
  uniformly for $\xi$ in a neighborhood of $\xi_o$.
\end{lemm}
Here $\fourier$ denotes the ``semi-classical Fourier transform''~:
\[ \fourier = \frac{1}{(2\pi h)^{\frac{n}{2}}}\int_{\RM^n}
e^{-\frac{i}{h}x.\xi} dx \] (we will sometimes use the usual Fourier
transform $\fouriero=\fouriero_1$). Note that $\fourier$ is \emph{a
priori} defined for tempered distributions, but the Fourier transform
of a compactly supported distribution in $\RM^n$ can be thought as a
continuous function (it is analytic). For the proof of this lemma, one
constructs a compactly supported \pdo\ with ``rectangular support'',
that is with symbol of the form $\chi(\xi)\phy(x)$.

Unfortunately, because of the lack of uniform estimates in the fibers
one should not expect an analogue of theorem 2.5.3 of the previously
cited article \cite{FIO1}. Indeed, if $U$ is an open set of $X$, then
$WF_h(u_h)\cap T^*U = \emptyset$ does \emph{not} imply that $u_h$ is
locally $O(h^\infty)$ on $U$. For instance, one easily sees from lemma
\ref{lem:local} that $u_h(x)=a(x)e^{ix/h^2}$ is microlocally zero on
$T^*\RM^n$, while obviously not locally $O(h^\infty)$.

The solution found in \cite{courscolin} for functionals on $\RM^n$ is
to extend the wave front set to a subset of the compactified cotangent
bundle obtained by adding a point at infinity to each direction $\xi$.
An admissible functional $u_h$ is then said to be microlocally zero in
the direction $(x_0,\RM^+\xi_0)$ if, for large $\xi$ in a conic
neighborhood of $\xi_0$, and for any cut-off function $\phy$, the
Fourier transform~:
\[ \fourier (\phy u_h)(\xi) \]
is of order $O(h^N/|\xi|^N)$ for all $N>0$. Then the following
statement is valid~:
\begin{lemm}
  \label{lem:global}
  An admissible functional $u_h$ is $O(h^\infty)$ uniformly for $x$
  near $x_o$ if and only if  there exists a function $\phy\in\Cinf_0(X)$
  with $\phy(x_o)\neq 0$ such that~:
  \[ \fourier(\phy u_h)(\xi) = O(h^N/(1+|\xi|)^N), \]
  for all $\xi\in\RM^n$, for all $N\in\NM$.
\end{lemm}
The condition expressed by this lemma is that there exists a
neighborhood $U$ of $x_o$ such that $u_h$ is microlocally zero in
$T^*U$ as well as in all directions of $T^*U$.

Similarly to lemma \ref{lem:local}, we have for instance the following
fact~: let $u_h$ be an admissible functional defined on $\RM^n$, so
that $\fourier u_h$ makes sense. If $\fourier u_h$ is a $O(h^\infty)$
function in a neighborhood of any $\xi$ in an open set $\Upsilon$,
then $u_h$ is microlocally zero in the whole
$\RM^n\times\Upsilon$. This allows the following construction~: let
$1=\chi_1+\cdots+\chi_k$ be a partition of unity subordinated to the
open cover $\bigcup_j\Upsilon_j$ in the $\xi$-space. Let $U$ be an
open subset of $\RM^n$ with compact closure. If $v_{j,h}$,
$j=1,\dots,k$ are admissible functionals on $U\times\Upsilon_j$ that
microlocally coincide on non-empty intersections
$U\times\Upsilon_{j_1}\cap\cdots\cap\Upsilon_{j_l}$, one can define an
admissible functional $u_h$ on $U\times\bigcup_j\Upsilon_j$ by cutting
off the $v_{j,h}$'s outside $\overline{U}$ in $\RM^n$ and requiring~:
\[ \fourier u_h = \sum_{j=1}^k \chi_j\fourier v_j. \] Then $u_h$
microlocally coincides with each $v_{j,h}$ on
$U\times\Upsilon_j$. Indeed, from the microlocal equality of the
$v_{j,h}$'s on non-empty intersections, and the compactness of their
support, it follows from lemma \ref{lem:local} that the Fourier transforms
$\fourier v_{j,h}$ are $\Cinf$ functions that coincide with each other
up to $O(h^\infty)$ around each point of such non-empty
intersections. Thus $\fourier u_h$ also has that common value around
these points, which yields the result. \cqfd

Finally, if $U(h)\in\Psi^m_0(\Omega_X,\Omega_Y,\chi)$, then for any
admissible $u_h$ in $\Omega_Y$, $Uu_h$ is admissible in $\Omega_X$,
and because of Egorov's theorem, $U$ transports the wave front set
by $\chi$.

\subsection{Microlocal solutions}
\label{sec:solutions}
\begin{defi}
  Let $\Omega$ be an open set in $T^*X$.  An admissible functional
  $u_h\in\mathcal{D}_h(\Omega)$ is called a \emph{microlocal solution}
  on $\Omega$ of a system $(P_1(h),\dots,P_k(h))$ of
  pseudo-differential operators if~:
  \[ \forall j=1,\dots,k, \qquad P_ju_h \sim 0 \textrm{ on } \Omega. \]  
\end{defi}
We allow here to choose a suitable set $H$ where $h$ is to vary.  The
space of all microlocal solutions of $(P_1,\dots,P_k)$ modulo
microlocal equality will be viewed as a $\bar{\CM}_h$-vector space.
Because of microlocal inversion of elliptic operators, any such
solution is ``microlocalized'' in $\bigcap_j p_j^{-1}(0)$, which means
that $WF_h(u_h)$ must lie in this set ($p_j$ denotes the principal
symbol of $P_j$). If two solutions have different non-empty wave front
sets, then they are independent.  In case the $p_j$'s are real-valued, the
``propagation of singularities'' theorem states that $WF_h(u_h)$ is locally
stable under the action of the Hamiltonian flow of $p_j$.

As we will see later, global problems often have no solution in the
previous sense.  Suppose that the $P_j$ depend continuously on an
additional parameter $E$ in a \emph{compact} topological space $V$
(cf. end of section \ref{sec:symbols}). Depending on the degeneracy of the
$p_j$'s, it might be impossible to find a solution (in the previous
sense) for some $E$, and depending on $E$, the appropriate sets $H_E$
where $h$ varies might be disjoint. Still, the whole theory being
asymptotic, we need to control things when $h$ tends to zero !
Therefore, we introduce the following~:
\begin{defi}
  \label{defi:admissibleE}
  Let $\Gamma$ be a subset of $]0,1]\times V$. Let $H$ be the
  projection of $\Gamma$ onto $]0,1]$. We assume that $0$ is an
  accumulation point of $H$. We say that a family
  $(u_{(h,E)})_{(h,E)\in\Gamma}$ is admissible if definition
  \ref{defi:admissible} holds uniformly for $(h,E)\in\Gamma$.
\end{defi}
Recall that $E$ varies in a compact.  Accordingly, $u_{(h,E)}$
microlocally vanishes at point $m$ if there is an elliptic \pdo\
$P(h)$ at $m$ such that $\|Pu_{(h,E)}\|_{L^2(U)}= O(h^\infty)$
uniformly for $(h,E)\in\Gamma$. Thus we are still able to define a
semi-classical wave-front set $WF_h(u_{(h,E)})$.

Now if $u_{(h,E)}$ is admissible and $P_j^E$ depend continuously
on $E$, then $P_j^Eu_{(h,E)}$ is admissible, and we say that
$u_{(h,E)}$ is a microlocal solution of the system
$(P_1^E,\dots,P_k^E)$ on $\Omega$ if $P_j^Eu_{(h,E)}\sim
0$ on $\Omega$ according to that definition.
Then we have~:
\begin{prop}
  \label{prop:wf}
  Let $u_{(h,E)}$ be a microlocal solution of the equation
  $P^Eu_{(h,E)}\sim 0$ on $\Omega$. The following inclusion holds~:
  \[ WF_h(u_{(h,E)})\subset \{ (x,\xi)\in T^*X, \quad \exists E_o\in V_o,
  \quad p_{E_o}(x,\xi)=0\}, \] where $V_o$ is the set of $E_o\in V$
  such that $(0,E_o)$ is an accumulation point of $\Gamma$. $p_{E}$ is
  the principal symbol of $P^E(h)$.
\end{prop}
\demo Let $(x,\xi)\in T^*X$ such that for any $E_o\in V_o$, the
principal symbol $p_{E_o}$ is non zero at $(x,\xi)$. Since $V_o$
is compact and $p_{E}$ is continuous in $E$, then there is an open
neighborhood $\mathcal{V}$ of $V_o$ and a constant $c>0$ such that
for $E\in\mathcal{V}$, $|p_ {E}|>c$. Then there is an open
neighborhood $\Omega$ of $(x,\xi)$ on which that remains valid. The
complement $\mathcal{V}^c$ being a compact set of non-accumulation
points, we can pick up a $h_o$ such that $([0,h_o]\times
\mathcal{V}^c)\cap\Gamma$ is empty. From now on, we will therefore
restrict $\Gamma$ to $[0,h_o]\times\mathcal{V}$. Because of the
uniform ellipticity of $P^E$ on $\Omega$, for $(h,E)$ in that
new $\Gamma$, we can construct a \pdo\  $Q^E(h)$ depending
continuously on $E\in\mathcal{V}$ such that $Q^EP^E\sim I$
on $\Omega$. It is not difficult to see that if an admissible
$v_{(h,E)}$ is microlocally zero at a point $m$, then for any \pdo\ 
$Q^E(h)$ depending continuously on $E$ in some compact, then
$Q^Ev_{(h,E)}$ is still microlocally zero. Back in our problem, we
see that~:
\[ u_{(h,E)} \sim 0 \textrm{ at } (x,\xi), \]
which is what we needed to prove. \cqfd

\subsection{Oscillatory integrals}
\label{sec:oscillatory}
The WKB ansatz for solving \schr -type equations consists of
restricting the scope of solutions to a subspace of admissible
functionals~: the space of oscillatory integrals, that we wish to
describe now. These functions have a long history, particularly in
Quantum mechanics, and it is well known since the treatise of Maslov
\cite{maslov}, made more rigorous and expounded by H\"ormander and
Duistermaat (see in particular \cite{duistermaat-oscillatory}), that
they are particularly fit to locally solve generic microlocal
pseudo-differential equations. In fact, as we shall see later, in some
cases all solutions are microlocally equal to oscillatory integrals.

We recall here the definitions of \cite{duistermaat-oscillatory}, with
some precisions.

We are given an compact immersed Lagrangean manifold $\Lambda\in 
T^*X$, endowed with a smooth half-density $\rho$.  Suppose that is 
fixed a covering of $\Lambda$ with simply connected embedded open 
subsets $\Lambda_{\phy_k}$ described by reduced phase functions 
$\phy_k(x,\theta)$ (``reduced'' means that the number of additional 
oscillatory variables $\theta$ is exactly the maximum dimension of the 
kernel of $d\pi$, where $\pi$ is the projection $\Lambda_\phy \fleche 
X$.  This dimension is known to be the minimum dimension of $\theta$ 
variables needed for a phase function defining $\Lambda$ -- \cite[theo 
3.1.4]{FIO1} or \cite[1.3.6]{duistermaat-oscillatory}).  $\phy$ gives 
rise to a function on $\Lambda_\phy$ -- that we still denote by $\phy$ 
-- whose exterior differential is the Liouville 1-form $\alpha_0$.
 
An \emph{oscillatory symbol} on $\Lambda_\phy$ is a half-density
$\sigma(h)$ on $\Lambda_\phy$ of the form~:
\[ \sigma(\lambda;h) = c(h)e^{\frac{i}{h}\phy(\lambda)} a(\lambda;h)
  \rho(\lambda), \] with $c(h)\in \CM_h$, and $a(h)$ a \emph{classical}
symbol in $\SY^*(\Lambda_\phy)$. Will will say that $\sigma$ is
classical if $c$ is.

If $a(h)$ is of order $0$ and elliptic on $\Lambda_\phy$ (as we will
generally assume), we will prefer the notation~:
\begin{equation}
  \label{equ:symbolexp}
  \sigma(\lambda;h) = e^{ic(h) + i\Phi(\lambda;h)}\rho(\lambda),
\end{equation}
where $c(h)\in\CM_h$, and $\Phi(h)$ is a $\SY^1$-valued classical symbol on
$\Lambda_\phy$ of order $-1$, whose principal symbol $\frac{\Phi_{-1}}{h}$
satisfies~:
\[ d\Phi_{-1} = \alpha_0. \]

As usual, we denote by $j_\phy$ the local diffeomorphism from $C_\phy$ to
$\Lambda_\phy$ given by $(x,\theta)\mapsto (x,\deriv{\phy}{x}(x,\theta))$,
where $C_\phy=\{(x,\theta)\in X\times R^N,
\deriv{\phy}{\theta}(x,\theta) = 0 \}$.

An \emph{oscillatory integral} on $\Lambda_\phy$ is a 
half-density on $X$ of the form~:
\begin{equation}
  \label{equ:integral}
  u_h(x) = \frac{c(h)}{(2i\pi h)^\frac{N}{2}} \int_{\RM^N}
e^{\frac{i}{h}\phy(x,\theta)}b(x,\theta;h) d\theta|dx|^{\frac{1}{2}},
\end{equation}
where $c(h)\in \CM_h$, and $b(h)$ is a classical symbol in
$\SY^*_0(X\times\RM^N)$. Then
\[\sigma={(j_\phy)^{-1}}^*(ce^{\frac{i}{h}\phy}{b}_{\restr
C_\phy})\] is an oscillating symbol.  $u_h$ is an admissible
functional, and because of the stationary phase formula, is uniquely
defined modulo microlocal equality by the asymptotic expansion of
$b(h)$ on any neighborhood of $C_\phy$. $\sigma$ is thus well
defined by $u_h$ modulo microlocal equality. 
%That this definition is independent of the choice of a
%reduced phase function is essentially a result of H\"ormander.
\begin{defi}
  An admissible functional $u_h$ is a \emph{Lagrangean distribution}
  on $\Lambda$ if there is a cover $\Lambda=\cup \Lambda_\phy$ such
  that $u_h$ is microlocally equal to an oscillatory integral on each
  $\Lambda_\phy$.
\end{defi}
On each $\Lambda_\phy$, such Lagrangean distributions have a symbol
$\sigma(h)$ defined modulo $O(h^\infty)$, but in general ($N\neq 0$)
the morphism $u_h\mapsto \sigma(h)$ is not injective; however, one can
prove that it is an isomorphism at the \emph{principal} level (if
$u_h$ can be written with $b(h)$ of order $m$ such that the $h^m$ part
of $\sigma(h)$ vanishes on a neighborhood of
$(x_o,\xi_o)\in\Lambda_\phy$, then near $j_\phy^{-1}(x_o,\xi_o)$,
$u_h$ can be rewritten with a ${b(h)}'$ of order $m+1$).

On the other way round, the symbol map is surjective, and one can
construct right inverses (``Maslov canonical operators'') in the
following way~: let $\sigma(h)$ be a oscillatory symbol near
$(x_o,\xi_o)\in\Lambda_\phy$. Then it defines the function $b(h)$ on a
neighborhood of $b_o=j_\phy^{-1}(x_o,\xi_o)$ in $C_\phy$. If $N>0$,
$C_\phy$ is a proper submanifold of $X\times\RM^N$. Let $B_o$ be a
sub-bundle of $T(X\times\RM^N)$ supplementary to $TC_\phy$ near $b_o$.
One defines the germs of $b$ on $C_\phy$ to
be zero along $B_o$, and smoothly extends $b$ to $X\times\RM^N$
accordingly. Then $b$ has the same $h$-order as $\sigma(h)$ and 
the oscillatory integral $u_h=Op(\sigma)$ so constructed has symbol
$\sigma(h)$ (mod $O(h^\infty)$).

Of course, since we have not made any mention to the Keller-Maslov
bundle yet, the symbol of a Lagrangean distribution depends on the
phase function $\phy$ considered. The result of H\"ormander-Maslov is
the following~: if $u_h$ and $v_h$ are zeroth order oscillatory
integrals on $\Lambda_\phy$ and $\Lambda_\psi$ respectively that are
equal modulo $O(h)$ at a point $(x,\xi)\in
\Lambda_\phy\cap\Lambda_\psi$, then there is an integer $\mu$ such
that, at the principal level,
%\footnote{the higher order terms needed in that correction are
%not easy to describe. However, one should probably be able to fix a
%suitable set of canonical operators $Op_\phy$ subordinated to a fixed
%cover $\Lambda=\cup\Lambda_\phy$ such that $Op_\phy(\sigma_\phy)\sim
%Op_\psi(\sigma_\psi)$ if and only if $\sigma_\phy=
%e^{i\mu\frac{\pi}{2}} \sigma_\psi + O(h^\infty)$.}
 $\sigma_{u_h} \sim
e^{i\mu\frac{\pi}{2}} \sigma_{v_h} + O(h)$ at $(x,\xi)$. Because of
the normalization in (\ref{equ:integral}), $\mu=\mu_\psi-\mu_\phy$, where,
for a phase function $g$, $2\mu_g$ is the signature of the
$\theta$-Hessian of $g$ at $(x,\xi)$ minus $N_g$.

With these transition functions, one can globally define the principal
symbol of a Lagrangean distribution as a section of the Keller-Maslov
bundle. In the sequel, Maslov indices will always be prompted very
clearly, so we don't insist any further on that point.

Note that if $u_h$ serves to construct a \fio\  $U(h)$ as
mentioned in section \ref{sec:FIO}, then the principal symbol of $u_h$ is
invariantly defined by $U$.

If the symbol is elliptic and we write it as (\ref{equ:symbolexp}), we
will refer to $\Phi(h)$ as the \emph{total phase} of $u_h$, and the
\emph{principal phase} is $\frac{\Phi_{-1}}{h}+\Phi_0$, both defined
modulo a constant $c(h)\in\CM_h$.

Finally, we admit the dependence of a Lagrangean distribution on an
additional parameter $E$, provided the involved estimates are uniform
in the sense of definition \ref{defi:admissibleE}. It is the case for
instance if the symbol $b$ depends continuously on $E$.

\section{Non-singular quantization conditions}
\label{sec:CBSreg}

Let $P_1(h),\dots,P_n(h)$ be a semi-classical integrable system, and
denote by $\Lambda_0$ the principal Lagrangean manifold associated to
it. $\Lambda_0$ may be critical, but in this section we shall work on
a submanifold $\Lambda\subset\Lambda_0$ consisting only of
non-singular points of the momentum mapping
$(p_1,\dots,p_n)$. Moreover, we
shall always assume that $\overline{\Lambda}$ is compact.

With these data, we express a necessary and sufficient condition for
the existence of microlocal solutions of the system on $\Lambda$
(proposition \ref{prop:cocycle}). Thanks to the WKB construction, we
prove that this condition involves, at least at the principal level,
geometrical characteristics of the system, namely the Liouville 1-form
$\alpha_0$ and the sub-principal form $\kappa$ (theorem
\ref{theo:quantization}). 

As an application, we recover the well-known Maslov-\BS\ quantization
conditions for $O(h^2)$ solutions in the case where
$\Lambda=\Lambda_0$ is a whole non-singular Liouville torus.

When the system depends smoothly on an $n$-dimensional parameter $E$
in a suitable sense, the above results are seen to be uniform on $E$,
which will allow us to apply the quantization condition to the study of
semi-classical spectra. In the case of a non-singular Liouville torus,
the study of the spectrum was also completely carried out, with a
different approach, in \cite{charbonnel}.

\subsection{Dimension of microlocal solutions}

The following result is a generalization of lemma 18 of
\cite{colin-p}.
\begin{prop}
\label{prop:dim1reg}
  Let $x\in\Lambda$ be a non-singular point of $\Lambda_0$. The space
  of microlocal solutions of the system $P_iu \sim 0$ on a
  neighborhood of $x$, modulo microlocal equality, has dimension 1.
  
  If $x\not\in\Lambda_0$, the space of solutions of this system has
  dimension 0.
\end{prop}
\demo
The second point comes from the fact that near $x$, the system is
elliptic, hence microlocally invertible (see proposition
\ref{prop:wf}).

Concerning the first one, we use the following microlocal normal form,
which is standard~:
\begin{lemm}
  \label{lem:fn0}
  There exists an elliptic \fio\ $U(h)$ associated to a
  Darboux-Carath\'eodory coordinate chart $(x,\xi)$ such that
  \[ {U}^{-1}P_jU \sim \frac{h}{i}\deriv{}{x_j} \]
  on a neighborhood of $x$. If every $P_j$ is formally
  self-adjoint, then $U$ can be chosen to be unitary (mod
  $O(h^\infty)$). 
\end{lemm}
Accordingly transformed, the system always admits the constant
solutions. It remains now to prove that these are the only ones, that
is, that any solution $u_h$ of $\deriv{}{x_j}u\sim 0$ is microlocally
constant. Let $u_h$ be such a solution. We can always assume that it
is compactly supported. We want now to ``cut-off in the frequencies'',
which we do in the following way~:

Let $v_h$ an admissible functional defined by $\fourier v_h =
\chi\fourier u_h$, where  the cut-off function $\chi$ is $\Cinf_0$ and
has value 1 on a neighborhood of $0$. Then $u_h$ and $v_h$ are
microlocally equal on a neighborhood of $0$. Moreover, $v_h$ is a
\emph{local} solution (modulo $O(h^\infty)$) of the system near $0$.
This comes from the fact that
$\fourier(\frac{h}{i}\deriv{v_h}{x_j})=\xi_j \fourier v_h$ satisfies
the hypothesis of lemma \ref{lem:global}.

We have thus reduced the problem to a system of the form~:
\[ d_hv_h = w_h = O(h^\infty) \] on a neighborhood of $0$ (here $d_h$
denotes the exterior differential operator $d_h =
(\frac{h}{i}\deriv{}{x_1}, \dots, \frac{h}{i}\deriv{}{x_n})$).
The solutions of this system are necessarily of the form 
\[ cste_h + \frac{i}{h}\int_0^1 w_h(tx).x dt. \]
The second term in the sum is $O(h^\infty)$, as $w_h$ was. Hence it is
 microlocally zero, which gives the result. \cqfd

Let $\phy_t=\phy_1(t_1)\circ\cdots\circ\phy_n(t_n)$ be the classical
flow associated to the momentum map of the system. By transport of
singularities, we obtain the following corollary~:
\begin{coro}
  \label{coro:dim1reg}
  If the action of $\phy_t$ on $\Lambda$ is transitive, then
  the space of microlocal solutions of the system $P_iu \sim 0$ on
  $\Lambda$ has dimension at most 1.
\end{coro}

If we restrict to a simply connected open subset $\Omega$ of
$\Lambda_0$, and if all trajectories of the flow remain open (which is
automatic if $\Omega$ contains no singular point), then local
microlocal solutions can be glued together to obtain a non-trivial
solution of the system on $\Omega$. This occurs for instance if
$\Omega$ is obtained by propagating a non-singular point of
$\Lambda_0$ by the flow during a short time (cf. theorem 1.4.1 of
\cite{duistermaat-oscillatory}).

As soon as we have closed trajectories, we get obstructions, which we
describe in the following paragraph. Later on (paragraph
(\ref{sec:BSconditions})), we will show how to interpret these
obstructions thanks to the WKB method, and get the so-called
``\BS\  quantization conditions''.

\subsection{The bundle of microlocal solutions}
\label{sec:bundle}
For $x\in\Lambda$, let $\mathcal{L}_h(x)$ be the 1-dimensional vector
space over $\bar{\CM_h}$ of microlocal solutions of the system at the
point $x$ (modulo microlocal equality). For any open subset $\Omega$
of $\Lambda$, the disjoint union $\mathcal{L}_h(\Omega)$ of the
$\mathcal{L}_h(x)$, $x\in\Omega$ has a natural structure of a
$\bar{\CM_h}$-line bundle over $\Omega$, defined by the following
local trivializations~: if $\Omega=\cup\Omega_\alpha$, where every
$\Omega_\alpha$ is chosen small enough so that the system admits a
non-trivial microlocal solution $u_h^\alpha$ on it, then we define
$\phy_{\Omega_\alpha}$ to be the isomorphism~:
\[  \phy_{\Omega_\alpha} : \Omega\times\bar{\CM_h} \fleche
\mathcal{L}_h(\Omega) \]
\begin{equation}
  \label{equ:trivialize}
  (x,c(h)) \mapsto c(h)u_h^\alpha.
\end{equation}
If we restrict to non-trivial solutions, this also defines a principal
bundle $\dot{\mathcal{L}}_h(\Lambda)$ over $\Omega$ with structure
group ${\bar{\CM_h}}^*= \bar{\CM_h}\setminus\{\bar{0}\}$. This bundle
has a natural flat connection, whose parallel sections are those that
are locally constant in the above trivializations (from the
infinitesimal viewpoint, its horizontal distribution is defined to be
the image of the standard horizontal spaces of
$\Omega\times\bar{\CM_h}$ by the trivialization maps
$\phy_{\Omega_\alpha}$). It is well defined because on intersections
$\Omega_\alpha\cap\Omega_\beta$, such trivializations differ by a
constant, due to proposition \ref{prop:dim1reg}. By the same argument, it
is also independent of the set of microlocal solutions $u_h^\alpha$
chosen to define the trivializations. The interest of this connection
is stressed by the following lemma~:
\begin{lemm}
  \label{lem:parallel}
  There is a natural identification of the space of non-trivial
  microlocal solutions of the system on an open subset
  $\Omega\subset\Lambda$ with parallel sections of
  $\dot{\mathcal{L}}_h(\Omega)$.
\end{lemm}
\demo
Let us first consider the case $\Omega=\Omega_\alpha$. The general case
will follow by a microlocal partition of unity. Now, using the
trivialization $\phy_{\Omega_\alpha}$, a section of
$\mathcal{L}_h(\Omega_\alpha)$ takes the form~:
\[ \sigma_h(x) = c(x;h)u_h^\alpha. \]
This section is parallel if and only if $c(x;h)=c(h)$ is a constant
(modulo microlocal equality, of course), so, by proposition
\ref{prop:dim1reg}, if and only if $\sigma_h(x)=\sigma_h$ is a microlocal
solution on $\Omega_\alpha$. It is non-trivial if and only if $c_h\neq
0$.

Now assume $\Omega$ is any open subset of $\Lambda$. One sense is
obvious~: if $u_h$ is a non-trivial solution on $\Omega$, then it can
be used to trivialize the line bundle $\mathcal{L}_h(\Omega)$ as in
(\ref{equ:trivialize}). Now, since the definition of the connection is
independent of the chosen solutions, the constant section $\sigma_h(x)
= u_h$ is a parallel section.

On the other way round, suppose $\sigma_h(x)$ is a parallel section of
$\dot{\mathcal{L}}_h(\Omega)$. Then we proved that its restrictions to
each $\Omega_\alpha$ are non-trivial microlocal solutions $v_h^\alpha$
on $\Omega_\alpha$. Let $P^\alpha(h)$ be a pseudo-differential
partition of unity on $\Omega$ subordinated to the cover
$\Omega_\alpha$ (it means that the Weyl symbols form a $\Cinf_0$
partition of unity of $\Omega$). Then each $P^\alpha v_h^\alpha$ is
an admissible functional on a neighborhood of $\Omega$, and the sum
\[ v_h = \sum_\alpha P^\alpha v_h^\alpha \] 
is microlocally equal to $v_h^\alpha$ (and hence to $\sigma_h$) on each
$\Omega_\alpha$. Therefore, $v_h$ is a microlocal solution on
$\Omega$, whose associated parallel section is $\sigma_h$. This
terminates the proof of the lemma. \cqfd

Thus, the obstruction to the existence of a microlocal solution on
$\Lambda$ only lies in the \emph{holonomy} of the connection, which,
in this case, describes the monodromy of the solution. This
holonomy can be represented by the cohomology class of a \v Cech
1-cocycle on $\Lambda$ with values in ${\bar{\CM_h}}^*$. Summarizing,
we get~:
\begin{prop}
  \label{prop:holonomy}
  There exists a non-trivial microlocal solution of the system on
  $\Lambda$ if and only if the holonomy cocycle of the
  ${\bar{\CM_h}}^*$-principal bundle $\dot{\mathcal{L}}_h(\Lambda)$
  with its natural flat connection has microlocally trivial
  cohomology.
\end{prop}
From now on, our principal task will be to explicit that
holonomy. Besides the above definition, we can give a direct and more
concrete description of it in terms of \v Cech cohomology, that shall
be more efficient for our further analysis.
 
We still let $\Omega_\alpha$ be a (finite) cover of $\Lambda$ by simply
connected open subsets such that each $\Omega_\alpha$ admits a
neighborhood in $M$ on which proposition \ref{prop:dim1reg} applies.
We thus get a family  $u_h^\alpha$ of non-trivial microlocal solutions
of the system on $\Omega_\alpha$. Moreover, on non-empty intersections
$\Omega_\alpha\cap\Omega_\beta$, the same proposition affirms the
existence of unique non-zero constants $c^{\beta\alpha}(h)\in\CM_h^*$
such that~:
\[ u_h^\beta \sim c^{\beta\alpha}u_h^\alpha \textrm{ on }
\Omega_\alpha\cap\Omega_\beta. \]
The $c^{\beta\alpha}$'s satisfy the cocycle condition~:
\[ c^{\beta\alpha} \sim c^{\beta\gamma}c^{\gamma\alpha} \]
on every non-empty intersection
$\Omega_\alpha\cap\Omega_\beta\cap\Omega_\gamma$. Furthermore, another
choice of solutions ${u'}_h^\alpha$ would give rise to constants
${c'}^{\beta\alpha}$ differing by a coboundary $d^\alpha(h)$~:
\[ {c'}^{\beta\alpha} \sim \frac{d^\beta}{d^\alpha}
c^{\beta\alpha}. \] This \v Cech cocycle $c^{\beta\alpha}$ is
exactly the holonomy form of the connection introduced above.  The
matching conditions are then the following~:

``there exists a non-trivial microlocal solution on $\Lambda$ if and
only if $c^{\beta\alpha}(h)$ is a coboundary (modulo $O(h^\infty)$).''

It is, of course, a simple re-formulation of proposition
\ref{prop:holonomy}. 
In the sequel, instead of the complex-valued \v Cech cocycle $[c(h)]$,
we will rather deal with the two real-valued cocycles obtained via the
exponential representation~:
\[ c^{\beta\alpha} = e^{\rho^{\beta\alpha}}
e^{i\lambda^{\beta\alpha}}, \] in which $\rho^{\beta\alpha}(h)\in\RM$
and $\lambda^{\beta\alpha}(h)\in \T$.
% (in some sense, we view the
%holonomy as a cocycle with values in the Lie algebra of
%$\bar{\CM}_h^*$).
 We denote by $\rho(h)$ and $\lambda(h)$ the
corresponding \v Cech cocycles.
\begin{prop}
  \label{prop:cocycle} Suppose that every $P_j$ is formally
  self-adjoint. Then $[\rho(h)]$ is microlocally trivial,
  and $[\lambda(h)]$ admits an asymptotic expansion of the form~:
  \[ [\lambda(h)] \sim \sum_{k\geq -1} [\lambda_k]h^k.\]
\end{prop}
\demo We cover $\Lambda$ by open subsets $\Omega_\alpha$ on which the
normal form of lemma \ref{lem:fn0} applies. This gives rise to a
family of microlocally unitary \fio s $U^\alpha(h)$ on
$\Omega_\alpha$. An associated set of solutions is given by $u_h^\alpha
= U^\alpha 1_\alpha$, where $1_\alpha$ denotes an admissible
functional microlocally equal to the constant $1$ on
$\Omega_\alpha$. Suppose $\Omega_\alpha\cap\Omega_\beta$ is a
non-empty intersection and let $x$ be an element of it. Finally, let
$1_x$ denote an admissible functional microlocally equal to $1$ on a
neighborhood of $x$. On $x$, we have $u_h^\alpha\sim U^\alpha 1_x$ and
$u_h^\beta\sim U^\beta 1_x$. Therefore, the constant
$c^{\beta\alpha}(h)$ is given by~:
\[  c^{\beta\alpha} 1_x \sim (U^\alpha)^{-1} U^\beta 1_x. \]
The unitarity of the \fio s implies that $|c^{\beta\alpha}|\sim 1$,
which gives the first statement of the proposition.

Now, $(U^\alpha)^{-1} U^\beta$ is an elliptic \fio\ which implies that
$(U^\alpha)^{-1} U^\beta 1_x$ is at $x$ a classical Lagrangean
distribution with a symbol of the form (\ref{equ:symbolexp}). This
gives the last statement of the proposition.  \cqfd

\begin{rema}
The description of the cocycle $c^{\beta\alpha}(h)$ involves no
regularity on $\Lambda$. And indeed, the notion of a flat line bundle with
holonomy is well-defined on any topological manifold. This will
enable us, in section \ref{sec:global}, to apply that description on a
singular Lagrangean. For the moment however, the $\Cinf$ structure of
$\Lambda$ makes it possible to apply the De Rham isomorphism and get a
smooth connection 1-form $\lambda(h)$ on $\Lambda$. In the next
sections, precisely, we show how the 1-forms $\alpha_0$ and $\kappa$
enter the picture. 
\end{rema}

\subsection{WKB method}
Although this might not seem obvious, propositions \ref{prop:holonomy}
and \ref{prop:cocycle} are essentially related to the celebrated WKB
construction.

The WKB method for solving pseudo-differential equations consists in
looking for a special kind of approximate solutions, called
oscillatory integrals (see section \ref{sec:oscillatory}). The use of
semi-classical oscillatory integral in mathematics was initiated by
Maslov, and pursued by different authors including Arnold, H\"ormander
and Duistermaat.  The result is that, provided we are given an
invariant half-density on $\Lambda$ (this will come from the Liouville
measure), we can always construct such solutions in a neighborhood of
a non-singular point of $\Lambda_0$. According to proposition
\ref{prop:dim1reg}, such a solution has to be the only one -- up to
multiplicative constant. This means that we can freely use the WKB
construction to look for general properties of our microlocal
integrable systems. This is the goal of the following paragraph.

Note that we previously constructed microlocal solutions via another
method, namely the use of a normal form. The link with oscillatory
integrals lies in the use of a \fio; more precisely, with the
notations of lemma \ref{lem:fn0}, the admissible functional $U1$ is a
Lagrangean distribution that is a microlocal solution of the system on
a neighborhood of $x$ ($1$ denotes a functional microlocally equal to
the constant 1 on a neighborhood of $x$).

The derivation of the transport equations in the WKB method for \cis s
is fairly standard now (see for instance \cite{colin-quasi-modes}),
and our aim here is only to briefly recall the idea. Of particular
importance is the entrance of the sub-principal form $\kappa$ into the
picture.

Let $u_h$ be a Lagrangean distribution on $\Lambda_0$ of order $0$. To
ensure that it is not trivial, we assume that it is elliptic at a point
$(x_0,\xi_0)\in\Lambda_0$, so that we use the notation
(\ref{equ:symbolexp}) for the symbol $\sigma(h)$ (we assume that we work
on a fixed $\Lambda_\phy$). The viewpoint here is to
find conditions on $\sigma$ in order that $u_h$ be a microlocal solution
of the \cis\  at $(x_0,\xi_0)$.

First of all, if one looks for a solution up to $O(h)$ terms, one
finds the condition that $\Phi_{-1}$ is indeed a phase function for
$\Lambda_\phy$, i.e.~: on $\Lambda_\phy$, 
\[ \frac{d\Phi_{-1}}{h} = \frac{\alpha_0}{h}.\]
 Note that
this first condition only involves the principal symbols of the
$P_j$'s.

Keeping $\Phi_{-1}$ fixed, one can look now for a solution up to
$O(h^2)$. Since we suppose that the sub-principal symbols are real
(which occurs for instance if the $P_j$'s are formally self-adjoint),
then this raises two necessary and sufficient conditions~: a) the
half-density $\rho$ must be \emph{invariant} with respect to the
Hamiltonian action of the system; b) on $\Lambda_\phy$, one must have
\[ d\Phi_0 = \kappa, \] where $\kappa$ is the
sub-principal form of the system.  Because of condition a), we can
assume that $\rho$ is the canonical Liouville half-density on
$\Lambda_0$. Condition b) implies that $\kappa$ must be closed, which
is already given by the fact that the $P_j$'s commute with each other
up to the order $O(h^2)$.  Finally, note that if the sub-principal
symbols are complex-valued, the conditions are still sufficient, but
not necessary.

The previous step generalizes by induction to provide solutions up to
any order. However, notice that the two conditions obtained above
depend only on the principal symbol of $u_h$ and on some closed
1-forms uniquely defined by the system. The next conditions involving
higher order terms of the symbol won't of course be invariantly
defined anymore. Still, there is a way to keep track of the total
symbol of $u_h$, and this is to fix a ``Canonical operator'' (see
section \ref{sec:oscillatory}) $u_h=Op(\sigma(h))$,
%, which amounts to restricting
%the subspace of solutions.
as stated in the following proposition~:
\begin{prop}
  \label{prop:WKB}
  Let $\Lambda_\phy$ be a neighborhood of $(x_0,\xi_0)$ in
  $\Lambda_0$ described by a phase function $\phy$. Let $Op$ be a
  Canonical operator on $\Lambda_\phy$. There exists
  closed real-valued 1-forms $\kappa_k$, $k=-1,0,1,\dots,\infty$ on
  $\Lambda_0$ such that~:

  The oscillatory integral $u_h=Op(e^{i\Phi(\lambda;h)}\rho(\lambda))$
  is a microlocal solution of the system up to the order $O(h^N)$ on
  $\Lambda_\phy$ if and only if the $\Phi_k$'s, $k=-1,\dots,N-2$
  satisfy~:
  \[ d\Phi_k = \kappa_k.\]

    $\kappa_{-1}$ is the Liouville form $\alpha_0$ restricted to
  $\Lambda_0$, and $\kappa_0$ is the sub-principal form $\kappa$. For
  $k>0$, $\kappa_k$ depends on the $Op$-quantization.
\end{prop}

A similar result holds if the $P_j$'s depend continuously on a
parameter $E$. Note that in \cite[p.38]{colin-quasi-modes} Y.Colin de
Verdi\`ere included the dependence on a parameter directly in the
definition of the oscillatory symbols.

\subsection{The regular quantization conditions}
\label{sec:BSconditions}
\begin{theo}
\label{theo:quantization}
  We still assume that the $P_j$'s are formally self-adjoint, and that
  $\Lambda$ is a regular open subset of $\Lambda_0$.

  Let $[\lambda(h)]\in H^1(\Lambda,\T)$ be, as in proposition
  \ref{prop:cocycle}, the holonomy phase of the \mi\ bundle
  $\dot{\mathcal{L}}_h(\Lambda)$. There exists a non-trivial \mi\
  solution of the system on $\Lambda$ if and only if 
  \[[\lambda(h)]\in H^1(\Lambda,2\pi\ZM) + O(h^\infty),\]
  and any other solution is then a multiple of it. Moreover,
  $[\lambda]$ admits the following expansion~:
  \[ [\lambda(h)] \sim [\mu]\frac{\pi}{2} + \sum_{k=-1}^\infty
  [\lambda_k]h^k,\] where $\mu$ is the Maslov cocycle,
  $\lambda_{-1}$ is the restriction of Liouville 1-form $\alpha_0$ on
  $\Lambda$, and $\lambda_0=\kappa_{\restr\Lambda}$ is the
  sub-principal form of the system.
\end{theo}
\demo The first points are direct consequences of propositions
\ref{prop:holonomy}, \ref{prop:cocycle} and \ref{prop:dim1reg}. 
It remains to prove
the statement about the principal terms of the asymptotic expansion of
$[\lambda]$. We cover $\Lambda$ by simply connected open subsets
$\Omega_\alpha$ such that each $\Omega_\alpha$ admits a representation
by a non degenerate phase function $\phy_\alpha$.
The most immediate consequence of proposition \ref{prop:WKB} is that
microlocal solutions as Lagrangean distributions always exist on
$\Omega_\alpha$; Let us denote by $u_h^\alpha$ such a set of
solutions, and by $\Phi^\alpha(h)$ their principal phases. Finally, we
let $\mu^{\alpha\beta}$ be the (integral) value of Maslov's cocycle
associated with that cover.

The equation
\[ u_h^\beta \sim c^{\beta\alpha}(h)u_h^\alpha \textrm{ on }
\Omega_\alpha\cap\Omega_\beta \]
implies that $\arg c^{\beta\alpha}$ is, mod $O(h)$, the difference of
the principal phases of $u_\psi$ and $u_\phy$, plus the value of the
corresponding Maslov cocycle. That is, 
\[ \arg c^{\beta\alpha}(h) = \Phi^\alpha(h) - \Phi^\beta(h)
+ \mu^{\alpha\beta}\frac{\pi}{2} + O(h). \] 
Because of proposition \ref{prop:WKB}, $\Phi^\alpha - \Phi^\beta$ is
the value of a \v Cech cocycle whose cohomology is represented, via the de Rham
isomorphism, by the cohomology of the smooth 1-form
$\frac{\alpha_0}{h} + \kappa$ restricted to $\Lambda$.
Therefore, $\arg c^{\beta\alpha}$ is represented, modulo $O(h)$, by
$\frac{[\alpha_0]}{h} + [\kappa] + [\mu]\frac{\pi}{2}$. Since $\arg
c^{\beta\alpha}$ is itself a representative of $[\lambda]$, the
theorem is proved. \cqfd

If $\Lambda=\Lambda_0$ is a whole Liouville torus, the principal
terms in the quantization condition yield the well-known
Maslov-\BS\  quantization condition~: ``For all closed paths
$\gamma$ in $\Lambda_0$ (or equivalently for the $n$ standard
generators of $H_1(\Lambda_0)$ obtained through action-angle
coordinates),
\begin{equation}
  \label{equ:BSreg}
  \frac{1}{h}\int_\gamma \alpha_0 + \int_\gamma \kappa +
\mu(\gamma)\frac{\pi}{2} \in 2\pi\ZM + O(h)\textrm{''}.
\end{equation}

\subsection{Spectral parameter dependence}
\label{sec:dependence}
The motivation of this paragraph is to make the above results suitable
for studying spectral problems, namely for investigating the
solutions of the microlocal system~:
\[ (P_1(h)-E_1)u_h \sim 0, \dots,  (P_n(h)-E_n)u_h \sim 0, \]
according to the appropriate definition of microlocal equality
(section \ref{sec:solutions}).

% begin change
In order to do this, the operators $P_j$ in the quantization
conditions have to be replaced by $P_j-E_j$.  In the following
analysis, we will need to state the problem in a slightly more general
setting, where we replace $P_j$ by an operator depending on $E$ in
essentially the same way as $P_j-E_j$~:
%end change
\begin{defi}
  \label{defi:parameter}
  Let $P_i(h)$ be a semi-classical integrable system on an 
  $n$-dimensional manifold $X$, with principal symbols $p_i$.  A 
  system $P_i^E(h)$ depending smoothly on the parameter $E\in\RM^n$ 
  will be called regular deformation of $P_i(h)$ if for each $E$, 
  $(P_1^E(h), \dots, P_j^E(h))$ is a semi-classical integrable system, 
  with principal symbols $p_i^E$, and there exists a $\Cinf$ function 
  $f^E(x)=f(E,x)$, $f:\RM^n\times\RM^n\fleche \RM^n$ such that 
  $f^0=id$ and $f(\cdot,0)$ is a local diffeomorphism of $(\RM^n,0)$, 
  such that:~ \[ (p_1^E,\dots,p_n^E) = f^E(p_1,\dots,p_n).  \]
\end{defi}
% We suppose here that a good microlocal functional calculus can be
% applied to $(P_1,\dots,P_n)$. With no assumptions on $P_j$, we shall
% restrict to polynomial functions $f^E$. When all $P_j$'s are
% self-adjoint however, we can use any $\Cinf$ functions.
% 
% This definition ensures that $P_i^E$ is still an integrable system,
% with momentum map $f^E(p_1,\dots,p_n)$. 
For such a regular deformation, the new Lagrangean
foliation $\Lambda_c^E$ is globally the same as the original one
$\Lambda_c$, with different labels~: $\Lambda_c=\Lambda^E_{f^E(c)}$.

Now we assume that there is a tubular neighborhood $\Omega$ of $\Lambda$
such that for $E$ in a sufficiently small compact neighborhood of $0$,
the intersections $\Lambda^E=\Lambda_0^E\cap\Omega$ are non-singular
submanifolds of $\Lambda_0^E$. Then it is easy to verify that every
result of this section extends to uniform solutions of the perturbed
system $P_i^E(h)$ on $\Omega$. In particular we obtain a perturbed
cocycle $[c^E(h)]$ that smoothly depends on $E$. It implies that the
terms in the asymptotic expansion in powers of $h$ of the cocycle
$[\lambda^E(h)]$, equal to the argument of $[c^E(h)]$, are locally
smooth but possibly multivalued functions of $E$.  At the principal
levels, this gives a parameter dependent version of the
\BS\  conditions (\ref{equ:BSreg})~: ``for any family of
loops $\gamma^E$ drawn on $\Lambda_0^E$,
\begin{equation}
  \label{equ:BSregE}
  \frac{1}{h}\int_{\gamma^E} \alpha_0 + \int_{\gamma^E} \kappa^E +
\mu(\gamma^E)\frac{\pi}{2} \in 2\pi\ZM + O(h)\textrm{''}.
\end{equation}

\begin{rema}
  \label{rem:dependence}
Let $P_i^E(h)$ be such a smooth deformation.  Then one can look at the 
family obtained by letting $E=ht$, where $t$ is a bounded parameter 
(let's assume for clarity in this remark that $n=1$; everything 
applies to the many-dimensional case as well).  Using that in every 
local coordinates, the Weyl symbols of $P_I^E$ is a smooth function of 
$E$, one can see that the operator \[\tilde{P}^t(h)\egdef P^{ht}(h)\] 
is still a classical \pdo, with uniform estimates for bounded $t$, but 
its principal and sub-principal symbols are of course \emph{not} the 
same as those of $P^E(h)$.  More precisely,  if we denote by $r^E$ the 
sub-principal symbol of $P^E$, then, because of Taylor's formula~:
\begin{equation}
  \label{equ:taylorE}
  f(E,x) = f(0,x) + E\deriv{f}{E}(0,x) + \cdots +
  E^N\deriv{\!^Nf}{E^N}(0,x) + O(E^{N+1}),
\end{equation}
which holds uniformly for bounded $E$, the principal and sub-principal 
symbols $\tilde{p}^t$ and $\tilde{r}^t$ of $\tilde{P}^t(h)$ satisfy~:
\begin{equation}
  \label{equ:symbolsE}
  \tilde{p}^t = p^0 = f^0(p), \quad\mathrm{and}\quad \tilde{r}^t =
r^0 + t\deriv{f^E}{E}_{\restr E=0}(p).
\end{equation}
Now, let us see what \BS\  conditions (\ref{equ:BSregE}) become.
First of all, (\ref{equ:symbolsE}) shows that the principal Lagrangean
$\tilde{\Lambda}_0$ of the new system is independent of $t$ and equal
to $\Lambda_0^0$. Let us call $\tilde{\kappa}^t$ the sub-principal
form on $\Lambda_0^0$. On the other hand, the family of Lagrangean
manifolds $\Lambda_0^E=\{f^E(p)=0\}$, with $E=ht$, depends smoothly on
$h$; let us call $\delta\kappa^t$ its deformation 1-form at $h=0$ (see
\ref{sec:deformation}), defined by~:
\[ \delta\kappa^t.\ham{f^0(p)} =  - t\deriv{f^E}{E}_{\restr
  E=0}(p). \]
Then paragraph \ref{sec:kappa} together
with (\ref{equ:symbolsE}) shows that, on $\Lambda_0^0$,
\[ \tilde{\kappa}^t = \kappa^0 + \delta\kappa^t. \]
Now, if $\gamma^E$ is a smooth family of loops drawn on $\Lambda_0^E$,
lemma \ref{lem:action} implies that, when $E=ht$,
\[ \int_{\gamma^E}\alpha_0 = \int_{\gamma^0}\alpha_0 + h\int_{\gamma^0}
\delta\kappa^t + O(h^2). \]
Since it is clear that $\int_{\gamma^E}\kappa^E =
\int_{\gamma^0}\kappa^0 + O(h)$, we get~:
\[ \frac{1}{h}\int_{\gamma^E} \alpha_0 + \int_{\gamma^E} \kappa^E =
\frac{1}{h}\int_{\gamma^0} \alpha_0 + \int_{\gamma^0} \tilde{\kappa}^t
+ O(h). \]
In other words, we exactly recover the \BS\  quantization
conditions we would have expected to find for the $t$-dependent system
$\tilde{P}^t(h)$. 

Note that this remark will mostly be used with $f(E,x)=x-E$, so that
$\tilde{r}^t=r^0-t$.
\end{rema}

\section{\BS\  for a \emph{focus-focus} singularity}
\label{sec:focus-focus}
\subsection{Integrable systems with a non-degenerate singularity}
Let $M$ be any symplectic manifold of dimension $2n$, and
$f_1,\dots,f_n$ a \cis\ on $M$.  The system is said to be singular at a
point $m\in M$ if $m$ is a critical point for the momentum map $F$. We
shall always suppose in this work that $m$ is of maximal corank as a
critical point, which means that each function $f_i$ is critical at
$m$. Without loss of generality, we shall also assume that for all
$i=1,\dots,n$, $f_i(m)=0$.

The map $\mathcal{H}$ that assigns to a smooth function its Hessian at
$m$ is a Lie algebra homomorphism with respect to the Poisson bracket
from the Lie algebra of smooth functions that are critical at $m$ to
the Lie algebra of quadratic forms on the tangent space $T_mM$.
Fixing any local set of Darboux coordinates, the latter is identified
with the space $\mathcal{Q}(2n)$ of quadratic forms on
$\RM^{2n}=\{(x,\xi)\}$. 

Every singular \cis\ thus gives rise to a linear sub-algebra $C_F$ of
$\mathcal{Q}(2n)$, namely the sub-algebra generated by
$\{\mathcal{H}(f_1),\dots,\mathcal{H}(f_n)\}$. Note that $C_F$ is
always Abelian.
\begin{defi}[Eliasson \cite{eliasson}]
  A singular \cis\ with momentum map $F$ admits $m$ as a
  \emph{non-degenerate} singularity if $C_F$ is a Cartan sub-algebra of
  $\mathcal{Q}(2n)$.
\end{defi}
Recall that a Cartan subalgebra of a semi-simple Lie algebra
$\mathfrak{a}$ is a commutative subalgebra $\mathfrak{c}$ that is equal
to its commutant and such that, for any element $H\in\mathfrak{c}$,
$ad_H$ is a semi-simple endomorphism of $\mathfrak{a}$.

From the work of Williamson \cite{williamson}, real Cartan
sub-algebras of $\mathcal{Q}(2n)$ are classified according to the
following scheme~:
\begin{theo}[Williamson]
  Let $\mathfrak{c}$ be a real Cartan sub-algebra of
  $\mathcal{Q}(2n)$.  Then there exists a symplectic linear change of
  coordinates in $\RM^{2n}$ and a basis $q_1,\dots,q_n$ of
  $\mathfrak{c}$ such that each $q_i$ is one of the following~:
  \begin{itemize}
  \item $q_i = x_i^2 + \xi_i^2$ (elliptic type)
  \item $q_i = x_i\xi_i$ (hyperbolic type)
  \item $\left.\begin{array}{c} q_i = x_i\xi_i + x_{i+1}\xi_{i+1} \\
      q_{i+1} = x_i\xi_{i+1} - x_{i+1}\xi_i \end{array}\right\}$
      (focus-focus type)
  \end{itemize}
\end{theo}

If $M$ has dimension 4, as we shall assume from now on, only four
combinations are possible.  We will restrict ourselves to the case of
a \emph{focus-focus} singularity.  Note that this is the only case
where $m$ is \emph{isolated} amongst the set of critical points of
$F$.  In particular, the singular foliation $\Lambda_c=F^{-1}(c)$ has,
for small $c$, only one singular leaf $\Lambda_0$, and if $\Omega$ is
a sufficiently small neighborhood of $m$ in $M$, then
$\Omega\setminus\{m\}$ is foliated by smooth Lagrangean submanifolds
of $M$.  If $\Lambda_0$ is connected, this also implies that the
critical value $0$ lies in the interior of the image of the momentum
map $F$.  The aim of the following paragraphs is to give a precise
description of the singular fibration on a full neighborhood of
$\Lambda_0$, generalizing the usual Arnold-Liouville theory for
regular fibers. The results are not new, since they are quoted in
\cite{zung2} and proposition \ref{prop:immersion} can in principle be
recovered from \cite{lerman-umanski-book}, but I believe that the
following presentation gives a clearer account of what is needed for
our purposes.

\subsection{Linear focus-focus}
\label{sec:linear}
We fix here the notations used when referring to a linear focus-focus
system on $\RM^4=\{(x,y,\xi,\eta)\}$. We shall always use Williamson
coordinates, in which the system is generated by the quadratic forms
$q_1=x\xi+y\eta$ and $q_2=x\eta-y\xi$. The Hamiltonian flows are
denoted by $\ham{1}$ and $\ham{2}$. It will be convenient to use polar
coordinates $(r,\theta)$ for $(x,y)$ and $(\rho,\alpha)$ for
$(\xi,\eta)$, in which $\ham{1}$ and $\ham{2}$ take a simpler form
(see figure \ref{fig:linearise}). Of interest are also the complex
representations $z_1=x+iy$ and $z_2=\xi+i\eta$. The flows of $\ham{1}$
and $\ham{2}$ are then respectively
\[ (t,(z_1,z_2))\mapsto (e^{t}z_1,e^{-t}z_2), \textrm{ and } \]
\begin{equation}
  \label{equ:flot}
  (s,(z_1,z_2))\mapsto (e^{is}z_1,e^{is}z_2).
\end{equation}
Similarly, we can identify the momentum map $F$ with the
complex-valued function $F=q_1+iq_2$, which gives
\[ F(z_1,z_2) = \bar{z_1}z_2. \]
The singular leaf $\Lambda_0=F^{-1}(0)$ is the union of the coordinate
planes $P_s=\{z_1=0\}$ and $P_u=\{z_2=0\}$, which are respectively the
stable and unstable manifolds for the $\ham{1}$-flow.
\begin{figure}[hbtp]
  \begin{center}
    \leavevmode
    \input{linearise.pstex_t}
    \macaption{the linearized flow}
    \label{fig:linearise}
  \end{center}
\end{figure}

\subsection{Geometry of the singular Lagrangean} 
\label{sec:geometrie}

Singular points of \emph{focus-focus} type are isolated, but it may
occur that the singular Lagrangean $\Lambda_0$ contains several of
them. Since $\Lambda_0$ is assumed to be compact, there will only be a
finite number of them. However, we will always assume that $\Lambda_0$
carries only \emph{one} critical point, which is certainly a ``generic''
assumption (in \cite{zung2}, it is claimed that one can always reduce
the number of \ff\ points by a small perturbation of the momentum
map). This point will always be denoted by $m$.

Recall that the usual theory of Arnold-Liouville gives action-angle
coordinates in a neighborhood of any connected compact regular
leaf. Here, we have compactness, but without smoothness. Another
viewpoint would be to consider the non-compact punctured leaf
$\Lambda_0\setminus\{m\}$; it is smooth, and invariant under the
action of the joint Hamiltonian flow of the system (this holds because
the critical point $\{m\}$ is a fixed point). But there is no hope for
the neighboring leaves to be diffeomorphic to
$\Lambda_0\setminus\{m\}$, since they are compact...

Nevertheless, the local structure of $\Lambda_0$ near the
singularity turns out to be sufficient data for its global
description.
\begin{prop}
  \label{prop:immersion}
  Let $\Lambda_0$ be a singular leaf of the momentum map, carrying a
  unique critical point $m$ of \emph{focus-focus} type.  Then the
  connected component of $m$ in $\Lambda_0$ is the image of a
  Lagrangean immersion of a 2-sphere with a double point. Deprived of
  $m$, it is an orbit of the Hamiltonian action of the system, with
  the structure of an affine infinite cylinder.
\end{prop}
\demo First of all,
note that we can suppose that $\Lambda_0$ is connected. Indeed, this
won't lead to any loss of information, since any other connected
components are necessarily regular tori, for which the usual theory
applies.

On $\Lambda_0\setminus\{m\}$, the action is locally free, hence by
standard arguments each connected component of
$\Lambda_0\setminus\{m\}$ (there are at most two of them) is an orbit,
on which all isotropy subgroups $I$ are conjugated. Therefore, this
orbit is diffeomorphic to $\RM^2/I$, which can be either $\RM^2$, or
the infinite cylinder $\RM^2/\ZM \simeq S^1\times\RM$ (it cannot be a torus
because it is not compact).

The local structure of the flow makes up for this ambiguity. Indeed,
according to Eliasson's theorem \cite{eliasson-these}, the action is
symplectically linearizable near $m$~: there exists a symplectic chart
in which the Hamiltonian vector fields of $f_1$ and $f_2$ are linear
combinations of the standard \emph{focus-focus} vector fields;
moreover, the involved coefficients form an invertible $2\times 2$
matrix $M_c$ which is locally \emph{constant} along each fiber
$\Lambda_c$. Now, the \emph{focus-focus} vector fields are by
assumption associated to the Hamiltonians $q_1$ and $q_2$. Since $q_2$
has periodic orbits in any neighborhood of $m$, the isotropy subgroup
$I$ is necessarily isomorphic to $\ZM$. The connected components of
$\Lambda_0\setminus\{m\}$ are therefore infinite cylinders.

On such a cylinder are \emph{globally} defined two infinitesimal
generators of the action, $\ham{1}$ and $\ham{2}$, by the requirement
that they are obtained from the initial generators $\ham{f_1}$ and
$\ham{f_2}$ by applying the matrix $M_c^{-1}$. In Eliasson's
coordinates, these vector fields are of course the standard
\emph{focus-focus} generators. The first one $\ham{1}$ is always
transversal to $\ham{2}$, and hence describes an ``axis'' of the
cylinder. $\ham{1}$ and $\ham{2}$ define the affine structure of the
cylinder.

Let us show now that $\Lambda_0\setminus\{m\}$ is necessarily 
connected.  From the connectedness of $\Lambda_0$ and the local 
structure at the critical point $m$, we \emph{a priori} know that 
$\Lambda_0\setminus\{m\}$ has at most two connected components, which 
correspond respectively to the stable and unstable manifolds of 
$\ham{1}$.  Let $x$ be a point on the unstable manifold, contained in 
a small $S^1$-invariant neighborhood $U$ of $m$.  We let the flow of 
$\ham{1}$ act on $x$ so that when time increases, the image $x(t)$ of 
$x$ goes out of $U$.  Now, we know from the dynamic on the infinite 
cylinder that $x(t)$ must go out of any compact subset in finite time.  
Since the manifold with boundary $\Lambda_0\setminus U$ is compact, 
the image of $x$ must necessarily enter $U$ again.  But it is clear 
from the local structure of the flow that $x(t)$ can approach $m$ only 
via the \emph{stable} manifold.  Therefore, the stable and the 
unstable manifold are connected to each other, and hence they are 
equal.

Thus, $\Lambda_0$ is homeomorphic to a cylinder whose extremities are
compactified in a unique point $m$, or a ``pinched torus'', or a
sphere with two points identified.
\begin{figure}[hbtp]
  \begin{center}
    \leavevmode
    \input{tore-pince2.pstex_t}
    \macaption{topology of $\Lambda_0$}
    \label{fig:tore-pince}
  \end{center}
\end{figure}
\begin{rema}
  We can readily deduce that $\pi_1(\Lambda_0)$ is isomorphic to $\ZM$,
  the integer coefficient counting the winding number along the ``big
  axis of the pinched torus''.
\end{rema}

From the differentiable viewpoint, we see in the local form that
$\Lambda_0\cap U$ is the union of two smooth submanifolds
diffeomorphic to open discs $D^2$ and transversally intersecting at
$m$. Therefore $\Lambda_0$ is the immersion of a smooth 2-manifold
with a unique double point. Since this manifold can be obtained by
gluing a $D^2$ at each extremity of the cylinder
$\Lambda_0\setminus U$, it is a sphere. \cqfd

\subsection{Monodromy of the Lagrangean fibration around 
  $\Lambda_0$}
\label{sec:fibration}
The first description of the monodromy invariant for singular
fibrations was carried out in \cite{duistermaat}, where the example of
the spherical pendulum was studied. Later, Maorong Zou in \cite{zou}
and Nguy{\^e}n Ti{\^e}n Zung in \cite{zung2} observed that the
monodromy at stake was typical of a \emph{focus-focus} singularity.
In \cite{zung2}, the computation of the monodromy is based on a
claimed normal form for the foliation. We present here another proof
of the result, relying on the local but Hamiltonian normal form of
Eliasson. The constructions explained in the proof will be used for
the quantization conditions (sections \ref{sec:regularisation} and
\ref{sec:global}).

At the same time as I was writing this article, Duistermaat and
Cushman were able to generalize the monodromy for non-Hamiltonian
systems \cite{cushman-duist2}.

Let us first briefly recall the definition, from a topological
viewpoint. For any \cis\ in dimension $2n$, a regular value
$c$ of the momentum map $F$ gives rise to a Liouville torus $\Lambda_c$,
%and the stabilizer of the action $\phy$ on
%$\Lambda_c$ is a $n$-dimensional lattice $I_c$ of $\RM^n$, called the
%\emph{period lattice}. $I_c$ is naturally isomorphic to
%$H_1(\Lambda_c,\ZM)$ by assigning to the n-uple $(t_1,\dots,t_n)$ the
%homology class of the loop $[0,1]\ni s\mapsto
%\phy(s.(t_1,\dots,t_n))$.  
%action coord -> basis -> extend to a neigh -> connection -> monodromy
and thus to a lattice $L_c=H_1(\Lambda_c,\ZM)$. Because $F$ is a
smooth fibration near $F^{-1}(c)$, the union of all $L_c$ is a locally
trivial smooth bundle $L$ over the open set $D_r$ of regular values of
$F$. Fixing a point $c_0$ and a smooth basis of $H_1(\Lambda_c)\simeq
\ZM^n$, a local section of $L$ near $c_0$ takes values in $\ZM^n$, and
thus is locally constant. We therefore have the notion of parallel
transport above a path in $D_r$. The holonomy associated with it is a
homomorphism from $\pi_1(D_r,c_0)$ to $\mathrm{Aut}(L_{c_0})\simeq
\mathrm{GL}(n,\ZM)$, called the \emph{monodromy} of $L$.

Back to our focus-focus system, we take $D_r$ to be a small ball
around the origin in $\RM^2\simeq \CM$, minus the origin. Because its
$\pi_1$ is $\ZM$, the problem is reduced to determining the monodromy
of a simple loop around $0$.

Recall from last paragraph that for any focus-focus system $(p_1,p_2)$,
there is a local normal form $(q_1,q_2)=f(p_1,p_2)$ near $m$, where
$f$ is a local diffeomorphism of $(\RM^2,0)$, and $q_i$'s are defined in
symplectic coordinates around $m$ by paragraph \ref{sec:linear}. This
yields two particular vector fields $\ham{1}$ and $\ham{2}$ defined on
a full neighborhood of $\Lambda_0$ by~:
\[ (\ham{1},\ham{2}) = df.(\ham{p_1},\ham{p_2}). \]
The flow of $\ham{2}$ on $\Lambda_c$ is, as we already saw, periodic
of period $2\pi$, so the orbit of any point $x$ is a simple loop
$\Gamma_x$ on $\Lambda_c$ depending smoothly on $x$.  Its homology
class $\gamma_2(c_0)$ is therefore an invariant of the monodromy
operator. Let us denote by $\gamma_1(c_0)$ the homology class of a
path starting on $\Gamma_x$ in the direction of $\ham{1}$, going back to
$\Gamma_x$ and closing up on $\Gamma_x$ in the direction of $\ham{2}$.

$(\gamma_1(c_0),\gamma_2(c_0))$ is a basis of $L_{c_0}$ depending
locally smoothly on $c_0$; the monodromy of the system is expressed on
this basis by the following proposition~:
\begin{prop}[\cite{zung2}]
  \label{prop:monodromie}
  Let $(p_1,p_2)$ be a \cis\ with a singular leaf carrying a unique
  singular point of \emph{focus-focus} type. Then for any $c_0$ close
  enough the the critical value $0\in\RM^2$, the monodromy matrix,
  expressed in the basis $(\gamma_1(c_0),\gamma_2(c_0))$ of $L_{c_0}$,
  is equal to
  \[ \left(\begin{array}{cc} 1 & 0 \\
      \epsilon & 1\end{array}\right).\] 
  Here $\epsilon$ is the sign of
  $\det M$, where $M\in GL(2,\RM)$ is the unique matrix such that~:
  \[ (\mathcal{H}(q_1),\mathcal{H}(q_2)) =
  M.(\mathcal{H}(p_1),\mathcal{H}(p_2)). \]
\end{prop}
\demo The idea of the proof, already apparent in the construction of
the basis $(\gamma_1(c_0),\gamma_2(c_0))$, is to show that all
contributions to the monodromy can be concentrated in a small
neighborhood of the critical point $m$.

Let $U$ be such a small ball around $m$, taken to be
$S^1$-invariant. The intersection of $\partial U$ with $\Lambda_0$ is
composed of two circles $\Gamma_u$ and $\Gamma_s$ respectively drawn
on the local unstable and stable manifolds at $m$. Let us take a point
$x$ on $\Gamma_u$, and let $\Sigma$ be a small section transversal to
$\Lambda_0$ at $x$. The tangent space to $\Lambda_0$ at $x$ is exactly
the kernel of $dF$. Thus, the momentum map $F$ realizes a local
diffeomorphism from $(\Sigma,x)$ onto $(\RM^2,0)$. We can then assume
that $\Sigma$ is a preimage of a small open disc $D\subset D_r$ around
$0$.

Let $\Omega$ be the invariant open subset of $M$ obtained as the union
of orbits of points in $\Sigma$. It is a neighborhood of $\Lambda_0$,
composed of the union of $\Lambda_0$ and invariant tori. Because of
the compactness of $\overline{\Omega\setminus U}$, any orbit of some
point $y\in\Sigma$ close to $x$ intersect $\partial U$ in two circles
close to $\Gamma_u$ and $\Gamma_s$. The restriction of $F$ to
$\Omega\setminus U$ is therefore a trivial fibration over $D$ whose
fibers are finite length $S^1$-invariant affine cylinders, where the
affine structure that defines the ``axis'' of each cylinder is given
by the vector field $\ham{1}$~:
\begin{equation}
  \label{equ:triv-fibration}
  \begin{array}{c} \Omega\setminus U \simeq D\times([0,1]\times
  S^1)\\ \flechebas{F} \\ D
\end{array}
\end{equation}
For instance, the circles $\Gamma_u$ and $\Gamma_s$ are represented in
this trivialization by $(0,(0,S^1))$ and $(0,(1,S^1))$ respectively
(see figure \ref{fig:fibration}). 
\begin{figure}[hbtp]
  \begin{center}
    \leavevmode
     \input{fibration.pstex_t}
    \macaption{The fibration $F$ on $\overline{\Omega\setminus U}$}
    \label{fig:fibration}
  \end{center}
\end{figure}
For any $y\in\Omega\setminus U$ corresponding to the point
$(c,(t,\tau))$ in the above trivialization, we denote by $[y_u,y_s]$
the $\ham{1}$-path corresponding to the axis $(c,([0,1],\tau))$. As we
can assume that $U$ lies inside an open neighborhood of $m$ where
Eliasson's normal form applies, $y_u$ and $y_s$ have complex
coordinates $(y_{u/s,1},y_{u/s,2})$ satisfying
\[ \overline{y_{u/s,1}}y_{u/s,2} = f(c). \] 
Note that the matrix $M$ in the statement of the proposition is equal
to $df(m)$.
\begin{rema}
  We don't actually need the
  full statement of Eliasson's normal form, because higher order terms
  in the normal form wouldn't be of any trouble for the local
  computation of the monodromy (see \cite{cushman-duist2}).  
\end{rema}

If $c\neq 0$, using formula
(\ref{equ:flot}) for the flow, $y_s$ can be joined to $y_u$ by letting the
flow of $q_1$ act during a time $t(y)=\ln(|y_{u,1}|/|y_{s,1}|)$ and
that of $q_2$ during a time \[s(y)= \arg(y_{u,1}/y_{s,1}).\] By this
construction, we get a locally smooth family of loops on $\Lambda_c$,
whose homology class at $c_0$, $\gamma_1(c_0)$, together with
$\gamma_2(c_0)$, form the basis referred to in the statement of the
proposition.

Let $[x_u,x_s]$ be the axis through $x$ (by definition, $x_u=x$ is on
the local unstable manifold $P_u$ and $x_s$ on the local stable
manifold $P_s$). If $y$ tends to $x$, then $y_u$ and $y_s$ tend to
$x_u$ and $x_s$ respectively. Therefore, since $x_{u,2}=x_{s,1}=0$,
$s(y)\sim \arg(x_{u,1}\overline{x_{s,2}}) + \arg(f(c))$ as $y$ tends
towards $x$.

Now let $\mathcal{C}$ be a small loop around $x$ in $\Sigma$ such that
its image by $F$ is a circle in $D$ through $c_0$, oriented in the
trigonometric sense~:
%For each
%$y\in \mathcal{C}$ with $F(y)=c$ we've just provided a construction of
%a simple loop $\gamma_1(y)$ on $\Lambda_c$ which, together with an
%$S_1$ orbit, form a basis of $\pi_1(\Lambda_c)$. 
letting $y$ run around $\mathcal{C}$, and with it, the loop
$\gamma_2(c)$, amounts to increasing $\arg(c)$. Therefore, after such
a revolution, $s(y)$ get increased or decreased by $2\pi$, depending
on the diffeomorphism $f$ being respectively orientation preserving or
reversing. This means that the loop obtained is equal to the initial
one composed by respectively \emph{plus} or \emph{minus} an $S^1$
orbit, which proves the proposition. \cqfd

\subsection{Microlocal focus-focus}
\label{sec:microloc}
We are now given a self-adjoint semi-classical \cis\  $P_1(h),P_2(h)$
whose associated classical system presents a focus-focus singularity
at a point $m$, as defined in the preceding section.

The aim of this section is the complete analysis of microlocal
solutions around the critical point $m$ of the system~:
\begin{equation}
  \label{equ:systeme}
  \left\{ 
  \begin{array}{l}
P_1^E(h) u_h^E \sim 0 \\
P_2^E(h) u_h^E \sim 0
  \end{array}\right.
\end{equation}
in which $P_j^E$ is a regular deformation of $P_j$ with a
parameter $E$ varying in a compact neighborhood of $0\in\RM^2$ (see
definition \ref{defi:parameter}).

The problem of globalizing these solution to a neighborhood of the
whole critical leaf $\Lambda_0$ will be dealt with in the next
section.

For now, the crucial point is the following~:
\begin{prop}
  \label{prop:dimension}
  The space of microlocal solution of (\ref{equ:systeme}) at
  $m$ has dimension $\leq 1$. In other words, if $u_h^E$,
  $(h,E)\in\Gamma$ is a non-trivial solution of  (\ref{equ:systeme}) at
  $m$, then for any other such solution $v_h^E$, there exists a
  constant $C^E(h)\in\CM_h$, depending smoothly on $E$, such that~:
  \[ v_h^E \sim C(h)^Eu_h^E \textrm{ at } m.\]
\end{prop}
To prove this, we use a parameter dependent version of the normal form of
\cite{san}~: let $Q_j(h)$ be the standard \emph{focus-focus} operators
on $\RM^2$, that is, written in polar coordinates $(r,\theta)$~:
\[ Q_1=\frac{h}{i}(r\deriv{}{r}+1) \quad \textrm{and} \quad
  Q_2=\frac{h}{i}\deriv{}{\theta}. \]
For notation convenience, we will write $\mbold{P}(h)$ for the vector
$(P_1(h),P_2(h))$, and similarly for other quantities.
% We have 
% \[ \mbold{P}^E = f^E(\vec{P}),\] and $df^E(0) \in GL_2(\RM)$
% for small $E$. 
Let $\mbold{p}$ be the momentum map for the system
$(p_1,p_2)$. The target space $\RM^2$ of $\mbold{p}$ will be identified 
with $\CM$ by $\mbold{p}=p_1+ip_2$. 
As we already saw, in suitable symplectic coordinates
near $m$, there is a local diffeomorphism $F$ of $\RM^2$ such that
$\mbold{p} = F(\mbold{q})$. Thus,
\[\mbold{p}^E=f^E\circ F(\mbold{q}).\] 
The semi-classical normal form that we get is the following~:
\begin{lemm}[\cite{san}]
  \label{lem:FNQ}
  There exists a unitary \fio\  $U^E(h)$, a microlocally invertible
  $2\times 2$ matrix $\mathcal{M}^E(h)$ of \pdo s, and 
  $\bep^E(h)\in\CM_h$, everything depending smoothly on $E$, such
  that, microlocally on a neighborhood of $m$~:
  \[ (U^E)^{-1}\mbold{P}^EU^E =
  \mathcal{M}^E.(\mbold{Q}-\bep^E(h)) + O(h^\infty). \]
  $\bep^E$ admits a semi-classical \das\ of the
  form~:
  \[ \bep^E(h) \sim \sum_{k\geq 0} h^k\bep_k^E. \]
  The first two terms satisfy~:
  \[ \bep_0^E = F^{-1}\circ (f^E)^{-1}(\vec 0), \]
  \[ \bep_1^0 = -M^{-1}.(\mbold{r}(m)). \]
  $\bullet$ $\mbold{r}(m)$ is the value at $m$ of the sub-principal
  symbol of $\mbold{P}$.\\ $\bullet$ $M$ is the value at $0$ of
  principal symbol of $\mathcal{M}^0$. It is also defined by
  $M=dF(0)$, or, equivalently~:
  \[d^2p^0(m) = M.d^2q(0). \]
\end{lemm}
\begin{rema}
  We will show later (paragraph \ref{sec:recover}) that the expansion 
  $\epsilon_2^E(h)=\im(\bep^E(h))$ is, up to $O(h^\infty)$, 
  independent on the choice of a \fio\ $U(h)$.  For the other 
  ``semi-classical invariant'' $\epsilon_1^E(h)=\re(\bep^E(h))$, this 
  might not hold; however, will will prove thanks to corollary 
  \ref{coro:global} that its formal infinite germ at $E=0$ is indeed 
  invariant by a change of \fio\ $U(h)$.
\end{rema}
\begin{rema}
  If we apply this lemma to $E=ht$ for a bounded $n$-dimensional
  parameter $t$, we can apply  remark \ref{rem:dependence} 
  and express everything in terms of $t$. We
  obtain an asymptotic expansion $\tilde{\bep}^t(h)\egdef
  \bep^{ht}(h)$ satisfying~:
  \[ \tilde{\bep}^t_0 = 0, \quad \tilde{\bep}^t_1 =
  -M^{-1}.\left(\mbold{r}^0(m) + \deriv{f^E}{E}_{\restr
      E=0}(0).t\right). \]
\end{rema}

Thanks to this lemma, the solutions of (\ref{equ:systeme}) get transformed
into those of simple differential operators on $\RM^2$ for which we
have explicit solutions. We can now turn to the proof of proposition
\ref{prop:dimension}. Assume that we know an exact solution of the system
$u^E_{\textrm{\tiny exact}}\in\mathcal{S}'(\RM^2)$ (this will always
occur as soon as the system admits a microlocal solution -- see
proposition \ref{prop:exacte}). Let $u_h^E$ be any other microlocal
solution near $0$.

Let $\Upsilon_R$ be the ``vertical strip'' $\{r<R\}$ in the cotangent
bundle $T^*\RM^2=((r,\theta),(\rho,\alpha))$. The Lagrangean leaves
of the system in $\Upsilon_R$, deprived of a small ball $B$ around
$0$, are non-singular, hence from proposition \ref{prop:dim1reg} there
exists a constant $C^E(h)$ such that microlocally on
$\Upsilon_R\setminus B$, 
\[ u_h^E \sim  C^E(h) u^E_{\textrm{\tiny exact}} \]
(remember that such a formula involves uniformity with respect to
$E$).

As in \cite{colin-p}, we shall use this exact solution to construct a
\emph{local} (mod $O(h^\infty)$) solution, microlocally equal to
$u_h^E$ near the origin. Let us define $v_h^E$ by~:
\[ \fourier v_h^E = \chi \fourier u_h^E + (1-\chi)\fourier C^E
u^E_{\textrm{\tiny exact}}, \] where $\chi$ is a compactly supported
$\Cinf$ function with value 1 in a neighborhood of $0$.
Then $v_h^E$  and $u_h^E$ are microlocally equal near the
origin. Furthermore, $v_h^E$ is a microlocal solution of the system at
each point of $\Upsilon_R$.

Now there exists differential operators $\hat{P}_j^E(h)$ such that for
any $w\in\mathcal{S}'(\RM^2)$, $\fourier(P_j^Ew) = \hat{P}_j^E
\fourier w$.  Differential operators preserve the support of
distributions; therefore, since $\fourier v_h^E = \fourier
u^E_{\textrm{\tiny exact}}$ far from the origin, we have there
$\hat{P}_j^E \fourier v_h^E = \fourier(P_j^Eu^E_{\textrm{\tiny
exact}}) = 0$. So we have
\[ P_j^E v_h^E \sim 0 \] near each point of $\Upsilon_R$, and
\[ \fourier (P_j^E v_h^E) = 0 \]
far from 0. We deduce from this that $P_j^E v_h^E$ satisfies the
hypothesis of lemma \ref{lem:global} on $\Upsilon_R$, that is, 
\[ P_j^E v_h^E = O(h^\infty) \textrm{ on }
\Upsilon_R\cap\{\rho=0\}. \]

We have thus at our disposal a solution $v_h^E$ of the system 
$P_j^E v_h^E = w_{jh}^E$ on $\RM^2$ with $w_{jh}^E=O(h^\infty)$ for
$r<R$. On $\RM^2\setminus\{0\}$ this is a regular linear differential
system of the first order, whose solutions are of the form
\begin{equation}
  \label{equ:particular}
  v_h^E = Cste^E(h).u^E_{\textrm{\tiny exact}}+ u^E_{\textrm{\tiny
particular}}; 
\end{equation} 
$u^E_{\textrm{\tiny exact}}$ is a solution of the
homogeneous system, and $u^E_{\textrm{\tiny particular}}$ is a
particular solution of the system that can be computed from
$u^E_{\textrm{\tiny exact}}$ by the method of ``variation of the
constant''. Explicitly, if we write 
\[  v^E_h = g^E.u^E_{\textrm{\tiny exact}}, \] we find that the
function $g^E$ must satisfy the following system, expressed in polar
coordinates~:
  \[ \left\{ \begin{array}{l} \deriv{g^E}{r} =
      \frac{i}{h}\frac{w^E_{1h}}{r u^E_{\textrm{\tiny exact}}} =:
      a_{1h}^E \\
      \vspace{1mm}\\
      \deriv{g^E}{\theta} =
      \frac{i}{h}\frac{w^E_{2h}}{u^E_{\textrm{\tiny exact}}} =: a_{2h}^E.
      \end{array} \right. \]
If $(r_0,\theta_0)$ is any point with $0<r_0<R$, $g^E$ is integrated
by the formula~:
\[ g^E(r,\theta) = \int_0^1 \mbold{a}^E_h (r_0+t(r-r_0),
\theta_0+t(\theta-\theta_0)) .(r,\theta) dt + g^E(r_0,\theta_0). \] It
is easy to check from the explicit formula of $u^E_{\textrm{\tiny
exact}}$ (see proposition \ref{prop:exacte}) that $\mbold{a}^E_h$ is
uniformly $O(h^\infty)$ in the disc $r<R$; so with $Cste^E(h) =
g^E(r_0,\theta_0)$ we get in (\ref{equ:particular}) that
$u^E_{\textrm{\tiny particular}}$ is uniformly $O(h^\infty)$ near the
origin.

The distribution $z_h^E = v_h^E - Cste^E. u^E_{\textrm{\tiny exact}}
    - u^E_{\textrm{\tiny particular}}$ is a solution of $P_j^E
    z_h^E = w_{jh}^E$ with support concentrated at the
    origin. Therefore, it is a finite sum of derivatives of the Dirac
    distribution. Now, the identity $P_2^E z_h^E = O(h^\infty)$
    readily shows that the the coefficients involved in $z_h^E$ are
    uniformly $O(h^\infty)$, and so $z_h^E$ is.

Finally we are left with $v_h^E = Cste^E. u^E_{\textrm{\tiny exact}}
+ O(h^\infty)$ near $0$, which implies $v_h^E \sim
Cste^E. u^E_{\textrm{\tiny exact}}$ near $0$. This, provided we
prove proposition \ref{prop:exacte}, achieves the proof of proposition
\ref{prop:dimension}. \cqfd

\begin{rema}
  The proof shows that the constants $C^E(h)$ and $Cste^E(h)$ are
  microlocally equal. Besides from the preceding arguments, this
  essentially comes from the local connectedness of the Lagrangean
  leaves $\Lambda_0^E$ for $E\neq 0$.
\end{rema}

It remains now to check~:
\begin{prop}
  \label{prop:exacte}
  We consider the system  (\ref{equ:systeme}) on $\RM^2$ with
  $\mbold{P}^E(h)=\mbold{Q}_j-\bep^E(h)$, as previously
  defined. The following conditions are equivalent~:
\begin{enumerate}
\item there exists a non-trivial microlocal solution $u_h^E$,
  $(h,E)\in\Gamma$ of (\ref{equ:systeme}) on a neighborhood of the
  origin;
\item up to $O(h^\infty)$ modifications of the $P_j^E$, there
  exists an exact solution  $u^E_{\textrm{\tiny exact}}$ of
  (\ref{equ:systeme})  in $\RM^2\setminus\{0\}$ for $(h,E)\in\Gamma$;
\item $\frac{\epsilon_2^E(h)}{h}\in \ZM + O(h^\infty)$ uniformly for
  $(h,E)\in\Gamma$.
\end{enumerate}
\end{prop}
\demo Equivalence between 2. and 3. comes from the fact that at each
point of $\RM^2\setminus\{0\}$ the solutions of the system are spanned
by the function~:
\begin{equation}
  \label{equ:solution}
u^E_{\textrm{\tiny exact}} = \frac{1}{r} e^{
  i\frac{\epsilon_1^E(h)}{h}\ln r} e^{
  i\frac{\epsilon_2^E(h)}{h}\theta}.
\end{equation}
Therefore, there exists a solution around the origin if and only if
$\frac{\epsilon_2^E}{h} \in \ZM$.

Since $u^E_{\textrm{\tiny exact}}$ as defined above is an admissible
functional in $\RM^2$ depending smoothly on $E$, and a solution of the
system on $\RM^2$, it is also a microlocal solution. Hence 2. implies
1.

Now let $u_h^E$ be a non-trivial microlocal solution on a
neighborhood $\Omega$ of $0$. From the uniqueness of microlocal
solutions at points different from the origin (proposition
\ref{prop:dim1reg}), we deduce that on every simply connected open subset
$\Omega_\alpha$ of $\Omega\setminus\{0\}$, there is a constant
$(C^E)^\alpha(h)$ such that
\[ u_h^E \sim (C^E)^\alpha u^E_{\textrm{\tiny exact}}. \]
Then, going around $0$, the $(C^E)^\alpha$'s must microlocally
coincide. Hence the $e^{ i\frac{\epsilon_2^E}{h}\theta}$ must all
microlocally coincide, uniformly for $(h,E)\in\Gamma$, which implies
3. \cqfd

Put together, the two propositions give the following result (``first
quantization condition'')~:
\begin{theo}
   \label{theo:quantif1}
   The system (\ref{equ:systeme}) has a non-trivial microlocal
   solution $u_h^E$ (with $(h,E)\in\Gamma$) on a neighborhood of $m$
   if and only if
   \[ \frac{\epsilon_2^E(h)}{h} \in \ZM + O(h^\infty), \quad
   (h,E)\in\Gamma. \]
   Any other solution is a multiple of it.
\end{theo}

\begin{rema} 
  \label{rem:quantif1}
  As we already mentioned, microlocal solutions should be viewed as
  half-densities on the base space. Here, we can rewrite the
  half-density $u^E_{\textrm{\tiny exact}}|dx\wedge dy|^{\frac{1}{2}}$
  as \[ e^{ i\frac{\epsilon_1^E}{h}\ln|r|} e^{
  i\frac{\epsilon_2^E}{h}\theta} \left|\frac{dx\wedge
  dy}{r^2}\right|^{\frac{1}{2}} , \] and for $r\neq 0$, its modulus
  $\left|\frac{dx\wedge dy}{r^2}\right|^{\frac{1}{2}}$ is exactly the
  pull down of the canonical invariant Liouville half-density
  $\rho(\lambda)$ on $\Lambda_0^E$. In this sense, $u^E_{\textrm{\tiny
  exact}}$ is a normalized classical oscillatory half-density on every
  open subset $\Lambda^E$ of $\Lambda_0^E$ that projects
  diffeomorphically onto $X$ (thus excluding, for $E=0$, the vertical
  $(\xi,\eta)$-space in $T^*X$).  Concerning the phase function
  $\phy_E=\epsilon_{1,0}^E \ln|r| + \epsilon_{2,0}^E \theta$, a handy
  tool for its expression is the use of complex coordinates $z_1=x+iy$
  and $z_2=\xi+i\eta$; we also write $\bep^E_0 = \epsilon_{1,0}^E +
  i\epsilon_{2,0}^E$, which leads to
  \[ \phy^E = \re(\bep^E_0\ln\bar{z_1}). \]
  (Note that $\bep^E_0\ln\bar{z_1}$ is not well-defined, but thanks to
  the assumption on $\epsilon^E_2(h)$, its real part is indeed
  well-defined modulo $2\pi h\ZM$.)  We readily deduce that the graph
  of $\phy^E$ in $T^*X$ is the set of $(z_1,z_2)$ such that
  $z_2=2\deriv{\phy^E}{\bar{z_1}} = \bep^E_0/\bar{z_1}$. Thus
  $\Lambda_0^E= \{ (z_1,z_2)\in T^*X,\quad \bar{z_1}z_2 = \bep^E_0
  \}$. Note that this also holds if $E=0$, for then $\bep^E_0=0$ and
  $\Lambda_0^0$ consists of the union of both coordinate planes
  $z_1=0$ and $z_2=0$.
\end{rema}

\subsection{Regularization of $[\kappa]$}
\label{sec:regularisation}
Before stating our principal result in the next section (theorem
\ref{theo:global}), we need to describe one more geometrical
ingredient, related to the divergence of the integral of the
sub-principal form $\kappa$ along paths through $m$.

We still let $(P_1(h), P_2(h))$ be a semi-classical \cis\  with a
singular principal Lagrangean carrying a unique singular point $m$ of
focus-focus type. Here again, we identify the target space of the
momentum map $p$ with the complex plane $\CM$. Thus $p=p_1+ip_2$.
For small $c\in\CM$, we denote by $\Lambda_c$ the
Lagrangean leaf $\Lambda_c=\{p^{-1}(c)\}$, and $\kappa_c$ is as usual
the sub-principal form on $\Lambda_c$.

Let $\gamma_c$ be a smooth family of loops on $\Lambda_c$ such that
$\gamma_0$ is a simple loop through $m$ starting on the local unstable
manifold of $(\Lambda_0,m)$, and coming back to $m$ via the local
stable manifold. We know that the non-triviality of the monodromy
implies that the map $c\fleche \gamma_c$ cannot be univalued, but
proposition \ref{prop:monodromie} gives a precise control over its
multivaluedness. The purpose of this paragraph is to control the
integral of the sub-principal form along such a family of paths, and in
particular to determine its divergence principal part as $c\fleche 0$.

Let $M$ be the $2\times 2$ real matrix such that
$\mathcal{H}(p)=M.\mathcal{H}(q)$, where $q=q_1+iq_2$ (for the
linearized flow, we always use the notations of \ref{sec:linear}). In
accordance with lemma \ref{lem:FNQ} we let
\[ \bep_1 = -M^{-1}.r(m), \]
where $r=r_1+ir_2$ is the joint sub-principal symbol of the system.
\begin{prop}
  \label{prop:regularisation}
  For any multivalued smooth family $\gamma_c$ of loops on
  $\Lambda_c$ such that $\gamma_0$ is a simple loop through $m$
  starting on the local unstable manifold of $(\Lambda_0,m)$, and
  coming back to $m$ via the local stable manifold, then 
  \begin{itemize}
  \item the limit~:
    \[ I_{\gamma_0}(\kappa) = \lim_{\dessus{c\fleche
        0}{\mathit{along\,a\,ray}}} 
    \left( \int_{\gamma_c} \kappa_c + \re(\bep_1\ln(M^{-1}c))\right) \]
    exists,
  \item and its class modulo $2\pi\im(\bep_1)\ZM$ actually does
    not depend on $\gamma_c$, provided $\gamma_0$ satisfies the
    hypothesis;
  \item finally, this class is also given by the formula~:
    \begin{equation}
      \label{equ:regularisation}
      I_{\gamma_0}(\kappa) = \lim_{(s,t)\fleche (0,0)} \left(
        \int_{A_0=\gamma_0(s)}^{B_0=\gamma_0(1-t)} \!\!\!\kappa
        +\re(\bep_1)ln(r_{A_0}\rho_{B_0}) +
        \im(\bep_1)(\theta_{A_0}-\alpha_{B_0})\right),
    \end{equation}
    where $(r,\theta,\rho,\alpha)$ are the polar symplectic
    coordinates near the origin defined in paragraph \ref{sec:linear}.
  \end{itemize}
  $I_{\gamma_0}(\kappa)$ will be called the
  \emph{principal value} of $\int_{\gamma_0}\kappa$.
\end{prop}
\demo Let $([\gamma_1(c)],[\gamma_2(c)])$ be the smooth basis of
$\pi_1(\Lambda_c)$ defined in proposition \ref{prop:monodromie}.
The homotopy class of $\gamma_c$ decomposes into
\[ [\gamma_c] = n_1(c)[\gamma_1(c)]+ n_2(c)[\gamma_2(c)], \]
for integers $n_1(c)$ and $n_2(c)$. On $\Lambda_0$, $\gamma_2(0)$
degenerates into a trivial loop, and by hypothesis we have
\[  [\gamma_0] = [\gamma_1(0)]. \]
But we know from proposition
\ref{prop:monodromie} that the group bundle
$\pi_1(\Lambda_c)/\gener{\gamma_2(c)}\fleche c$ is trivial near $c=0$,
with fibers isomorphic to $\ZM$. Therefore, $n_1(c)$ must be
identically equal to $1$ for small $c$.

If $c$ is constrained into a sector $\Delta$ of the plane, then the
family $\gamma_c$ must be continuous at zero, which implies that
$n_2(c)$ is constant for small $c$. Since the integral of $\kappa_c$
along the $\ham{2}$-circle $\gamma_2(c)$ is equal to $-2\pi r_2(m) +
O(c)$, we have
\[  \lim_{\dessus{c\fleche 0}{c\in\Delta}}
\int_{n_2(c)[\gamma_2(c)]}\kappa_c = 2\pi\im(\bep_1). \]

This, provided we prove the first point of the proposition, shows the
second point. It remains to show the convergence of
$\int_{\gamma_1(c)}\kappa_c$, and the formula
\ref{equ:regularisation}. We shall still assume that $c$ stays in a
sector $\Delta$ and by $\gamma_c$ we denote a smooth representative of
the first basis element $[\gamma_1(c)]$.

Let $f$ be the function on $(\RM^2,0)$ such that, in some
neighborhood $\Omega$ of $m$ and some local symplectic coordinates
around $m$, $p=f(q)$, so that $M=df(0)$. If we let $\tilde{c}=f(c)$,
the problem amounts to investigating the limit of
\[ \int_{\gamma_{\tilde{c}}}
  \kappa_{\tilde{c}} + \re(\bep_1\ln(c)), \]
because $\ln(M^{-1}f(c)) = \ln(c) +O(c)$.
Now, $\kappa_{\tilde{c}}$ is defined by 
\[ \kappa_{\tilde{c}}.\ham{p_j} = - r_j,\]
which also reads, in $\Omega$~:
\[ (\kappa_{\tilde{c}}.\ham{q_1},\kappa_{\tilde{c}}.\ham{q_2}) = -
(df(q))^{-1}.(r_1,r_2). \] This means that $\kappa_{\tilde{c}}$ is, in
$\Omega$, the sub-principal form on the $c$-level set of a system with
principal symbol $q$ and sub-principal symbol $(df(q))^{-1}.r$ (which
we call $r$ again). So let us work with these data and forget the
tilde. 

As we shall see, we can even reduce the problem to a
\emph{constant} sub-principal symbol $r(0) (= -\bep_1)$. This is due
to the following ``critical Poincar\'e lemma''(see \cite{san})~: there
exists a function $g$ on a neighborhood of $0\in T^*\RM^2$ and smooth
functions $h_j$ on $(\RM^2,0)$ such that
\[ r_j = dg.\ham{q_j} + h_j(q).\] 
This implies that on $\Omega$, for non-zero $c$, the function
\[ K_c \equiv  -g  -h_1(c)\ln(r) -h_2(c)\theta \quad (\mod
2\pi h_2(c) ) \] is a primitive of $\kappa_c$ on
$\Lambda_c=\{r\rho e^{i(\alpha-\theta)}=c\}$.  Suppose now that $A_0$
and $B_0$ are points on $\gamma_0$ lying respectively on the local
unstable and stable manifolds of $\Lambda_0$. Then there exist smooth
families of points $A_c$ and $B_c$ on $\Omega\cap\gamma_c([0,1])$ such
that $\lim_{c\fleche 0}(A_c,B_c) = (A_0,B_0)$. We can write~:
\[ \int_{\gamma_c} \kappa_c =  \int_{A_c}^{B_c}\kappa_c +
K_c(A_c)-K_c(B_c) \]
\[ =  \int_{A_c}^{B_c}\kappa_c + g(B_c)-g(A_c)
-h_1(c)ln(r_{A_c}/r_{B_c})  -h_2(c)(\theta_{A_c}-\theta_{B_c}), \]
which can also be written 
\[ \int_{A_c}^{B_c}\kappa_c + g(B_c)-g(A_c)
-h_1(c)ln(r_{A_c}\rho_{B_c}) -h_2(c)(\theta_{A_c}-\alpha_{B_c}) \]
\[ + \left(h_1(c)ln(|c|) -
  h_2(c)\arg(c)\right). \]
Everything but the last parenthesis has a smooth behavior as $c\fleche
0$ in $\Delta$, and that last parenthesis is nothing else than
$\re((h_1(c)+ih_2(c))\ln c)$. So there is a smooth function
$I(A_c,B_c)$ such that~:
\[ I(A_c,B_c) =  \int_{\gamma_c} \kappa_c -
\re((h_1(c)+ih_2(c))\ln c) .\]
\begin{rema}
  Because of the monodromy matrix $\left(
    \begin{array}{cc}
      1 & 0 \\ 1 & 1
    \end{array}\right)$ (proposition \ref{prop:monodromie}), as $c$
  loops in the positive sense around the origin, $\int_{\gamma_c}
  \kappa_c$ gets increased by
  \[ \int_{\gamma_2(c)} \kappa_c = -2\pi h_2(c) \]
  whereas $\re((h_1(c)+ih_2(c))\ln c) = h_1(c)ln(|c|) - h_2(c)\arg(c)$
  gets increased by the same quantity. This shows that $I(A_c,B_c)$ is
  actually a \emph{univalued} smooth function around $c=0$.
\end{rema}
\begin{figure}[hbtp]
  \begin{center}
    \leavevmode
    \input{kappa.pstex_t}
    \macaption{regularization of $\kappa$}
    \label{fig:kappa}
  \end{center}
\end{figure}
Because  $(h_1(c)+ih_2(c))= -\bep_1 + O(c)$, we deduce that if we fix
a determination of $\ln c$ and make $c$ tend to $0$
in $\Delta$, then
\[\re\left((h_1(c)+ih_2(c))\ln c +
\bep_1\ln c\right) \fleche 0\]
 Therefore,
$I_{\gamma_0}(\kappa)$ exists and is equal to $I(A_0,B_0)$. Note that
it is independent of $A_0$ and $B_0$, and hence equal to the limit of
$I(A_0,B_0)$ as $A_0$ and $B_0$ tend to the critical point $m$. That
means that the function $g$, being smooth at $m$, plays no role,
justifying our assertion that the sub-principal symbols can be taken
to be constant near $m$. This remark also yields the useful expression
of $I_{\gamma_0}(\kappa)$ involving only objects living on
$\Lambda_0$ given by formula (\ref{equ:regularisation}). \cqfd

%Let us turn now to the last assertion of the proposition.  For each
%$c$, the integral of $\kappa_c$ along a $\ham{2}$-circle $\delta_2(c)$
%is equal to $-2\pi h_2(c)$. So, modulo $2\pi h_2(c)\ZM$, the
%map $\gamma_c\fleche \int_{\gamma_c}\kappa_c$ is well-defined on
%$\pi_1(\Lambda_c)/\gener{\delta_2(c)}$. We know from proposition
%\ref{prop:monodromie} that the group bundle
%$\pi_1(\Lambda_c)/\gener{\delta_2(c)}\fleche c$ is trivial near $c=0$,
%with fibers isomorphic to $\ZM$. Therefore, the homotopy class in
%$\pi_1(\Lambda_c)/\gener{\delta_2(c)}$ of any smooth family $\gamma_c$
%is, for small $c$, entirely determined by that of $\gamma_0$ on
%$\Lambda_0$. This shows that the map $\gamma_0\fleche
%I_{\gamma_0}(\kappa)$ is well defined and homotopically invariant.

%Because $\pi_1(\Lambda_0)=\ZM$ and $\gamma_0$ is supposed to be a
%simple, oriented, non trivial loop, this map has only one possible
%value. This proves last point of the proposition.

\subsection{Global solutions}
\label{sec:global}

Here is finally the main result of our work, namely the condition
under which there exists a microlocal solution of the system
(\ref{equ:systeme}) on a neighborhood of the whole critical Lagrangean
$\Lambda_0$. The hypothesis are the same as in paragraph
\ref{sec:microloc}. From now on, $n$ is a generic integer, and not the 
dimension of the ambient symplectic manifold, which is $4$. 

\begin{theo}
  \label{theo:global}
  Let $\epsilon_j^E=\epsilon_j^E(h)$ be the constants defined by lemma
  \ref{lem:FNQ}. \\ $\bullet$ There exists a $\T$-valued constant
  $(\lambda^E)^{\out}=(\lambda^E)^{\out}(h)$
  admitting an asymptotic expansion in $\geq -1$ powers of $h$ whose
  terms are smooth functions of $E$, such that the system
  (\ref{equ:systeme}) admits a microlocal solution $u_h^E$,
  $(h,E)\in\Gamma$, on a full neighborhood of $\Lambda_0$ if and only
  if the following two conditions are satisfied~:
  \begin{enumerate}
  \item ``first quantization condition''~:
    \[  \frac{\epsilon_2^E(h)}{h} \in \ZM + O(h^\infty);\]
  \item ``second quantization condition''~:
    \[ (\lambda^E)^{\out} +
    n\frac{\pi}{2} -
    \frac{\epsilon_1^E}{h}\ln(2h) - 2\arg
    \Gamma\left(\frac{i\epsilon_1^E/h+1+n}{2}\right) \in
    2\pi\ZM + O(h^\infty), \]
    where $\frac{\epsilon_2^E(h)}{h}\sim n\in\ZM$.
  \end{enumerate}
  $\bullet$ Moreover, the first two terms in the asymptotic expansion
  of $(\lambda^E)^{\out}(h)$ are the following~:
  \begin{itemize}
  \item for non-zero $E$,
    \begin{eqnarray}
      \lefteqn{\frac{1}{h}\left(\int_{\gamma_1^E} \alpha_0^E -
          \epsilon^E_{1,0} + \epsilon^E_{1,0}ln|\bc|
          -\epsilon^E_{2,0}\arg \bc\right) + {}} \nonumber\\ &
      &{} + \left(\int_{\gamma_1^E} \kappa^E +
        \epsilon^E_{1,1}ln|\bc| -\epsilon^E_{2,1}\arg
        \bc + \mu(\gamma_1^E)\frac{\pi}{2} \right), \nonumber
    \end{eqnarray}
    where $\gamma_1^E$ is a smooth family of loops on $\Lambda_0^E$
    such that $\gamma_1^0$ satisfies the hypothesis of proposition
    \ref{prop:regularisation}, $\bc= \epsilon^E_{1,0} + i
    \epsilon^E_{2,0}$, and $\mu$ is the Maslov cocycle on
    $\Lambda_0^E$. For small $E$, $\mu(\gamma_1^E)$ is a constant
    integer.
  \item for $E=0$,
    \[ \frac{1}{h}\int_{\gamma_1^0} \alpha_0 + I_{\gamma_1^0}(\kappa) +
    \mu(\gamma_1^0)\frac{\pi}{2}, \]
    where $\kappa=\kappa^0$ is the sub-principal form of the system
    $(P_1^0,P_2^0)$.
  \end{itemize}
\end{theo}
The last point of this theorem follows from the expression for
non-zero $E$ (recall that $\Lambda^E$ is then smooth) via proposition
\ref{prop:regularisation}~: indeed, if we let
$\tilde{\bc}=(f^E)^{-1}(0)$, we have
$\Lambda_{\tilde{\bc}}=\Lambda_0^E$, and the family of loops
$\gamma_{\tilde{\bc}}$ satisfies the hypothesis of proposition
\ref{prop:regularisation}. The 1-form $\kappa^E$ is defined by~:
\[ (\kappa^E.\ham{p_1^0},\kappa^E.\ham{p_2^0}) =
(df^E)^{-1}.(r^E_1,r^E_2). \] Writing $E$ as a function of
$\tilde{\bc}$, this shows that $\kappa^E$ is the sub-principal form of
a system with principal symbol $(p_1^0,p_2^0)$ and sub-principal
symbol $(df^E)^{-1}.(r^E_1,r^E_2)$. The value at $m$ of this
sub-principal symbol is $r^0(m)$, hence proposition
\ref{prop:regularisation} implies that
\[ \int_{\gamma_1^E} \kappa^E + \re(\bep_1^0\ln M^{-1}\tilde{\bc}) \]
has a limit $I$ as $E\fleche 0$. Because $\bep_1^E$ tends to
$\bep_1^0$ and $\bc=\bep_0^E$ is tangent to $M^{-1}\tilde{\bc}$ at
$\tilde{\bc}=0$, this limit is equal to that of
\[ \int_{\gamma_1^E} \kappa^E + \re(\bep_1^E\ln \bc). \]
Moreover, we saw in (\ref{equ:regularisation}) that it only depends on
the value of the sub-principal form on $\Lambda_0$. But, on
$\Lambda_0$, $(df^0)^{-1}.(r^0_1,r^0_2) = r^0$. Therefore, $I$ is
equal to the principal value $I_{\gamma_0}(\kappa)$ of a system with
principal symbol $p^0$ and sub-principal symbol $r^0$. \cqfd

%The uniqueness statement about the semi-classical constants
%$\epsilon^E_j(h)$, also annouced after lemma \ref{lem:FNQ}, will be
%treated in paragraph \ref{invariant}.

The rest of the proof of this theorem will be split into two
propositions \ref{prop:Chfunction} and \ref{prop:outer}.

We know from the geometry of the Lagrangean fibration (section
\ref{sec:fibration}) that a neighborhood of the principal Lagrangean
$\Lambda_0$ also contains regular tori $\Lambda_0^E$, $E\neq 0$ ($E$
small) for which we already derived the quantization conditions
(theorem \ref{theo:quantization}). Therefore, one of our further tasks
will be to show how these regular quantization conditions are
contained in the ``singular'' quantization conditions announced
above. This link is provided in paragraph \ref{sec:recover}

Finally, this theorem can be applied to ``zoom'' to a $O(h)$ scale~:
let,  as in  remark \ref{rem:dependence}, 
$(\tilde{P}_1^t(h),\tilde{P}_2^t(h))$ be the semi-classical \cis\
defined by
\[ \tilde{P}_j^t(h)\egdef P_j^{ht}(h), \] where $t$ is a
2-dimensional bounded parameter.
From the constant $\bep^E(h)=(\epsilon_1^E(h),\epsilon_2^E(h))$ given by lemma
\ref{lem:FNQ} we define $\tilde{\bep}^t(h) \egdef \frac{1}{h}\bep^{ht}(h)$. This
$\tilde{\bep}^t$ admits an asymptotic expansion in positive powers of
$h$ whose coefficients are smooth functions of $t$. The principal term
is given by~:
%\[ \tilde{\bep}_0^t = \bep_0^0 = 0, \quad \textrm{and} \]
\[ \tilde{\bep}_0^t = \deriv{\bep_0^E}{E}_{\restr E=0}.t + \bep_1^0 \]
\[ \quad = -M^{-1}.\left(\deriv{f^E}{E}_{\restr E=0}(0).t + r(m)\right), \]
that is to say, in view of (\ref{equ:symbolsE}), 
\[ \tilde{\bep}_0^t = -M^{-1}.\tilde{r}^t(m). \]
Here $\tilde{r}^t$ denotes the joint sub-principal symbol of the system
$(\tilde{P}_1^t,\tilde{P}_2^t)$. Recall that its principal symbol is
$\tilde{p}=p^0$ and does not depend on $t$.

Similarly, the holonomy coefficient $(\lambda^E)^{\textrm{\tiny
  out}}(h)\sim \sum_{k=-1}^\infty (\lambda^E)^{\textrm{\tiny
  out}}_k$ has a $t$-dependent version $(\tilde{\lambda}^t)^{\textrm{\tiny
  out}}(h)$ whose first terms are~:
\[ \frac{1}{h}(\lambda^0)^{\out}_{-1} + \left(
  \deriv{(\lambda^{ht})^{\out}_{-1}}{h}_{\restr h=0} +
  (\lambda^0)^{\out}_0 \right) + O(h). \] We know from the
theorem that $\frac{1}{h}(\lambda^0)^{\out}_{-1} =
\frac{1}{h}\int_{\gamma_1^0}\alpha_0$. For the other term, using the
identification $\bep^E=\epsilon_1^E+i\epsilon_2^E$, we have, for
  non-zero $E$~:
\[ \deriv{(\lambda^E)^{\out}_{-1}}{E} = \deriv{}{E}
  \int_{\gamma_1^E} \alpha_0^E - \deriv{}{E} \re(\bep_0^E -
  \bep_0^E\ln\bep_0^E) \]
\[ = \deriv{}{E} \int_{\gamma_1^E} \alpha_0^E +
  \re(\ln(\bep_0^E)\deriv{\bep_0^E}{E}).\]
Therefore, 
\[ \deriv{(\lambda^{ht})^{\out}_{-1}}{h}_{\restr h=0} +
(\lambda^0)^{\out}_0 = \]
\[ = \lim_{h\fleche 0}
(\deriv{}{h}\int_{\gamma_1^E}\alpha^E_0) +\int_{\gamma_1^E}\kappa^E
+
\re\left((\bep_1^E+\deriv{\bep_0^E}{E}.t)\ln\bep_0^E\right)
+ \mu(\gamma_1^E)\frac{\pi}{2},\] with $E=ht$.
Using lemma \ref{lem:action}, one easily shows that, for fixed $h=h_0\neq
0$, $\deriv{}{h}_{\restr h=h_0}\int_{\gamma_1^E}\alpha^E_0$
is equal to $\int_{\gamma_1^E}\delta\kappa^E$, where the 1-form
$\delta\kappa^E$ is given by
\[ \delta\kappa^E.\ham{p^E} = -\deriv{f^E}{E}.t .\]
Therefore, the sum $\delta\kappa^E+\kappa^E$ is the sub-principal form
of a system with principal symbol $p^E$ and sub-principal symbol
$\deriv{f^E}{E}.t + r^E$, whose value at $m$ is equal to
$\tilde{r}^t(m)$. It follows that
\[ \lim_{h\fleche 0} \int_{\gamma_1^E}(\delta\kappa^E+\kappa^E) +
\re\left((\bep_1^E+\deriv{\bep_0^E}{E}.t)\ln\bep_0^E\right)
= I_{\gamma_1^0}(\delta\kappa^0+\kappa). \]
But $\delta\kappa^0+\kappa$ is the 1-form on $\Lambda_0^0$ that takes
the value $\deriv{f^E}{E}_{\restr E=0}.t + r^0$ on
$(\ham{p_1^0},\ham{p_2^0})$~: it is exactly the sub-principal form of
the system $(\tilde{P}^t_1,\tilde{P}^t_2)$. We have proved the
following result~:
\begin{coro}
  \label{coro:global}
  The system $(\tilde{P}^t_1,\tilde{P}^t_2)$ defined above admits a
  microlocal solution $u_h^t$ on a neighborhood of the critical
  Lagrangean $\Lambda_0$ if and only if the following two conditions
  are satisfied~:
  \begin{enumerate}
  \item ``first quantization condition''~:
    \[  \tilde{\epsilon}_2^t(h) \in \ZM + O(h^\infty);\]
  \item ``second quantization condition''~:
    \[ (\tilde{\lambda}^t)^{\out} +
    n\frac{\pi}{2} - \tilde{\epsilon}_1^t\ln(2h) -
    2\arg
    \Gamma\left(\frac{i\tilde{\epsilon}_1^t+1+n}{2}\right)
    \in 2\pi\ZM + O(h^\infty), \]
    where $\tilde{\epsilon}_2^t(h)\sim n\in\ZM$.
  \end{enumerate}
  $\bullet$ Moreover, the first two terms in the asymptotic expansion
  of $(\tilde{\lambda}^t)^{\out}(h)$ are the following~:
  \[ \frac{1}{h}\int_{\gamma_1^0} \alpha_0 +
  \ I_{\gamma_1^0}(\tilde{\kappa}^t) + \mu(\gamma_1^0)\frac{\pi}{2},\]
  where $\gamma_1^0$ is a loop on $\Lambda_0$ that satisfies the
  hypothesis of proposition \ref{prop:regularisation} and $\tilde{\kappa}$
  is the sub-principal form of the system.
\end{coro}

\begin{rema}
  This corollary makes it easy to show that, provided the first
  quantization condition is satisfied, then the ``semi-classical
  invariant'' $\tilde{\epsilon}_1^t(h)$ is intrinsically defined by
  the microlocal behavior of the system around the critical point.
  This gives an \emph{a posteriori} justification of the claim
  following lemma \ref{lem:FNQ}.
  
  Indeed, let $n$ be the integer such that
  $\tilde{\epsilon}_2^t(h)\sim n$ (we shall prove in paragraph
  \ref{sec:recover} that $\tilde{\epsilon}_2^t(h)$ -- and even
  $\epsilon_2^E(h)$ -- are intrinsically defined by the system around
  $m$, hence so is $n$).

  For bounded $t$, $\tilde{\epsilon}_1^t(h)$ is uniformly bounded as
  $h\fleche 0$. As will be shown in the proof of proposition
  \ref{prop:Chfunction}, for any integer $n\in\ZM$, the function
  $\frac{\Gamma(\frac{1}{2}(ix+1+n))}{\Gamma(\frac{1}{2}(-ix+1+n))}$
  is analytic for $x\in\RM$, which implies that the function $2\arg
  \Gamma(\frac{1}{2}(ix+1+n))$ is locally analytic. Therefore,
  plugging in $x=\tilde{\epsilon}_1^t(h)$ we obtain a
  semi-classical expansion in non-negative powers of $h$, whose terms
  are smooth $\T$-valued functions of $t$.

  This in turn implies that the value of the cocycle
  $[\tilde{\lambda}^t(h)]$ on $\gamma_1^0$, which gives the ``second
  quantization condition'' of corollary \ref{coro:global}, can be
  written under the following form~:
  \[ [\tilde{\lambda}^t(h)](\gamma_1^0) \equiv C^t(h) -
  \tilde{\epsilon}_1^t(h)\ln h + O(h^\infty),\] where $C^t(h)$ is a
  semi-classical expansion in $\geq -1$ powers of $h$ whose terms are
  smooth functions of $t$. Therefore, $-\tilde{\epsilon}_1^t(h)$ is
  uniquely defined as the coefficient of $\ln h$ in the expansion of
  $[\tilde{\lambda}^t(h)](\gamma_1^0)$.
  
  Now, if the system is changed outside a neighborhood of the critical 
  point, then only $(\tilde{\lambda}^t)^{\out}$ is perturbed, which 
  means only a semi-classical perturbation that has an asymptotic 
  expansion in powers of $h$, which does not affect the $\ln h$ 
  coefficient.
\end{rema}

\begin{rema}
  \label{rema:colin-p3}
  In dimension 1, the singular Bohr-Sommerfeld conditions of Colin de
  Verdi\`ere and Parisse has been recently extended \cite{colin-p3} to
  handle the case of multiple critical points of saddle type on the
  critical Lagrangean.  In our case as well, if we restrict to \ff\
  singularities, I believe that the results presented here can be
  extended to treat $k>1$ \ff\ points on $\Lambda_0$.  This basically
  holds because the topology of such a $\Lambda_0$ is just a $k$-times
  pinched torus, and the monodromy if equal to the $k$-th power of the
  simple monodromy.  However, including other kinds of singularities
  should require a more difficult analysis, starting from the
  topological description of $\Lambda_0$.  Understanding the work of
  Fomenko (e.g. \cite{fomenko}), should probably help coping with this
  difficulty.
\end{rema}

\subsubsection{The \mi\  bundle}
As a starting point, we shall assume that the system
fulfills the first quantization condition (theorem
\ref{theo:quantif1}). This is of course a necessary condition for the
existence of a global solution.

From this assumption, there exists a \mi\ solution $(u_h^E)^0$ on a
neighborhood of the critical point $m$. As in paragraph \ref{sec:bundle},
one can thus form a \mi\ line bundle $\mathcal{L}_h(\Lambda_0^E)$ over
the topological manifold 
$\Lambda_0^E$, with well-defined holonomy, by means of the trivializations
(\ref{equ:trivialize}). Of course this yields a \emph{family} of line
bundles depending on the 2-dimensional parameter $E$, but gathering
the results of paragraph \ref{sec:dependence} and proposition
\ref{prop:dimension} we see that the trivialization functions in
(\ref{equ:trivialize}) can be chosen to depend smoothly on $E$.

As a consequence, the holonomy of $\mathcal{L}_h(\Lambda_0^E)$
depends smoothly on $E$, in the sense that the elements
$(c^E)^{\alpha\beta}$ of the \v Cech
cocycle $c^E(h)$ representing it  are $\Cinf$ functions of $E$ with
values in ${\bar{\CM_h}}^*$.

From this picture, we already know the abstract form of the
quantization condition we are looking for~: it is exactly the
parameter dependent version of proposition \ref{prop:holonomy}. In other
words, there exists a microlocal solution $u_h^E$ on a neighborhood
of $\Lambda_0$ if and only if for any family $\gamma^E$ of loops in
$\Lambda_0^E$ that lie in that neighborhood, the value of $c^E(h)$ on
$\gamma^E$ is uniformly $1+O(h^\infty)$ (here, we are viewing $c^E$
as a homomorphism from $\pi_1(\Lambda_0^E)$ to ${\bar{\CM_h}}^*$). Let
us now explicit this condition.

Let $\Omega$ be a neighborhood of $m$ in $T^*X$ in which the normal
form of lemma \ref{lem:FNQ} applies. Then, as the system admits a
microlocal solution on $\Omega$, any loop contained in
$\Omega\cap\Lambda_0^E$ has microlocally trivial holonomy. But for
every small $E\neq 0$, $\Omega$ contains a generator $\gamma_2^E$ of
$\pi_1(\Lambda_0^E)$, namely the one associated with the $S^1$
rotational symmetry of the fibration (it is an integral curve of the
vector field $\ham{2}$ of paragraph \ref{sec:geometrie}). Such a
$\gamma_2^E$ also exists for $E=0$ but it is then homotopic to
$\{m\}$. Anyhow, this yields~:
\[ c^E(\gamma^E_2) = 1 + O(h^\infty), \] uniformly for small $E$
(this is of course due to the fulfillment of the first quantization
condition).

Thus $c^E$ restricts to the quotient group
$\pi_1(\Lambda_0^E)/(\gamma^E_2)$. For fixed $E$, every
$\pi_1(\Lambda_0^E)/(\gamma^E_2)$ is isomorphic to $\ZM$. What's more,
proposition \ref{prop:monodromie} implies
that the group bundle $\pi_1(\Lambda_0^E)/(\gamma^E_2)\fleche U\ni E$ over
a neighborhood $U$ of $0$ in $\RM^2$ is trivial~:
\[ \pi_1(\Lambda_0^E)/(\gamma^E_2) \simeq U\times \ZM. \] Let
$\bar{\gamma}_1^E$ be the generator of
$\pi_1(\Lambda_0^E)/(\gamma^E_2)$ equal to $(E,1)$ in the above
trivialization. The ``second quantization condition'' is now the
requirement that $c^E(\bar{\gamma}_1^E) = 1 + O(h^\infty)$.

Let us now examine that statement more closely. Recall that the line
bundle $\mathcal{L}_h(\Lambda_0^E)$ was defined through a set of local
trivializations subordinated to a cover $\cup\Omega_\alpha$ of a
neighborhood of $\Lambda_0$ in $T^*X$, such that on every
$\Omega_\alpha$ existed a non-trivial \mi\  solution
$(u_h^E)^\alpha$ of the system. We can assume that the critical point
$m$ is contained in a unique element of that cover, which we will
denote by $\Omega_0$. Let us pick up a simple loop $\gamma_1^0(t)$,
$t\in[0,1]$, in $\Lambda_0$ representing the quotient class
$\bar{\gamma}_1^0$, and enumerate $\Omega_0$,
$\Omega_1,\dots,\Omega_\ell,\Omega_0$ the $\Omega_\alpha$'s encountered
by $\gamma_1^0(t)$ on its way from $m$ back to $m$ again. We will here
choose $\gamma_1^0$ so that it goes away from $m$ via the local stable
manifold of $\Lambda_0$ and returns back to $m$ via the local unstable
manifold. Also, we will suppose that
$\Omega_0\cup\Omega_1\cup\Omega_\ell\subset \Omega$.

Now for every fixed $E$, we choose a representative $\gamma_1^E$ of
$\bar{\gamma}_1^E$. Because of the triviality of the fibration
outside $m$ (\ref{equ:triv-fibration}), the paths $\gamma_1^E$ can be
homotopically deformed so as to be close to $\gamma_1^0$, in the sense
that their images enter $\bigcup_{i=0}^\ell \Omega_i$.

Then by definition, the value of $c^E$ on $\gamma_1^E$ is the
product~:
\[ c^E(\gamma_1^E) =
(c^E)^{0,1}(c^E)^{1,2}\cdots(c^E)^{\ell-1,\ell}(c^E)^{\ell,0}. \]
We are going to calculate this expression by splitting it in two parts~: a
local one -- the product $(c^E)^{0,1}(c^E)^{\ell,0}$ -- and an outer
one -- the rest. For this purpose, we are free to fix the solutions
$(u_h^E)^i$ on $\Omega_i$ used for the local trivialization of
$\mathcal{L}_h(\Lambda_0^E)$; let us make the following choice~:

\noindent $\bullet$ choice 1~: for $i\neq 0,1,\ell$, $(u_h^E)^i$ is any
classical oscillatory integral constructed as in proposition
\ref{prop:WKB}, with invariant half density given by the canonical
Liouville density on $\Lambda_0^E$.

For $i=0,1,\ell$, we make use of the local analysis of the preceding
paragraph \ref{sec:microloc}. To begin with, we transform the system on
$\Omega$ to the normal form of lemma \ref{lem:FNQ} by a unitary \fio\
$U(h)$. The \mi\ solutions of the transformed system on $\Omega$ are
spanned by~:
\[ u^E_{\textrm{\tiny exact}} = \frac{1}{r}
e^{i\frac{\epsilon_1^E}{h}\ln r}e^{i\frac{\epsilon_2^E}{h}\theta}. \]
Because of the first quantization condition, we know
$\frac{\epsilon_2^E}{h}$ to be microlocally equal to an integer
$n\in\ZM$. For any real $\varepsilon$, let $u_{\varepsilon,n}$ be
the tempered distribution on $\RM^2$~:
\[ u_{\varepsilon,n} = \frac{1}{r} e^{in\theta}e^{i\varepsilon\ln r},
\] so that 
$u^E_{\textrm{\tiny exact}} =
u_{\frac{\epsilon_1^E}{h},\frac{\epsilon_2^E}{h}}$.
The distribution $u_{\varepsilon,n}$ is a solution of the system
\begin{equation}
  \label{equ:homogeneous}
  \left\{\begin{array}{l} \frac{1}{i}\deriv{}{\theta} u_{\varepsilon,n}= 
  n u_{\varepsilon,n} \\
      \frac{1}{i}(r\deriv{}{r} + 1) u_{\varepsilon,n} = \varepsilon 
      u_{\varepsilon,n} \end{array}\right.
\end{equation}
 on $\RM^2$, and it is easy to see that its Fourier transform is a
solution of the same system with $\varepsilon$ changed to
$-\varepsilon$. Therefore
$u_{\frac{\epsilon_1^E}{h},\frac{\epsilon_2^E}{h}}$ and
$\fourier^{-1}
(u_{\frac{-\epsilon_1^E}{h},\frac{\epsilon_2^E}{h}})$ are two
solutions of our transformed system on $\Omega$. This allows us to
complete our choice of local solutions by the following~:

\noindent $\bullet$ choice 2~: $\left\{\begin{array}{lr} (u_h^E)^0 \sim
U(u_{\frac{\epsilon_1^E}{h},\frac{\epsilon_2^E}{h}}) &
\textrm{ on } \Omega_0 \\ (u_h^E)^1 \sim
U(u_{\frac{\epsilon_1^E}{h},\frac{\epsilon_2^E}{h}}) &
\textrm{ on } \Omega_1 \\ (u_h^E)^\ell \sim U(\fourier^{-1}
u_{\frac{-\epsilon_1^E}{h},\frac{\epsilon_2^E}{h}}) & \textrm{
on } \Omega_\ell \end{array}\right.$

\subsubsection{The local holonomy at $m$}

We give here an explicit expression for the product
$(c^E)^{1,0}(c^E)^{0,\ell}$, which is defined to be the constant
$C^E(h)$ such that, on $\Omega_0\cap\Omega_\ell$, $(u_h^E)^0 \sim C^E
(u_h^E)^\ell$.
\begin{prop}
  \label{prop:Chfunction}
  The unique microlocal constant $C^E(h)$ such that $(u_h^E)^0 \sim
  C^E (u_h^E)^\ell$ on $\Omega_0\cap\Omega_\ell$ is given by~:
  \[ C^E(h) \sim i^{-n} {(2h)}^{i\epsilon_1^E/h}
  \frac{\displaystyle
  \Gamma\left(\frac{i\epsilon_1^E/h+1+n}{2}\right)} {\displaystyle
  \Gamma\left(\frac{-i\epsilon_1^E/h+1+n}{2}\right)}, \] with
  $\frac{\epsilon_2^E}{h}\sim n\in\ZM$.
\end{prop}
\demo Applying the \fio\  $U_h$, we see that $C_h^E$ is also defined to
be the constant such that, microlocally on $\Omega_0\cap\Omega_\ell$,
\[  \fourier
u_{\frac{\epsilon_1^E}{h},\frac{\epsilon_2^E}{h}} \sim C_h^E 
u_{\frac{-\epsilon_1^E}{h},\frac{\epsilon_2^E}{h}}. \] It turns out
that the constant $C_h^E$ thus defined is a microlocal version of an
exact problem concerning the homogeneous distributions $u_{\varepsilon,n}$ on
$\RM^2$.

We know that the tempered distributions $u_{\varepsilon,n}$ are solutions
of the system (\ref{equ:homogeneous}) in $\RM^2$. They are in fact the
only ones -- up to multiplicative constants, naturally. To see this,
one first restricts to $\RM^2\setminus\{0\}$, where the result is
standard, and then prove that the system (\ref{equ:homogeneous}) cannot
admit linear combinations of derivatives of the Dirac distribution  as
solutions. Now, $\fouriero^{-1}u_{-\varepsilon,n}$ is also a tempered
distribution and a solution of the same system, which leads to the
first point of the
following proposition~:
\begin{prop}[\cite{cassels,guelfand}]
  \begin{enumerate}
  \item For any $(\varepsilon,n)\in\RM\times\ZM$, there exists a
    unique constant $C(\varepsilon,n)$ such that~:
    \[ \fouriero u_{\varepsilon,n} = C(\varepsilon,n) u_{-\varepsilon,n}. \]
  \item $C(\varepsilon,n)$ has the following expression~:
    \begin{equation}
      \label{equ:Cfunction}
      C(\varepsilon,n) = i^{-|n|} 2^{i\varepsilon}
      \frac{\Gamma(\frac{i\varepsilon+1+|n|}{2})}
      {\Gamma(\frac{-i\varepsilon+1+|n|}{2})}.
    \end{equation}
  \end{enumerate}
\end{prop}
The computation of $C(\varepsilon,n)$ is done by testing
$u_{\varepsilon,n}$ on derivatives of Gaussian functions
$(\deriv{}{z})^ne^{\frac{z\overline{z}}{2}}$ and
$(\deriv{}{\overline{z}})^ne^{\frac{z\overline{z}}{2}}$
($z=re^i\theta$).

\begin{rema}
  Apparently, this formula appeared for the first time in Tate's
  thesis, during the year 1950. However, it has not been published
  until 1967, when it appeared in the book~\cite[chapter
  XV]{cassels}. But at the same time as Tate wrote his thesis, the
  idea was in the air. Gelfand was studying homogeneous distributions
  in the real and published a closely related formula in
  \cite{guelfand-schapiro}. The complex version, that is, exactly the
  result mentioned in the above proposition, appeared in the
  \emph{addendum} of the French edition of the book by Gelfand \&
  al. on distributions \cite{guelfand}, where it is claimed to be
  published for the first time; the author was certainly unaware of
  Tate's work. Also interestingly enough, few years before (1951), a
  slightly different version of the same formula can be found in
  Bochner's work \cite{bochner}. It is actually a generalization to
  homogeneous distributions on $\RM^k$.
\end{rema}

Now, let us (momentarily) denote by $D(\varepsilon,n)$ the quantity~:
\[   D(\varepsilon,n) = i^{-n} 2^{i\varepsilon}
\frac{\Gamma(\frac{i\varepsilon+1+n}{2})}
{\Gamma(\frac{-i\varepsilon+1+n}{2})}.
\] 
The Gamma functions involved are analytic functions of $\varepsilon$
that, for negative odd $n$, have simple poles at
$\varepsilon=0$. The quotient $D(\varepsilon,n)$ is therefore analytic
-- with no pole -- on the real line (we won't consider here analytic
continuation to complex $\varepsilon$). More precisely, for $n=-2m-1$,
$m\geq 0$, as the residues of the numerator and the denominator are
the same, we have~: $D(0,-2m-1) = \mathbf{-}i^{-n}$.  For non-zero
$\varepsilon$, we can use the complement relation to prove that
$D(\varepsilon,n)=D(\varepsilon,-n)$. Then of course this also holds
for $\varepsilon=0$. Since $D(\varepsilon,n)$ and $C(\varepsilon,n)$
coincide for $n\geq 0$, we deduce that they are equal for all
$n\in\ZM$.

Finally, using the relation $\fourier u(\xi,\eta) =
\frac{1}{h}\fouriero u (\frac{\xi}{h},\frac{\eta}{h})$, we terminate
the proof of proposition \ref{prop:Chfunction}. \cqfd

\begin{rema}
  One of the important features of the constants $C(\varepsilon,n)$
  (and hence of $C^E$) is that they have modulus one. It might be
  interesting to see that, without any explicit calculation, this
  comes as a consequence of the unitarity of the Fourier transform
  with respect to the $L^2$ norm. As this does not appear -- to my
  knowledge -- in the literature, here is the argument~:

  Acting on functions of the Schwartz space $\mathcal{S}(\RM^2)$,
  $\fouriero$ commutes with the operator
  $\frac{1}{i}\deriv{}{\theta}$, and hence leaves stable its
  eigenspaces $\mathcal{E}_n$, associated to the eigenvalues $2\pi
  n$. These are of course the subspaces of the Fourier series
  decomposition in $\theta$ variable~: $\mathcal{E}_n =
  \{e^{in\theta}f(r)\}$. Here, $f$ is a rapidly decreasing function.
  Therefore, there exists a transformation $\fouriero^n$ on
  $\mathcal{E}_n$ such that
  \[ \fouriero(e^{-in\theta}f(r))(\rho,\alpha) =
  e^{-in\alpha} \fouriero^n(f)(\rho) .\] $\fouriero$ being an isometry
  for the norms $L^2(\RM^2,r dr d\theta)$ and $L^2(\RM^2,\rho d\rho
  d\alpha)$, we deduce that $\fouriero^n$ is an isometry from
  $L^2([0,\infty[,rdr)$ into $L^2([0,\infty[,\rho d\rho)$.
%Note also
%that $\fouriero^n(\bar{f})=(-1)^n\overline{\fouriero^n f}$.

  On the other hand, let $\mathcal{G}$ be the transformation assigning
  to a function  $f(r)\in L^2([0,\infty[,rdr)$ the function
  \[\mathcal{G}(f)(\varepsilon) =
  \int_0^\infty r^{i\varepsilon-1}f(r)rdr. \]
  \begin{lemm}
    $\mathcal{G}$ is an isometry from $L^2([0,\infty[,rdr)$ into
    $L^2(\RM,d\varepsilon)$.
  \end{lemm}
  \demo The change of variable $x=\ln r$ allows us to reduce the problem
  to the unitarity of the Fourier transform~: let $g(\varepsilon) =
  \mathcal{G}(f)(\varepsilon)$. We find
  $g(\varepsilon)=\mathcal{F}(f(e^x)e^x)(-\varepsilon)$, which gives
  $\|g\|_2=\|f(e^x)e^x\|_2 = \int_0^\infty |f(r)|^2 rdr = \|f\|_2$.
  \cqfd
  
  We know from the definition of $C(\varepsilon,k)$ that for any
  function of the form
  \[ \phy=\fouriero(e^{-in\theta}f(r)),\]
  we have~:
  $\pscal{\fouriero(u_{(\varepsilon,k)})}{\phy} = C(\varepsilon,k)
  \pscal{u_{(-\varepsilon,k)}}{\phy}$, that is~:
  \[ \pscal{e^{ik\theta}r^{i\varepsilon-1}}{\fouriero \phy} = C(\varepsilon,k)
  \pscal{e^{ik\alpha}\rho^{-i\varepsilon-1}}{\phy}, \]
  or~:
  \[  \pscal{r^{i\varepsilon-1}} {e^{i(k-n)(\theta+\pi)}f(r)} =
  C(\varepsilon,k) 
  \pscal{\rho^{-i\varepsilon-1}}
  {e^{i(k-n)\alpha}{(\fouriero^n f)(\rho)}}.\]
  
  Both sides of the equality vanish if $n-k\neq 0$. So let us assume
  $k=n$. We obtain~:
  \[ \mathcal{G}(f)(\varepsilon) = (-1)^nC(\varepsilon,n)
  \mathcal{G}\fouriero^n f(-\varepsilon), \] which is best depicted by
  the following commutative diagram~:
  \[ \begin{array}{ccc} L^2([0,\infty[,rdr) & \flechedroite{\mathcal{G}} &
    L^2(\RM,\varepsilon) \\ 
    \flechebas{\fouriero^n} & & \flechehaut{(-1)^nC(\varepsilon,n)} \\ 
    L^2([0,\infty[,\rho d\rho) &
    \flechedroite{\check{\mathcal{G}}} 
    & L^2(\RM,-\varepsilon)
  \end{array} \]
  All transformation involved ($\mathcal{G}$, $\fouriero^n$ and
  $\varepsilon\fleche -\varepsilon$) being unitary ones, one can indeed
  deduce that $\varepsilon\mapsto C(\varepsilon,n)$ is a function whose
  modulus is 1. \cqfd
\end{rema}

\begin{rema}
  The above introduced transformations $\mathcal{G}$ and $\fouriero^n$
  are related to Mellin's and Hankel's transforms $\mathcal{M}$ and
  $\mathcal{H}_n$ respectively (see \cite{colombo} for their
  definitions). More precisely, we have~:
  \[ (\mathcal{G}f)(\varepsilon) = (\mathcal{M}f)(s=i\varepsilon+1)
  \textrm{ and } \fouriero^n= i^{-n}\mathcal{H}_n. \]
  (see for instance \cite{dautray}).
  The previous remark thus yields the following remarkable relation
  between these transformations~:
  \begin{equation}
    \label{equ:hankel}
    \forall f, \qquad \mathcal{M}f(i\varepsilon+1) = i^n
    C(\varepsilon,n).(\mathcal{M}\mathcal{H}_n f)(-i\varepsilon+1).
  \end{equation}
\end{rema}
\begin{rema}
  Formula (\ref{equ:hankel}) allows us to give another proof of the
  expression of $C(\varepsilon,n)$. First of all, since
  $\mathcal{H}_{-n}= (-1)^n\mathcal{H}_n$, it readily gives
  $C(\varepsilon,n) = C(\varepsilon,-n)$. Now for $n\geq 0$ choose
  $f_n=r^ne^{-\frac{r^2}{2}}$. $f_n$ is fixed by Hankel's transform of
  order $n$. On the other hand, it is easy to see that
  $\mathcal{M}f_n(s)=2^{\frac{s+n}{2}-1}\Gamma(\frac{s+n}{2})$. This
  yields~:

  \[ C(\varepsilon,n)=
  i^{-n}\frac{2^{s/2} \Gamma((s+n)/2)}
  {2^{\overline{s}/2}\Gamma((\overline{s}+n)/2)},
  \textrm{ with } s=i\varepsilon+1, \overline{s}=-i\varepsilon+1, \]
  which is easily identified with formula (\ref{equ:Cfunction}). \cqfd
\end{rema}

\subsubsection{The outer holonomy}
Let us turn now to the study of the ``outer holonomy'', that is, the
product
\[ (C^E)^{\out} = (c^E)^{1,2}\cdots(c^E)^{\ell-1,\ell}. \]
Recall that its quotient by the local holonomy $C^E$ defined in the
previous paragraph gives the value of the \v Cech cocycle $c^E$ on
$\gamma_1^E$.
\begin{prop}
  \label{prop:outer}
  \begin{itemize}
  \item $|(C^E)^{\out}| = 1 + O(h^\infty)$;
  \item $\arg( (C^E)^{\out})\in\T$ admits an asymptotic
    expansion in $\geq -1$ powers of $h$, whose terms are smooth
    functions of $E$;
  \item for non zero $E$, the first two terms in that expansion are
the following~:
\begin{eqnarray}
  \lefteqn{ \frac{1}{h}\left(\int_{\gamma_1^E} \alpha_0^E -
      \epsilon^E_{1,0} + \epsilon^E_{1,0}ln|\bc|
      -\epsilon^E_{2,0}\arg \bc\right) +{}}\nonumber\\ & &
      {}+\left(\int_{\gamma_1^E} \kappa^E +
      \epsilon^E_{1,1}ln|\bc| -\epsilon^E_{2,1}\arg \bc
      + \mu(\gamma_1^E)\frac{\pi}{2} \right), \nonumber
\end{eqnarray}
where $\bc= \epsilon^E_{1,0} + i  \epsilon^E_{2,0}$.
  \end{itemize}
\end{prop}
\demo
For non zero $E$, $\Lambda_0^E$ is a regular Lagrangean torus, and we
know from proposition \ref{prop:cocycle} that $|c^E(\gamma_1^E)|\sim
1$. This, together with the fact that $|C^E|\sim 1$ too, implies
that $|(C^E)^{\out}| \sim 1$ as well. Then of course,
since we already knew that $(C^E)^{\out}$ is a smooth
function of $E$, the results also holds if $E=0$.

The second point in the proposition is essentially due to our choice
of the microlocal solution $(u_h^E)^\ell$ on $\Omega_\ell$. The
distribution 
\[ \frac{1}{2\pi ih} \int 
e^{\frac{i}{h}\left( \pscal{(x,y)}{(\xi,\eta)} -
    \epsilon_1^E\ln\rho + \epsilon_2^E\alpha \right) }
    \frac{|d\xi\wedge d\eta|}{\rho} \] (by $(\rho,\alpha)$ we always
    denote polar coordinates for $(\xi,\eta)$ in the fiber of
    $T^*\RM^2$) is microlocally on $\Omega_\ell$ equal to a classical
    oscillatory integral on $\Lambda_0^E\cap\Omega_\ell$, associated
    with the canonical Liouville half density on $\Lambda_0^E$, and
    smoothly depending on $E$.

\begin{rema}
  This is of course the pendant of remark
  \ref{rem:quantif1}. Here as well, the phase function $\phy^E=
  \pscal{(x,y)}{(\xi,\eta)} - \epsilon_{1,0}^E\ln\rho +
  \epsilon_{2,0}^E\alpha$ is best expressed in terms of complex
  coordinates $z_1$ and $z_2$~: we have
  \[ \phy^E (z_1,z_2) = \re(\bar{z_1}z_2 - \bep_0^E\ln z_2). \]
\end{rema}

Now, since $(u_h^E)^\ell$ was defined to be the image of that
distribution by an classical \fio, it is also a classical Lagrangean
distribution, as are all the other solutions $(u_h^E)^i$, $i\neq
0$. We deduce from proposition \ref{prop:cocycle} and paragraph
\ref{sec:dependence} that for each $i=1,\dots,\ell-1$, the argument of
$(c^E)^{i,i+1}$ admits a semi-classical expansion in $\geq -1$ powers
of $h$ whose terms smoothly depend on $E$, as functions with values in
$\T$. Hence so does their product $(C^E)^{\out}$.

We turn now to the last point of the proposition, which requires a
closer look at the principal phases of the solutions near the critical
point $m$. Indeed, we know from the non-singular analysis of section
\ref{sec:CBSreg} that the principal phases together with a section of the
Keller-Maslov line bundle form an invariantly defined object that
gives the first two terms of the \v Cech cocycle $\lambda^E(h)$.

First, let us extend the open sets $\Omega_1$ and $\Omega_\ell$ to be
both equal to $\Omega\setminus\{m\}$. This is possible as the
microlocal solutions $(u_h^E)^1$ and $(u_h^E)^\ell$ are actually
well-defined on all $\Omega$. This ensures that for every non zero
$E$, $\bigcup_{i=1}^\ell\Omega_i$ is a full cover of the image of
$\gamma_1^E$.  Still assuming that $E$ is non zero, we get from
theorem \ref{theo:quantization} that $\arg(c^E(\gamma_1^E))$ is, up to
$O(h)$ (and modulo $2\pi$), equal to
\begin{equation}
  \label{equ:outerphase}
\frac{1}{h}\int_{\gamma_1^E} \alpha_0^E + \int_{\gamma_1^E} \kappa^E +
\mu(\gamma_1^E)\frac{\pi}{2}.
\end{equation}

On the other hand, $c^E(\gamma_1^E)$ is by definition the value of the
product \[(C^E)^{\out}(c^E)^{\ell,1}.\] The problem is
thus reduced to the computation of $\arg (c^E)^{\ell,1}$. It is
unchanged by conjugation by a \fio\ on $\Omega$.

Let us denote by $\phy_1^E(z_1)$ and $\phy_\ell^E(z_1,z_2)$ the phase
functions of the oscillatory functions
$u_{\frac{\epsilon_1^E}{h},\frac{\epsilon_2^E}{h}}$ and $\fourier^{-1}
u_{\frac{-\epsilon_1^E}{h},\frac{\epsilon_2^E}{h}}$. They have the
following expressions (modulo $2\pi h\ZM$)~:
\[ \phy_1^E(z_1) \equiv  \re(\bep^E_0\ln\bar{z_1}) \textrm{ and }
\phy_\ell^E(z_1,z_2) \equiv \re(\bar{z_1}z_2 -
\bep_0^E\ln z_2). \] As we already saw, for non zero $E$,
they represent -- of course -- the same Lagrangean manifold
$\Lambda_0^E=\{\bar{z_1}z_2 = \bep_0^E\}$ (the terms
$\re(\bep^E_0\ln\bar{z_1})$ and
$\re(\bep_0^E\ln z_2)$ are actually Legendre transforms
of each other). Their difference is a function of $c=\bar{z_1}z_2$~:
\[ \phy_\ell^E(z_1,z_2)-\phy_1^E(z_1) \equiv
\re(\bc) -\re(\bep^E_0\ln \bc), \] whose
real Hessian is the $2\times 2$ matrix whose lines are respectively
identified with $\bar{\bep^E_0}/\bar{\bc}^2$ and
$-i\bar{\bep^E_0}/\bar{\bc}^2$. Its determinant is
$-|\frac{\bep^E_0} {\bc^2}|^2$. If $E\neq 0$ then
$\bep^E_0\neq 0$ so the signature of that matrix is
zero. This means that the value of Maslov's cocycle between
$u_{\frac{\epsilon_{1h}^E}{h},\frac{\epsilon_{2h}^E}{h}}$ and
$\fourier^{-1}
u_{\frac{-\epsilon_{1h}^E}{h},\frac{\epsilon_{2h}^E}{h}}$ is zero.

The principal phases of those oscillatory functions are expressed by a
formula similar to the one expressing the phase functions, but with a
contribution of the sub-principal terms of $\epsilon^E_{ih}$
included. The two principal terms in the argument of $(c^E)^{\ell,1}$ is
their difference~:
\[ \arg (c^E)^{\ell,1} \equiv \frac{1}{h}\left( \re(\bep^E_0) -
\re(\bep^E_0\ln\bep^E_0) \right) - \re(\bep^E_1\ln\bep^E_0) +
O^E(h). \] Here, we have fixed $E\neq 0$, so the $O^E(h)$ depends on
$E$.

Finally, together with equation (\ref{equ:outerphase}), we get that the
first two terms in the expansion of $(C_h^E)^{\out}$
are, for non zero $E$~:

\[ \frac{1}{h}\left(\int_{\gamma_1^E} \alpha_0^E -  \re(\bep^E_0)
  + \re(\bep^E_0\ln\bep^E_0) \right) + \left(\int_{\gamma_1^E}
\kappa^E + \mu(\gamma_1^E)\frac{\pi}{2} +
\re(\bep^E_1\ln\bep^E_0)\right), \]
and the proof of proposition \ref{prop:outer} is completed. \cqfd

Of course, this also terminates the proof of theorem
\ref{theo:global}. \cqfd

\subsubsection{Recover the regular conditions}
\label{sec:recover}
The aim of this paragraph is to show how the quantization conditions
of theorem \ref{theo:global}, when restricted to a non-singular
Lagrangean, are equivalent to the regular quantization conditions of
theorem~\ref{theo:quantization}.

Let us fix $E\neq 0$. We get a smooth holonomy 1-form $[\lambda_h^E]$
whose action on a basis of the $H_1$ of the torus $\Lambda_0^E$ gives
the two quantization conditions. Let $\gamma_1^E$ and $\gamma_2^E$ be,
as previously, loops on $\Lambda_0^E$ representing such a basis. 

The loop $\gamma_2^E$ can be assumed to lie in the neighborhood
$\Omega$ where the normal form applies. In the linearized coordinates,
one can choose $\gamma_2^E$ to be a simple orbit of the $\ham{2}$
vector field. On each point of that orbit we saw that the solutions
were spanned by 
\[\frac{1}{r} e^{
  i\frac{\epsilon_1^E}{h}\ln r} e^{
  i\frac{\epsilon_2^E}{h}\theta}.
\]
Therefore, the value of $[\lambda_h^E]$ on $\gamma_2^E$ is exactly the
monodromy coefficient  $\frac{2\pi\epsilon_2^E(h)}{h}$. The
quantization condition given by theorem \ref{theo:quantization} is
thus
\[ \frac{2\pi\epsilon_2^E}{h} \in 2\pi\ZM + O(h^\infty),\]
which is indeed the same as the first singular quantization condition 
of theorem \ref{theo:global}.

\begin{rema} 
  \label{rem:recover}
  This simple observation gives rise, for non zero $E$, to an
  intrinsic definition of the semi-classical expansion
  $\epsilon_2^E(h)$ introduced in lemma \ref{lem:FNQ}, namely~:
  \[ \epsilon_2^E(h) = \frac{1}{2\pi}\int_{\gamma_2^E} \lambda^E(h). \]
  Since $\epsilon_2^E$ is smooth at $E=0$, this implies that it is
  intrinsically defined by the \mi\ system on a neighborhood of $m$,
  that is, independently of the \fio\ $U$. This will be restated in
  the next paragraph.
\end{rema}

Let us turn now to the second quantization condition. We know that the
action of $[\lambda^E]$ on $\gamma_1^E$ is given by the argument of
the quotient
$(C^E)^{\out}/C^E$. Using Stirling's formula for
$\arg(\Gamma(z))\equiv \im(\ln \Gamma(z))$, with
$z=\frac{i\bar{\bep}(h)}{2h} +\frac{1}{2}$
($\bep(h)=\epsilon_1^E(h) + i \epsilon_2^E(h)$; recall
that if $E\neq 0$, the principal term $\bep_0$ of
$\bep(h)$ is non-zero), one can perform a very nice
cancellation of undesirable terms in the local holonomy $C^E$.

Indeed, one finds that
\[ 2\arg(\Gamma(z))\equiv
\frac{1}{h}\left(\re(\bep\ln\bep) -
  \re(\bep) + \im(\bep)\frac{\pi}{2} -
  \re(\bep\ln(2h))\right) +O(h). \] Here of course,
  $O(h)$ depends on $E$. This gives
\[ \arg (C^E) \equiv
  \frac{1}{h}\left(\re(\bep\ln\bep) -
  \re(\bep) \right) + O(h), \] that is,
\[ \arg (C^E) \equiv
\frac{1}{h}\re(\bep_0\ln\bep_0-\bep_0)
+ \re(\bep_1\ln\bep_0) + O(h). \] We thus find that,
if $E\neq 0$,
\[ \int_{\gamma_1^E}\lambda^E(h) =  \frac{1}{h}\int_{\gamma_1^E}
\alpha_0^E + \int_{\gamma_1^E} \kappa^E + \mu(\gamma_1^E)\frac{\pi}{2}
+ O(h),\] which gives the regular quantization condition of
theorem~\ref{theo:quantization}.

%\subsubsection{The semi-classical invariants}
%\label{invariant}
%\begin{prop}
%  the semi-classical expansions $\epsilon^E_1(h)$ and
%  $\epsilon^E_2(h)$ defined in lemma \ref{lem:FNQ} are independent on the
%  \fio\ $U(h)$ used to perform the microlocal normal form, and
%  therefore completely characterize the \mi\ solutions of the system
%  (\ref{equ:systeme}) near the critical point. They are called the
%  \emph{semi-classical invariants} of the system.
%\end{prop}
%\AAAA
%\demo the uniqueness of $\epsilon^E_2(h)$ has just been proved in the
%remark of the previous paragraph \ref{sec:recover}.

%\subsection{exemple: le pendule sphÚrique}
%\input{pendule.tex}

%-------------------------------------------------------------------------

\section{Structure of the joint spectrum}
\label{sec:spectrum}
This section is devoted to the application of our microlocal study to
the precise study of the joint spectrum of 2-degree of freedom
quantum integrable systems around the critical value of a \ff\
singularity.

More precisely, we shall assume in all this section that we are given
two commuting essentially self-adjoint $h$-\pdo s of order zero
$P_1(h),P_2(h)$ on a 2-dimensional manifold $X$, with the following
properties~:
\begin{itemize}
\item the joint principal symbol $p=(p_1,p_2)$ is proper with
  connected leaves;
\item $p$ admits a critical point $m$ of \ff\ type, with critical
  value $0\in\RM^2$, such that the critical set $\Lambda_0=p^{-1}(0)$
  has $m$ as its only critical point;
\item at least one of the $P_j(h)$'s is a classical elliptic \pdo\  in
  the sense of H\"ormander (for instance, $P_1(h)$ is a \schr\ 
  operator~: 
\[ P_1(h) = \frac{h^2}{2}\Delta_g + V(x),\]
for some Riemannian metric $g$ on $X$).
\end{itemize}
The last point implies that the solutions $u_h$ of $P_j(h)u_h=0$ ($j=1,2$)
cannot have a wave front set at infinity (see \cite{courscolin}),
which allows us to microlocalize $u_h$ on the Lagrangean level set
$\Lambda_0=p^{-1}(0)$.

Note also that $P_1(h)$ and $P_2(h)$ really commute (and not only
modulo $O(h^\infty)$). This assumption is perhaps unnecessary,
but we take it in order to use the definition of the joint spectrum
from \cite{charbonnel}.

A typical example of this situation would be a \schr\  
operator on a surface of revolution with a radial potential, like 
the quantum spherical pendulum \cite{duist-cushman} or the quantum 
Champagne bottle \cite{child}.  As a matter of fact, we will 
illustrate our results with the latter example, thanks to numerics 
kindly provided by M.Child.
%The organisation of this section goes as follows. After recalling that
%the joint spectrum is discrete near the origin, we prove that, for $h$ small
%enough, all the joint eigenvalues are simple. We can then focus on the
%geometry of the joint spectrum, which splits into two parts~:
%\begin{itemize}
%\item far from the critical value $0$, the main property is the
%  non-triviality of the quantum monodromy (defined in
%  \cite{san?}). .....
%\item ....
%\end{itemize}

\subsection{The joint spectrum}
Let $\Sigma(h)$ be the joint spectrum of the operators $P_1(h)$ and
$P_2(h)$.
We know from \cite{charbonnel} that the assumption of properness for
the momentum map $p$ implies that the intersection of $\Sigma(h)$ with
any compact $K\subset \RM^n$ is a discrete spectrum~: each
$E(h)=(E_1(h),E_2(h))\in\Sigma(h)\cap K$ is isolated
and the joint eigenspace 
\[ F_{E(h)}=\{\phy(h)\in L^2(X), (P_1(h)-E_1(h))\phy(h) =
(P_2(h)-E_2(h))\phy(h) = 0 \}\] 
has finite dimension.
\begin{theo}
  \label{theo:simple}
  There is a compact neighborhood $K$ of $0\in\RM^n$ and an $h_0>0$
  such that all the joint eigenvalues $E(h)\in\Sigma(h)\cap K$ for
  $h<h_0$ are simple. Furthermore, the distance between two distinct
  joint eigenvalues is bounded below by some finite power of $h$.
\end{theo}
\demo 
This is a consequence of proposition \ref{prop:dimension}, applied to
the operators $P_j^E(h)\egdef P_j(h)-E_j$. Indeed, let
$K\subset\subset\RM^n$ be a neighborhood of zero such that this
proposition holds uniformly for $E=(E_1,E_2)\in K$. Let $\psi_1(h)$ and
$\psi_2(h)$ be two normalized joint eigenfunctions for the joint
eigenvalue $E(h)$, such that $\pscal{\psi_1}{\psi_2}=0$. If we can
find such eigenfunctions for a set of values of $h$ accumulating at
zero, then $\psi_1(h)$ and $\psi_2(h)$ are both microlocal solutions
(see section \ref{sec:solutions}) to the equation
\[ (P_j(h)-E_j(h))\psi_k(h)\sim 0, \qquad j=1,2\quad k=1,2. \] 
Therefore, their wave front sets are included in a compact $W\subset
T^*X$. Let $\pi:T^*X\fleche X$ be the standard projection and let $U$
be an open neighborhood of $\pi(W)$ such that $\overline{U}$ is
compact. Then,
\[ WF_h(\psi_k)\cap (T^*(X\setminus U)) = \emptyset \qquad(k=1,2), \]
and because $\psi_k$ cannot have any wave front set at infinity, $\psi_k$
is actually $O(h^\infty)$ on $X\setminus U$. 
Consequently if we let $\chi$ be a
smooth function with support in $U$ and equal to 1 on $\pi(W)$ then
\[ \|\chi\psi_k\|_{L^2(X)} = \|\psi_k\|_{L^2(X)}+O(h^\infty) = 1 +
O(h^\infty),\]
\[ \textrm{and } \pscal{\chi\psi_1}{\chi\psi_2} = O(h^\infty).\]
That allows us to replace, in the following, $\psi_k$ by $\chi\psi_k$.

Because of proposition \ref{prop:dimension}, there exists a constant
$C(h)\in\CM_h$ such that 
\begin{equation}
  \label{equ:twoeigfun}
  \psi_1(h)\sim C(h)\psi_2(h) \textrm{  on } T^*U.
\end{equation}
Combined as before with the fact that nothing happens at infinity,
equation (\ref{equ:twoeigfun}) implies, by taking scalar products, that
$C(h)$ simultaneously satisfies $|C(h)|=1+O(h^\infty)$ and $C(h)=0 +
O(h^\infty)$, which is impossible. This proves that for small $h$ all joint
eigenvalues $E(h)$ must be simple.

Now, if $E(h)$ and $E'(h)$ are distinct joint eigenvalues but
differing only by $O(h^\infty)$, the same argument applies to the
eigenfunctions $\psi_1(h)$ and $\psi_2(h)$ associated to $E(h)$ and
$E'(h)$, and thus proves the last statement of the theorem. \cqfd

This theorem extends the results of \cite{charbonnel}, where it is
shown that the eigenvalues are simple as long as $K$ does not meet any
critical value of the momentum map $p$. Note that this latter statement
comes as a consequence of our proposition \ref{prop:dim1reg} as well.

Finally, remark that the proof of theorem \ref{theo:simple}
actually shows that joint eigenfunctions $\psi(h)$ corresponding to
joint eigenvalues $E(h)\in K$ are in one-to-one
correspondence with microlocal solutions of the system
\[ (P_j(h)-E_j(h))u_h\sim 0 \] in a neighborhood of $\bigcup_{E\in
  K}\Lambda_E$. In other words, these eigenfunctions are determined by
the quantization conditions of theorem \ref{theo:global}.

\subsection{Quantum monodromy}
Away from the critical value of the momentum map $p$, the joint
spectrum locally looks like a lattice $h\ZM^2$. However, simply
patching these lattices together around $0$ is not possible. This
phenomenon is called \emph{quantum monodromy} and is described in
details in \cite{san-mono}. Let us take advantage of section \ref{sec:CBSreg}
to give a short description of it.

If $c_0$ is regular value of $p$, then on a neighborhood $B$ of $c_0$
there exists a map $f(h):B\fleche \RM^2$ that ``asymptotically
identifies $\Sigma(h)\cap B$ with the straight lattice $h\ZM^2$'', in
the following sense~: $f(h)$ is an elliptic semi-classical symbol of
order zero such that for any family $E(h)=(E_1(h),E_2(h))\in B$, $h\in
H$, 
\begin{equation}
  \label{equ:affinechart1}
         E(h)\in\Sigma(h)\cap B + O(h^\infty)\textrm{ if and only if } 
        f(E(h);h)\in h\ZM^2 + O(h^\infty).
\end{equation}
 Indeed, for $c\in B$ let
$\lambda_c(h)$ be the semi-classical holonomy 1-form on $\Lambda_c$
defined in proposition \ref{prop:holonomy}; then because of theorem
\ref{theo:quantization}, one can choose
\begin{equation}
  \label{equ:affinechart2}
  f(c;h) \egdef 
\left(\frac{1}{2\pi}\int_{\gamma_1(c)}\lambda_c(h),
\frac{1}{2\pi}\int_{\gamma_2(c)}\lambda_c(h)\right),
\end{equation}
where $(\gamma_1(c),\gamma_2(c))$ is any smooth basis of
$H_1(\Lambda_c)$.

It turns out that up to conjugation of $f(h)/h$ 
by an element of the affine group
with integer coefficients $GA(2,\ZM)$ such a map $f(h)$ is unique
modulo $O(h^\infty)$. If $\delta$ is a simple loop through
$c_0$ enclosing the critical value $0$ and small enough not to 
enclose any other critical value of the momentum map,
 we can cover its image by a
finite number of open ball $B_0,\dots,B_\ell$. Assuming that $B_0=B$,
we can construct a corresponding map $f_1(h)$ on $B_1$ that coincide
with $f(h)$ on $B_0\cap B_1$. Continuing this way, we end up with a
map $f_\ell(h)$ on $B_\ell$, and an element $\mu(\delta)\in GA(2,\ZM)$
such that 
\[ \frac{1}{h}f_\ell(h) = \mu(\delta)\circ (\frac{1}{h} f(h)) + O(h^\infty). \]
$\mu(\delta)$ depends only on the homotopy class of $\delta$ and is
called the \emph{quantum monodromy} of the joint spectrum. Because the
$f_j(h)$'s can be chosen as in equation (\ref{equ:affinechart2}), this
$\mu(\delta)$ is actually dual to the classical monodromy of the
bundle $\pi_1(\Lambda_c)\fleche c$~:
\begin{prop}[\cite{san-mono}]
  \label{prop:QM}
  If $(\gamma_1(c),\gamma_2(c))$ and $\epsilon$ 
  are chosen as in proposition \ref{prop:monodromie}
  then
  \[ \mu(\delta) = \iota \left(
    \begin{array}{cc}
      1 & -\epsilon \\
      0 & 1
    \end{array}\right), \]
where $\iota$ is the inclusion $GL(2,\RM)\hookrightarrow
GA(2,\RM)$ such that for any $M\in GL(2,\RM)$, the origin $0\in\RM^2$
is fixed by $\iota(M)$.
\end{prop}

The article \cite{san-mono} shows how this quantum monodromy can be told
from a picture of the spectrum by detecting the change in the lattice
structure of the joint spectrum around the critical value. One of the 
ways of doing this is to \emph{unwind} the spectrum onto $\RM^2$, 
in the following way. Suppose $\delta$ is actually drawn on the universal 
cover of the punctured plane~:
$\tilde{W}\flechedroite{\pi}\ \RM^2\setminus\{0\}$ 
and let $\tilde{B}_j$ be the open cover of its image such that 
$\pi(\tilde{B}_j)=B_j$. Then the maps $f_j(h)$ on $B_j$ define a 
global map $\hat{f}(h)$ on $\tilde{B}_0\cup \tilde{B}_1\dots\cup \tilde{B}_\ell$. 
By definition we 
have
\[  \hat{f}(\delta(1);h)/h = \mu(\delta)\circ \hat{f}(\delta(0);h)/h. \]
Since the restriction of $\hat{f}(h)$ to the points belonging to 
the spectrum $\pi^{-1}(\Sigma(h))$ takes values in $h\ZM^2 + 
O(h^\infty)$, there is for 
small $h$ a unique map 
\[ \tilde{f}(h) : \pi^{-1}(\Sigma(h))\cap\bigcup_j \tilde{B}_j \fleche
h\ZM^2\]
that coincides with $\hat{f}(h)$ modulo $O(h^\infty)$.  This map is of
course as well defined on the intersection with $\pi^{-1}(\Sigma(h))$
of any domain $\tilde{\D}$ of $\tilde{W}$ that projects into a small
(but fixed) annulus $\D$ around the critical value $0$, and is called
the \emph{unwinding} (or \emph{developing}) map of the spectrum.

For instance, if one uses this map to unwind a sequence of points
$\tilde{A}_j(h)\in\pi^{-1}(\Sigma(h))\cap \tilde{B}_j$ with
$\pi(\tilde{A}_\ell(h))=\pi(\tilde{A}_0(h))$ (see figure
\ref{fig:unwinding}), then one obtains an integral open polygonal line
$(A_0,\dots,A_\ell)$ in $h\ZM^2$, whose extremities satisfy
\[ A_\ell/h =  \mu(\delta)(A_0/h).\]
\begin{figure}[hbtp]
  \begin{center}
    \leavevmode
    \input{unwinding3.pstex_t}
    % begin change
    \macaption{unwinding of a sequence of points in $\Sigma(h)\cap\D$,
      in the case of the spectrum of the Champagne bottle (numerics by
      M.Child). In the case $b)$, the starting point belongs to
      $L_0(h)$, which implies that the unwound line is closed.}
    % end change
    \label{fig:unwinding}
  \end{center}
\end{figure}
For short, we shall allow us to refer to such a sequence of points
$\tilde{A}_j$ as a \emph{closed polygonal line in $\Sigma(h)$}
whenever these points are close enough to each other so that the line
segments $[A_j,A_{j+1}]$ lie in $\hat{f}(\tilde{\D};h)$. The
``segment'' between two consecutive points $\tilde{A}_j(h)$ and
$\tilde{A}_{j+1}(h)$ is by definition the preimage in $\tilde{\D}$ of
$[A_j,A_{j+1}]$.

\subsection{The exact counting function}
The purpose of this paragraph is to count the number of joint
eigenvalues inside a simple closed polygonal line around the critical
value $0$.  Let $L_0(h)$ be the subset of $\Sigma(h)\cap {\D}$ of
joint eigenvalues $E(h)$ that are unwound on the horizontal line
$\{y=0\}\subset\RM^2$ (in other words,
$\tilde{f}_2(\pi^{-1}(E(h));h)=0$). Because of proposition
\ref{prop:QM}, $L_0(h)$ is pointwise fixed by the quantum monodromy
$\mu(\delta)$. That means that any closed polygonal line
$\tilde{\A}(h)=(\tilde{A}_0(h),\dots,\tilde{A}_\ell(h))$ on
$\Sigma(h)$ starting on $L_0(h)$ is unwound as a \emph{closed}
polygonal line $\A(h)=(A_0(h),\dots,A_\ell(h))$ in $\RM^2$ (see figure
\ref{fig:unwinding}).

It is easy to see that if we have another such polygonal line on
$\Sigma(h)\cap {\D}$ that lies \emph{inside} the previous one, then
the number of joint eigenvalue \emph{between} these two is exactly the
number of $h$-integral lattice points between the two unwound polygons
(see lemma \ref{lemm:simplyconnected}). As a consequence, there should
be a universal index $\nu\in\ZM$ such that the number
$N_h(\tilde{\A}(h))$ of joint eigenvalues inside $\tilde{\A}(h)$ be
\[  N_h(\tilde{\A}(h)) = N_h(\A(h)) + \nu,\]
where $N_h(\A(h))$ denotes the number of $h$-integral points inside 
the polygon $A(h)$. Note that the latter is given by Pick's formula
\cite{pick}~:\\
\begin{equation}
  \label{equ:pick}
  N_h(\A(h)) = \textrm{Area}(\A(h)) +
\frac{\ZM\textrm{Length}(\A(h))}{2} + 1,
\end{equation}
where $\ZM\textrm{Length}(\A(h))$ is the number of $h$-integral points
on the boundary.

The result is actually that $\nu=0$~:
\begin{theo}
  \label{theo:compter}
  Let $\tilde{\A}(h)=(\tilde{A}_0(h),\dots,\tilde{A}_\ell(h))$ be a 
  simple closed polygonal line on $\Sigma(h)$ starting on $L_0(h)$ 
  and enclosing one the origin, and let 
  $\A(h)=(A_0(h), \dots, A_\ell(h)=A_0(h))$ be the unwound polygonal 
  line in $\RM^2$.  Denote by $N_h(\tilde{\A}(h))$ the number of 
  joint eigenvalues inside $\tilde{\A}(h)$ (counting the points in 
  the boundary), and similarly denote by $N_h(\A(h))$ the number of 
  $h$-integral points inside the polygon $A(h)$, counting the 
  boundary.
  
  Then
  \[  N_h(\tilde{\A}(h)) = N_h(\A(h)).  \]
  In other words, the recipe to calculate $N_h(\tilde{\A}(h))$ is
  simply to unwind $\tilde{\A}(h)$ and then apply Pick's formula
  (\ref{equ:pick}).
\end{theo}
\demo The idea is to cut the ``polytope'' delimited by $\tilde{\A}(h)$
into two halves separated by $L_0$; then on each part the control over
the joint eigenvalues given by theorem \ref{theo:global} is used to
shift the problem into a region where the result will be obvious, due
to the following fact~:
\begin{lemm}
  \label{lemm:simplyconnected}
  if $B$ is a simply connected open subset of $\D$, then for any closed 
  polygonal line $\tilde{\C}(h)$ in $\Sigma(h)\cap B$, the unwound
  line $\C(h)$ is also closed and we have
  \[ N_h(\tilde{\C}(h)) = N_h(\C(h)).   \]
\end{lemm}
This lemma also proves the existence of the index $\nu$ claimed before
the statement of the theorem, and this implies that we can choose any
particular $\tilde{\A}(h)$ to prove the theorem. So let us assume that
the intersection of the polytope delimited by $\tilde{\A}(h)$ with the
horizontal lines $\epsilon_2=hn$ are segments whose extremities are
vertices of $\tilde{\A}$.% (figure \ref{fig:segments}).

Now, let us apply theorem \ref{theo:global} to the operators
$P_j^E(h)\egdef P_j(h)-E_j$, and keep the notations from the statement
of this theorem. We know from lemma \ref{lem:FNQ} that the map
\[ E=(E_1,E_2) \fleche 
\bep^E(h)=(\epsilon_1^E(h),\epsilon_2^E(h)) \] is an elliptic
semi-classical symbol of order zero. It is therefore for small $h$ a
diffeomorphism of a compact neighborhood $K$ of zero into its image;
this implies that the determination of the joint spectrum in $K$
amounts to that of the corresponding values of $\bep^E(h)$, which we
shall hence simply denote by $\bep=(\epsilon_1,\epsilon_2)$. From now
on, we shall assume that the annulus $\D$ defined in the previous
paragraph is a subset of $K$.

Let us look for solutions satisfying $\epsilon_2\geq 0$. Let $g_1(h)$
and $g_2(h)$ be the following functions~:
\[ \frac{2\pi}{h}g_1(\epsilon;h) = 
(\lambda^E)^{\out} +
\epsilon_2\frac{\pi}{2} -
\frac{\epsilon_1}{h}\ln(2h) - 2\arg
\Gamma\left(\frac{i\epsilon_1/h+1+\epsilon_2}{2}\right), \]
\[ \textrm{and }g_2(\epsilon;h) = \epsilon_2.\]
Because of the analyticity properties of $\Gamma$, these functions are
holomorphic in $\epsilon_2\geq 0$.  The map $f(E;h)\egdef
g(\epsilon^E(h);h)$ satisfies equation (\ref{equ:affinechart1}) thanks
to theorem \ref{theo:global}. Because of the $\ln h$ term in it, it is
\emph{not} a standard semi-classical symbol globally; but we know from
paragraph \ref{sec:recover} that it is so for $E$ in a small compact
at finite distance from the origin. Moreover, because of remark
\ref{rem:recover}, it must satisfy equation (\ref{equ:affinechart2})
as well, where $(\gamma_1(c),\gamma_2(c))$ is the basis given in
proposition \ref{prop:monodromie}.

Therefore, the joint spectrum is spread on the horizontal lines
$\epsilon_2=hn$, $n\in\ZM$, and the set $L_0(h)$ lies in the line
$\epsilon_2=0$. Our problem is now to count the number of joint
eigenvalues inside the intersection $\tilde{\A}^+(h)$ of the half
plane $\epsilon_2\geq 0$ with the polytope delimited by
$\tilde{\A}(h)$ (figure \ref{fig:segments}).
%\begin{figure}[hbtp]
%  \begin{center}
%    \leavevmode
%    \input{halfplane.pstex_t}
%    \macaption{the upper half of the polytope.}
%    \label{fig:halfplane}
%  \end{center}
%\end{figure}
 Let us denote by $\tilde{I}_n(h)$ the set of joint eigenvalues
 contained in $\tilde{\A}^+(h)$ and satisfying $\epsilon_2=hn$.
\begin{figure}[hbtp]
  \begin{center}
    \leavevmode
    \input{segments.pstex_t}
    \macaption{the upper half of the polytope.}
    \label{fig:segments}
  \end{center}
\end{figure}

From now on, we forget the negative half of our polytope, and label
\[\tilde{A}_0,\dots,\tilde{A}_\ell\]
 the consecutive vertices of
$\tilde{\A}^+(h)$, starting and ending on $L_0(h)$. We denote by
$N_h(\tilde{\A}(h))$ the number of joint eigenvalues inside
$\tilde{\A}^+(h)$, and by $N_h(\A(h))$ the number of $h$-integral
points inside the polytope delimited by the unwound points
$A_0,\dots,A_\ell,A_0$. The theorem will be proved provided we show
that $N_h(\tilde{\A}(h))=N_h(\A(h))$.

The next step of the proof is now to translate $\tilde{\A}^+(h)$ far
from the origin, as follows (figure \ref{fig:translate}).  The domain
$\D$ where the unwinding map is defined can always be chosen such that
there exists an integer $k(h)\in\ZM$ such that the origin is outside
of the polytope
$\tilde{\C}(h)=(\tilde{C}_0(h),\dots,\tilde{C}_\ell(h))$ defined by
\[  \tilde{f}(\tilde{C}_j(h);h) = \tilde{f}(\tilde{A}_j(h);h) + 
h(k(h),0).\]
\begin{figure}[hbtp]
  \begin{center}
    \leavevmode
    \input{translate.pstex_t}
    \macaption{translation of the polytope}
    \label{fig:translate}
  \end{center}
\end{figure}

Of course, the set of joint eigenvalues inside 
$\tilde{\C}(h)$ satisfying  $\epsilon_2=hn$ is in bijection with 
$\tilde{I}_n(h)$, which implies that 
\[ N_h(\tilde{\C}(h)) = N_h(\tilde{\A}(h)). \]

On the other hand, because $f(\cdot;h)$ is a valid unwinding function
in a region containing both $\tilde{\A}(h)$ and $\tilde{\C}(h)$, the
unwinding $\C(h)$ of $\tilde{\C}(h)$ is by definition the translation
of $\A(h)$ by the horizontal vector $h(k(h),0)$.  It has therefore the
same number of interior points~:
\[  N_h(\C(h)) = N_h(\A(h)). \]

But lemma \ref{lemm:simplyconnected} yields
\[ N_h(\tilde{\C}(h)) = N_h(\C(h)),\]
which finally proves the theorem. \cqfd

\subsection{Shape of the spectrum near the critical value}
We investigate here the shape of the joint spectrum in a region of
size $O(h)$ near the origin, in terms of the $(\epsilon_1,\epsilon_2)$
variables. We already know that the joint eigenvalues are distributed
on the horizontal lines $\epsilon_2=hn$, $n\in\ZM$. It remains to
study their repartition on each of these lines. For this purpose, we
fix $n\in\ZM$ and denote by $x=\epsilon_1/h$ the generic variable on
$\{\epsilon_2=hn\}$. Of course, the results will be based upon
corollary \ref{coro:global}, whose notations we use here, and which
we express in the following way~:
\begin{prop}[corollary \ref{coro:global}]
  \label{prop:zoom}
  Let $K$ a compact neighborhood of the origin in $\RM^2$. For any
  family $E(h)=(E_1(h),E_2(h))\in(hK)$ such that
  $\epsilon_2^{E(h)}(h)=hn + \oh$ ($n\in\ZM$), the following
  characterization holds~:
  \[ E(h)\in \Sigma(h) +O(h^\infty) \quad \textrm{if and only if} \]
  \[  |n|\frac{\pi}{2} -  x(h)\ln(2h) -
  2\arg\Gamma\left(\frac{ix(h)+1+|n|}{2}\right) +
  f_n(x(h);h) \in 2\pi\ZM + O(h^\infty),\]
  with $x(h)=\epsilon^{E(h)}_1(h)/h$. Here $f_n(h)$ is a symbol of the
  form~:
  \[ f_n(x;h) = \frac{A}{h} + xB + nC + D + \mu\frac{\pi}{2} + hf_{n,1}(x) +
  h^2f_{n,2}(x) + \cdots, \]
  for some real constants $A,B,C,D$.
\end{prop}
\demo Of course, in view of corollary \ref{coro:global} -- applied to
the system $\tilde{P}^t_j=P_j-ht_j$, $j=1,2$ --, $A$ is just the
action integral on the critical Lagrangean~: $A=\int_{\gamma_1^0}
\alpha_0$, and it only remains to prove that there exists constants
$B$, $C$, $D$ such that
\[ I_{\gamma_1^0}(\tilde{\kappa}^t) = x(h)B + nC + D + O(h), \]
where $t$ and $x(h)$ are linked by 
\[ x(h)=\tilde{\epsilon_1}^t(h)=\frac{1}{h}\epsilon_1^{ht}(h).\]
Because the sub-principal symbols of $\tilde{P}^t_j$ are affine
functions of $t$~: $\tilde{r}^t=r-t$, it is easy to deduce from
formula (\ref{equ:regularisation}) of proposition
\ref{prop:regularisation} that $I_{\gamma_1^0}(\tilde{\kappa}^t)$ is
affine in $t$. The result then follows from the fact that the
principal part of $\tilde{\bep}^t(h)=(x(h),n +O(h^\infty))$ is also an
-- invertible -- affine function of $t$. \cqfd

We study the distribution of joint eigenvalues on the line
$\epsilon_2=hn$ by investigating the \emph{gaps} between two
consecutive eigenvalues.
\begin{theo}
  \label{theo:gaps}
  For any $n\in\ZM$, let $E_k(h)$ be the sequence of joint eigenvalues
  in $\{\epsilon_2=hn\}\cap (hK)$, such that the corresponding
  variables $x_k(h)$ form (for each $h$) a strictly increasing finite
  sequence.

  The gap between $E_k(h)$ and $E_{k+1}(h)$ is given by~:
  \begin{equation}
    \label{equ:gaps}
    \left\|\frac{E_{k+1}(h)-E_k(h)}{h}\right\| =  
    \frac{2\pi a}{|ln h| + B - \ln 2 -
      \Psi_n'(x)} (1+O(h))
  \end{equation}
  \[ = a\left(\frac{2\pi}{|\ln h|} +
    \frac{2\pi}{(\ln h)^2}\left(\Psi_n'(x)-B+\ln 2\right)\right) +
  O\left(\frac{1}{(\ln h)^3}\right),\] for an $x\in ]x_k,x_{k+1}[$,
  and $a=\|M\cdot(1,0)\|$; we have denoted
  \[ \Psi_n(x) = 2\arg\Gamma\left(\frac{ix+1+|n|}{2}\right).\]
\end{theo}
Note that $\Psi_n'$ is related to the logarithmic derivative $\psi$ of
the Gamma-function~:
\[ \Psi_n'(x) = \re\left(\psi\left(\frac{ix+1+n}{2}\right)\right). \]
A consequence of this theorem is that the number of $E_k(h)$'s inside
$hK$ increases as $h$ tends to zero at a rate equivalent to
$const.|\ln h|$, which easily yields a Weyl type asymptotic formula
for the counting function (see theorem \ref{theo:weyl} in next
paragraph). A more precise result is that the function $\Psi_n'$
controls the shape of the graph ``gaps versus $x$''. It turns out that
$\Psi_n'$ has a non-degenerate minimum at $x=0$, which is particularly
sharp for $n=0$ (figure \ref{fig:cusp-gamma}).
\begin{figure}[hbtp]
  \begin{center}
    \leavevmode
    \input{cusp-gamma.pstex_t}
    \macaption{the function $\Psi_0'(x)$ (here $\gamma$ denotes Euler's
    constant $\gamma \simeq 0,5772156649\dots$).}
    \label{fig:cusp-gamma}
  \end{center}
\end{figure}
This might have some interesting applications in numerical
analysis. First, it should refrain people from thinking that the gap
curve has a cusp; secondly, however, this ``near-cusp'' might be
useful to precisely locate the critical value $x=0$. Numerics for the
example of the Champagne bottle suggest that both remarks apply for
reasonable value of the small constant $h$ (figures \ref{fig:cusp} and
\ref{fig:cusp-z}).
\begin{figure}[hbtp]
  \begin{center}
    \leavevmode \input{cusp.pstex_t}
    \macaption{spectral gaps for the quantum Champagne bottle, $n=0$,
      $h=10^{-4}$ (numerics).}
    \label{fig:cusp}
  \end{center}
\end{figure}
\begin{figure}[hbtp]
  \begin{center}
    \leavevmode
    \input{cusp-z.pstex_t}
    \macaption{spectral gaps for the quantum Champagne bottle, $n=0$,
      $h=10^{-5}$ (numerics).}
    \label{fig:cusp-z}
  \end{center}
\end{figure}
Finally, of course, this theorem can be used to find a rigorous
approximation for the spectral gaps. In the example of the Champagne
bottle, formula (\ref{equ:gaps}) can be made completely explicit. The
Hamiltonian $H$ is the following~:
\[ H(x,y,\xi,\eta) = \frac{1}{2}(\xi^2+\eta^2) - r^2 + r^4, \]
with $r^2=x^2+y^2$. It commutes with the angular momentum
$I=x\eta-y\xi$, hence we can let $p=(p_1,p_2)=(H,I)$. The
corresponding $h$-operators are
\[ \hat{H}=-\frac{h^2}{2}\Delta -r^2+r^4, \quad
\hat{I}=\frac{h}{i}\deriv{}{\theta}.\]
One can compute the quantities involved by formula (\ref{equ:gaps}), and one finds 
\[ M = \left(
  \begin{array}{cc}
    \sqrt{2} & 0\\ 0 & 1
  \end{array}\right) \textrm{ so } a=\sqrt{2}, \quad \textrm{ and }
B=(5/2)\ln 2. \] 
Thus, using formula (\ref{equ:gaps}) to approximate the
\emph{smallest} gap (obtained at $x=0$ for $n=0$) we get~:
\[  \min_k\left(\frac{E_{k+1}-E_k}{h}\right) \sim
\frac{2\pi\sqrt{2}}{|\ln h| + (9/2)\ln 2 + \gamma} + O(h). \]
The comparison between the principal part of this formula as a
function of $\ln h$ and the ``experimental'' result from the numerical
computation of the spectrum is depicted in figure
\ref{fig:gaps-formule}, and tends to show that the approximation is
quite accurate.
\begin{figure}[hbtp]
  \begin{center}
    \leavevmode
    \input{gaps-formule.pstex_t}
    \macaption{The smallest gap; comparison between numerics and the
    top order term in formula (\ref{equ:gaps}).}
    \label{fig:gaps-formule}
  \end{center}
\end{figure}

\demo[of theorem \ref{theo:gaps}] Because $(E_{k+1}-E_k)/h =
(M+O(h))\cdot(x_{k+1}-x_k,0) + \oh$, it suffices to estimate
$|x_{k+1}-x_k|$.

For $n\in\ZM$, let
\[ g_n(x;h) \egdef \frac{1}{2\pi}\left( |n|\frac{\pi}{2} -  x\ln(2h) -
  2\arg\Gamma\left(\frac{ix+1+|n|}{2}\right) +
  f_n(x;h) \right).\]
The function $x\fleche g_n(x;h)$ is smooth and we have
\[ g'_n(x;h) = B - \ln(2h) - \Psi_n'(x) + O(h).\]
For small $h$, $g_n(h)$ is therefore strictly increasing, with slope
of order $|\ln h|$. The solutions $x_k(h)$ to $g(x_k(h);h)=k\in\ZM +
O(h^\infty)$ are therefore, for small $h$, in one-to-one
correspondence with the exact solutions to $g(x_k(h);h)=k$ (except
maybe for the first and the last ones, sitting near the boundary of
the compact interval in which we are looking for
solutions). Moreover the differences between the exact and the
$O(h^\infty)$ solutions are of order $\frac{1}{|\ln h|}O(h^\infty)$,
which is \emph{a fortiori} $O(h^\infty)$. In the following $x_k$
denotes the exact solution.

By Rolle's theorem, there is an $x\in]x_k,x_{k+1}[$ such that
\[ \frac{1}{|x_{k+1}-x_k|} = g'_n(x;h). \]
This means
\[ |x_{k+1}-x_k| = \frac{2\pi}{|ln h|+B-\ln 2 - \Psi_n'(x) + O(h)}, \]
which proves the theorem. \cqfd

\subsection{Weyl's formula}
Let $K$ be a compact convex subset of $\RM^2$, and let $\tilde{K}(h)$
be its image by the change of variables:
\[ (t_1,t_2)\fleche
(\tilde{\epsilon}^t_1(h),\tilde{\epsilon}^t_2(h)). \]
Recall that the principal part of this map is
\[ \tilde{\bep}_0^t= M^{-1}\cdot(t-r(m)),\]
where $r(m)=(r_1(m),r_2(m))$ is the value at $m$ of the joint
sub-principal symbol of the system $P_1(h),P_2(h)$, and $M=dF(0)$,
where $(p_1,p_2)=F(q_1,q_2)$ for some symplectic coordinates near $m$
(see section \ref{sec:microloc}).

Denote by $\tilde{K}_0$ the principal part of $\tilde{K}(h)$; it is a
convex compact of $\RM^2$. Theorem \ref{theo:gaps} easily yields the
following estimate for the counting function~:
% containing the critical value $0$ in its interior.  
\begin{theo}
  \label{theo:weyl}
  Denote by $N_h(K)$
  the number of eigenvalues belonging to 
  $hK$~: 
  \[ N_h(K)\egdef \#\{(E_1,E_2)\in hK\cap \Sigma(h)\}.\]
  The following estimate holds, as $h$ tends to zero~:
  \[ N_h(K) = \frac{|\ln h|}{2\pi} \int_{\tilde{K}_0} |dx\wedge dn| +
  O(1).\]
  The measure  $|dn|$ is the counting measure on $\ZM$, and $|dx|$ is
  the standard one-dimensional Lebesgue measure.
\end{theo}
We present in figure \ref{fig:weyl} a numerical illustration of this
result.
\begin{figure}[hbtp]
  \begin{center}
    \leavevmode
    \input{log-law.pstex_t}
    \macaption{The counting function $N_h(K)$ for the Champagne
      bottle, where $K$ is a compact selecting only joint eigenvalues
      on the line $\epsilon_2=0$ (numerics).}
    \label{fig:weyl}
  \end{center}
\end{figure}

The question that this theorem arouses is whether the formula for
$N_h$ can be related to a Weyl-type heuristic involving volume in
phase space. The answer is positive (proposition \ref{prop:DH}), up to
the difference between the measures of $\tilde{K}_0$ with respect to
$|dx\wedge dn|$ and with respect to the standard Lebesgue measure on
$\RM^2$, which of course strongly depends on the geometry of
$\tilde{K}_0$.

Let $\mu$ be the push-forward on $\RM^2$ of the symplectic (or
Liouville) measure by the momentum map $p=(p_1,p_2)$~: for any Borel
$\mathcal{B}$ of $\RM^2$,
\[ \mu(\mathcal{B})=\int_{p^{-1}(\mathcal{B})}|\omega\wedge\omega|/2. \]
Let $|dc|=|dc_1\wedge dc_2|$ be the standard Lebesgue measure on 
$\RM^2$. The Lebesgue volume of $\mathcal{B}$ will be denoted by 
$|\mathcal{B}|$.
\begin{prop}
  \label{prop:DH}
  As $h$ tends to zero, 
  \[ \frac{\mu(hK)}{(2\pi h)^2}  = \frac{|\ln h|}{2\pi}
  |\tilde{K}_0| + O(1). \]
\end{prop}
\demo Let $\Omega$ be a neighborhood of the critical point $m$ in
which the system $(P_1,P_2)$ can be brought to the normal form of
lemma \ref{lem:FNQ}, so that there is a local diffeomorphism $F$
defined on $U=p(\Omega)$ such that, with the same notations as in
section \ref{sec:linear}, $(p_1,p_2)=F(q_1,q_2)$.

Replacing $K$ by $h_0K$ with $h_0$ small enough, if necessary, we can 
always assume that $K\subset U$. 
%Denote now $\tilde{K}=F^{-1}(K)$, 
Let $\varepsilon>0$ be such that the ball 
\[ B_\varepsilon\egdef  \{(r,\theta,\rho,\alpha)\in M, \quad 
r<\varepsilon, \rho<\varepsilon\}\]
is included in $\Omega$.
Then $\mu(hK)$ can be split into two parts~:
\begin{equation}
  \label{equ:splitmu}
  \mu(hK) = \int_{B_\varepsilon\cap\{(q_1,q_2)\in F^{-1}(hK)\}}
  |dxdyd\xi d\eta| + 
  \int_{B_\varepsilon^{\mathsf{c}}\cap\{(p_1,p_2)\in hK\}}
  |\omega\wedge\omega|/2.
\end{equation}

If $h$ is small enough, each fiber $p^{-1}(c)$, $c\in hK$ intersects 
$B_\varepsilon$, and we know from section \ref{sec:fibration} that the 
fibration $p_{\restr B_\varepsilon^{\mathsf{c}}}$ over $hK$ is trivial with 
fibers diffeomorphic to compact cylinders.  Moreover, $p$ naturally 
defines on each cylinder a locally invariant measure $\brho$ such that 
the symplectic measure splits according to~:
\[ |\omega\wedge\omega|/2 = \brho\otimes|dp_1\wedge dp_2|. \]
Therefore, the second term of (\ref{equ:splitmu}) can be written as~:
\[ \int_{B_\varepsilon^{\mathsf{c}}\cap\{(p_1,p_2)\in hK\}}
|\omega\wedge\omega|/2 = \int_{(c_1,c_2)\in hK}
\brho(B_\varepsilon^{\mathsf{c}}\cap p^{-1}(c)) |dc_1\wedge dc_2|.\] 
Because $\brho(B_\varepsilon^{\mathsf{c}}\cap p^{-1}(c))$ is bounded for
bounded $(c_1,c_2)$, the following estimate holds, for some $C>0$~:
\[  \int_{B_\varepsilon^{\mathsf{c}}\cap\{(p_1,p_2)\in hK\}}
|\omega\wedge\omega|/2 < C.\int_{(c_1,c_2)\in hK} 
|dc_1\wedge dc_2| = C.h^2|K|.\]
        
Let us turn now to the first term in the sum (\ref{equ:splitmu})~:
\[ I_\varepsilon(h) \egdef  
\int_{B_\varepsilon\cap\{(q_1,q_2)\in F^{-1}(hK)\}}
        |dxdyd\xi d\eta|.\]
Using the new coordinates $(r,\theta,\ell,\lambda)$ with $\ell=r\rho$ 
and $\lambda=\alpha-\theta$, the symplectic measure is transformed 
into \[ \frac{\ell}{r}|d\lambda\wedge d\ell \wedge d\theta \wedge dr|
= \frac{1}{r}|dq_1\wedge dq_2|\otimes |d\theta \wedge dr|. \]
This gives~:
\[ I_\varepsilon(h) = 2\pi\int_{(c_1,c_2)\in F^{-1}(hK)}
1_{|c|\leq\varepsilon^2}
\left(\int_{\frac{|c|}{\varepsilon}}^\varepsilon  \frac{1}{r} |dr|
\right) |dc_1\wedge dc_2|\]
\[ = 2\pi\int_{F^{-1}(hK)}1_{|c|\leq\varepsilon^2} 
\ln(\frac{\varepsilon^2}{|c|}) |dc_1\wedge dc_2|.\]
Note that for $h$ small enough, the condition $|c|\leq\varepsilon^2$ 
is always satisfied for $c\in F^{-1}(hK)$.

Now, the change of variables $\hat{c}=F(c)/h$ yields
\[ I_\varepsilon(h) = 2\pi h^2\int_{\hat{c}\in K} 
\ln(\frac{\varepsilon^2}{|F^{-1}(h\hat{c})|}) J(h\hat{c})|d\hat{c}|, \]
where $J$ is the Jacobian of $F^{-1}$.
Taylor's formula for $J$ and $F^{-1}$ easily shows that the integrand 
can be reduced to its principal part~:
\[ I_\varepsilon(h) = 2\pi h^2\int_{\hat{c}\in K} 
\ln(\frac{\varepsilon^2}{|M^{-1}.h\hat{c}|}) J(0)|d\hat{c}| + O(h^3). \]
Since $J(0)=\det M^{-1}$, the right hand-side is equal to~:
\[ 2\pi h^2\int_{\tilde{c}\in M^{-1}(K)} 
\ln(\frac{\varepsilon^2}{|h\tilde{c}|}) |d\tilde{c}| + O(h^3), \]
with $\tilde{c}=M^{-1}.\hat{c}$.
Because the measure $\ln|c||dc|$ is locally integrable at the origin, 
this finally gives~:
\[ I_\varepsilon(h) = 
-2\pi h^2\ln h|M^{-1}(K)| + O(h^2),\]
thus proving proposition \ref{prop:DH}. \cqfd

\bibliographystyle{plain}
\bibliography{bibli}
\end{document}